%% file: thesis_with_title.tex
\documentclass[12pt]{amsbook}
\usepackage{graphicx}
\usepackage{amsmath}
\usepackage{amssymb}
\usepackage{amscd}
\usepackage{amsthm}
\usepackage{latexsym}

\makeindex

\numberwithin{section}{chapter}
\numberwithin{equation}{chapter}
\newtheorem{theorem}{Theorem}[chapter]
\newtheorem{lemma}[theorem]{Lemma}
\newtheorem{corollary}[theorem]{Corollary}
\newtheorem{proposition}[theorem]{Proposition}

\theoremstyle{definition}
\newtheorem{definition}[theorem]{Definition}

\newtheorem{postulate}{Postulate}

\theoremstyle{remark}
\newtheorem{remark}[theorem]{Remark}

\newtheorem{summary}[theorem]{Summary}

\newcommand{\set}[1]{\left\{#1\right\}}
\newcommand{\wt}[1]{\widetilde{#1}}
\newcommand{\Real}{\mathbb R}

\newcommand{\Si}{\Sigma}
\newcommand{\To}{\longrightarrow}
\newcommand{\x}{\mathbf{x}}

\newcommand{\n}{\mathbf{n}}
\newcommand{\N}{\mathbf{N}}
\newcommand{\K}{\mathcal{K}}

\newcommand{\G}{\Gamma}
\newcommand{\tans}{T_{p}\Sigma}
\newcommand{\tanm}{T_p M}
\title{General Relativistic Shock-Waves\\
Propagating at the Speed of Light}
\author{Michael B. Scott}
\date{July 26, 2002}
\begin{document}
%\pagenumbering{roman}
\frontmatter
% ---Title Page---------------------------------------------------
\begin{center}
{\large \bf General Relativistic Shock-Waves\\
Propagating at the Speed of Light}

\vspace{12pt}
\normalsize
By

MICHAEL BRIAN SCOTT

B.S. (California State University, Northridge) 1992

M.S. (California State University, Northridge) 1996

\vspace{12pt}

DISSERTATION

Submitted in partial satisfaction of the requirements for the degree of

DOCTOR OF PHILOSOPHY

in

MATHEMATICS

in the

OFFICE OF GRADUATE STUDIES

of the

UNIVERSITY OF CALIFORNIA

Davis

Approved:

\vspace{10pt}

\rule{3.0in}{.01in} \vspace{10pt}

\rule{3.0in}{.01in} \vspace{10pt}

\rule{3.0in}{.01in} \vspace{6pt}

Committee in Charge

\vspace{0.12in}

2002

\end{center}
% ----------------------------------------------------------------
%
% ----------------------------------------------------------------
%\newpage
%{\Large \bf ACKNOWLEDGEMENTS}
\chapter*{Acknowledgements}
First and foremost I would like to thank my wife, Monika. Without
her support and encouragement it is difficult to imagine
completing this dissertation. I also owe a great deal to my thesis
advisor Professor Blake Temple. Besides giving his support and
encouragement, I have also benefitted greatly from his influence
on my approach to research mathematics, and I look forward to
working with him in the future.

It has been a long road for me to get to this point in my career,
and I am certain I would not have made it this far without the
patience, and support of my parents, Richard and Jean Scott. I am
also grateful to them for letting me pursue a career of my
choosing, and offering nothing but encouragement. My only wish is
that my mother were still with us to witness this accomplishment.

I would also like to thank my wife's parents Chris and Eleni
Poulos for their support throughout my time at U.C. Davis even
before my wife and I were married.

% ----------------------------------------------------------------

\chapter*{\underline{Abstract}}

%\normalsize
\noindent We investigate shock-wave solutions of the Einstein
equations in the case when the speed of propagation is equal to
the speed of light. The work extends the shock matching theory of
Smoller and Temple to the lightlike case. After a brief
introduction to general relativity, we introduce a previously
known generalization of the second fundamental form by
Barrab\`{e}s and Israel. Then we use this to develop an extension
of a shock matching theory, which characterizes solutions of the
Einstein equations when the spacetime metric is only Lipschitz
continuous across a hypersurface, to include the case when the
hypersurface is lightlike. The theory also demonstrates an
unexpected result that the matching of the generalized second
fundamental form alone is not a sufficient condition for
conservation conditions to hold across the interface. Using this
theory we then construct a new exact solution of the Einstein
equations that can be interpreted as an outgoing spherical shock
wave that propagates at the speed of light. This is done by
matching a Friedman Robertson Walker (FRW) metric, which is a
geometric model for the universe, to a Tolman Oppenheimer Volkoff
(TOV) metric, which models a static isothermal spacetime. Then our
theory is used to show that the matched FRW, TOV metric is a
solution. The pressure and density are finite on each side of the
shock throughout the solution, the sound speeds, on each side of
the shock, are constant and subluminous. Moreover, the pressure
and density are smaller at the leading edge of the shock which is
consistent with the Lax entropy condition in classical gas
dynamics. However, the shock speed is greater than all the
characteristic speeds. The solution also yields a surprising
result in that the solution is not equal to the limit of
previously known subluminous solutions as they tend to the speed
of light.

% ---Table of Contents--------------------------------------------
\tableofcontents
% ---Notation and Conventions-------------------------------------
%\newpage
%{\Large \bf Notation and Conventions}
\chapter*{Notation and Conventions}
\input{notation.tex}

% ---Begin Dissertation-------------------------------------------
\mainmatter
% ---Chapter 1----------------------------------------------------
\chapter{Introduction}
\label{introduction}
\input{intro.tex}

% ---Chapter 2----------------------------------------------------
\chapter{General Relativity and the Second Fundamental Form}
\label{background}
\input{back.tex}

% ---Chapter 3----------------------------------------------------
\chapter{Lightlike Shock-Wave Solutions of the Einstein Equations}
\label{main}
\input{main.tex}

% ---Chapter 4----------------------------------------------------
\chapter{An Exact Lightlike Shock-Wave Solution of the Einstein Equations}
\label{example}
\input{exam.tex}

% ---Chapter 5----------------------------------------------------
\chapter{Conclusions and Summary of Contributions}
\label{conclusion}
\input{conc.tex}

% ---Back Matter--------------------------------------------------
%\renewcommand{\baselinestretch}{1}

\backmatter

\nocite{weinberg, wald, schutz}

\bibliographystyle{amsplain}

\bibliography{master}

\printindex
\end{document}

%% file: notation.tex
% ---Notation-----------------------------------------------------

\noindent $R^{\mu}_{\hspace{4pt}\alpha\nu\beta}$ denotes the
Riemann curvature tensor.

\noindent $R_{\alpha\beta}=R^{\mu}_{\hspace{4pt}\alpha\mu\beta}$
denotes the Ricci curvature tensor.

\noindent $R=g^{\alpha\beta}R_{\alpha\beta}$ denotes the Ricci
scalar.

\noindent $G$, $G_{\alpha\beta}=
R_{\alpha\beta}-(1/2)g_{\alpha\beta}R$ denotes the Einstein
curvature tensor.

\noindent $M$ denotes the $n$- or four-dimensional spacetime
manifold.

\noindent $\Si$ denotes a hyper- or shock-surface in $n$- or
four-dimensional spacetime.

\noindent $\tans$ denotes the tangent space of $\Si$ at a point
$p$ in $\Si$.

\noindent $\tanm$ denotes the tangent space at $p$ in the
spacetime manifold $M$.

\noindent $X_a$ denotes a basis vector of $\tans$.

\noindent $X$ denotes a tangent vector in $\tans$.

\noindent $\n$ denotes the vector normal to $\tans$

\noindent $\N$ denotes a vector transverse to $\Si$, which cannot
be in $\tans$.

\noindent $[\,\cdot\,]$ denotes the jump in a quantity across
$\Si$. For example, $[g(p)]=g^L(p)-g^R(p)$ for $p$ in $\Si$ where
$g=g^L$ on the left side of $\Si$, and $g=g^R$ on the right side
of $\Si$.

\noindent $K(X)=-\nabla_X \n$ denotes the second fundamental form.

\noindent $\K(X)=-\nabla_X \N$ denotes a generalized second
fundamental form which depends $\N$.

\noindent $C^k$ denotes a function which is at least $k$ times
differentiable and it $k$th derivative is continuous.

\noindent $C^{k,1}$ is a $C^k$ function whose $k$th derivative is
Lipschitz continuous.

\noindent $\mbox{diag}(-1,1,\ldots,1)$ denotes an $n\times n$
diagonal matrix with the first entry equal to $-1$ and the
remaining non-zero entries equal to $1$.

\noindent Notation for derivatives: For partial derivatives
$\partial F/\partial x^a=F_{,a}$. For covariant derivatives
$\nabla_a F=F_{;a}.$

%% file: intro.tex
% ---Introduction-------------------------------------------------
\label{chp:1}

In this dissertation we give a general theory of shock matching in
the lightlike case, and use this theory to construct a new exact
solution of the Einstein equations that can be interpreted as an
outgoing spherical shock wave that propagates at the speed of
light. The general theory extends the shock matching theory of
Smoller and Temple~\cite{st94,st95} to the case of lightlike
interfaces. Based on this new theory we construct our exact
solution by matching an Friedman Robertson Walker (FRW) metric to
a Tolman Oppenheimer Volkoff (TOV) metric across an outgoing
radial, lightlike shock wave. In this exact solution matter
crosses the interface, but nothing propagates at the speed of
light except the shock wave. As far as we know this is the first
such exact solution in general relativity.

In this exact solution the shock wave emerges from the FRW origin
at the instant of the Big Bang, and propagates all the way out to
infinity. The pressure and density are finite on each side of the
shock throughout the solution, the sound speeds, on each side of
the shock, are constant and subluminous. Moreover, the pressure
and density are smaller at the leading edge of the shock which is
consistent with the entropy condition in classical gas
dynamics~\cite{lax}. However, the shock speed is greater than all
the characteristic speeds, see~\cite{st95}. Subluminous shocks
with this characteristic condition were constructed
in~\cite{st95}.

In the general theory we translate generalized notion of the
second fundamental form of Barrab\`{e}s and Israel~\cite{isr91}
into the shock matching framework of Smoller and
Temple~\cite{st94,st95}. We base the analysis on a modified
Gaussian Skew (MGS) coordinate system on lightlike surfaces which
we introduce in this dissertation.

One surprise is that the relation between the sound speeds on the
front and back sides of the shock in our new exact solution, are
not equal to the speeds in the limit as the shock speed in the
Smoller-Temple solutions, obtained in~\cite{st95}, tend to the
speed of light.

Another surprise is that the general theory shows there exist
gravitational metric components which are $C^{1,1}$ across a
lightlike surface for which $\mbox{div}\,G=0$ does not hold in the
weak sense. This implies that the matching of the (generalized)
second fundamental form alone is not a sufficient condition for
conservation conditions to hold across the interface.

\section{The Central Problem: Lightlike Shock-Waves}
The goal of the work here is to construct shock-wave solutions of
the Einstein equations which move at the speed of light. In
Einstein's general relativity all physically possible spacetimes
correspond to solutions of Einstein's Equations
\begin{equation}\label{eq:einsteinintro}
G=\kappa T,
\end{equation}
which are metrics that describe the geometric structure of
spacetime. Equation~\eqref{eq:einsteinintro} is a tensor equation
representing 10 nonlinear partial differential equations where $G$
describes the geometric structure of spacetime, and $T$ represents
the matter which is the source of the gravitational field.

Shock-waves were first studied in compressible, non-viscous gas
flow as discontinuities that form in the fluid quantities
pressure, density, etc. The mathematical theory of shock-waves is
contained in the study of hyperbolic conservation laws, and
applies to much more that just gas dynamics. Mathematical shock
wave theory not only models phenomena such as the sonic boom
created by a fast moving plane, but the same theory also
incorporates the propagation of ``gridlock'' in traffic flow, the
leading edge of a nuclear explosion, flame fronts in combustion,
and separation of boundaries between chemical species in
chromatography~\cite{temp99}. If we take the covariant divergence
of equation~\eqref{eq:einsteinintro}, then it turns out that
\begin{equation}\label{eq:divintro}
\mbox{div}\, T=0.
\end{equation}
In the limit of low velocities and weak gravitational fields
equation~\eqref{eq:divintro} reduces to the classical compressible
Euler equations which is a hyperbolic conservation law in gas
dynamics.

The shock-wave solutions considered in this dissertation are
constructed by matching two metrics (geometries) across a surface
embedded in spacetime. In their well-known 1939 paper~\cite{os39}
describing gravitational collapse of a star, Oppenheimer and
Snyder gave the first example, which had interesting dynamics, of
a solution of the Einstein equations obtained by matching two
solutions across a surface~\cite{st94}. With their simplifying
assumption that the pressure is zero, the surface across which the
metrics are matched is not shock-wave, but a \emph{contact
discontinuity}\index{contact discontinuity} which is a
discontinuous solution where neither mass or momentum
cross~\cite{st94}. Smoller and Temple extended the
Oppenheimer-Snyder model of gravitational collapse to the case of
non-zero pressure in their 1994 paper~\cite{st94}. Their theory
was based on work done by Israel in his 1966 paper~\cite{isr66},
which related the second fundamental form across a shock surface
to the Rankine-Hugoniot jump conditions across the shock. Then, in
later work~\cite{st95}, Smoller and Temple constructed a family of
exact, spherically symmetric, shock-wave solutions of the Einstein
equations by matching two spherically symmetric metric across a
surface. The shock surfaces in Smoller and Temple's work are
assumed to be non-lightlike, that is, moving slower than the speed
of light. The contribution here will be to extend their results
given in~\cite{st94} and~\cite{st95} to also incorporate the case
when the surface is lightlike.

\subsection{The Two Difficulties of the Lightlike Case}
In the lightlike case the mathematical machinery used in the
sub-luminal case breaks down in two areas. First, the induced
metric on a lightlike surface is degenerate, and the second is
that the second fundamental form,
\begin{equation}\label{eq:2formintro}
K=-\nabla_X \n,
\end{equation}
cannot be used to describe the dynamics of the surface in the
ambient spacetime. A lightlike surface is characterized by the
length of the normal vector being zero. Note the in this case the
normal vector is non-zero, but the metric is not positive
definite. A normal of zero length is orthogonal to itself, hence
it lies in the tangent space of the surface. Since the normal is
no longer transverse to the surface, the second fundamental form
no longer gives information about how the surface in embedded in
the spacetime manifold.

The degeneracy problem is dealt with by considering the problem in
the context of the whole spacetime manifold where the metric is
not degenerate. Unfortunately, the failure of the second form is
not fixed so easily. To rectify this failure we make use of an
idea originated by Barrab\`{e}s and Israel in~\cite{isr91} in
which they define a generalized second fundamental form in terms
of a vector $\N$ transverse to the surface.

The work of Barrab\`{e}s and Israel in~\cite{isr91} focuses on the
dynamics of \emph{surface layers}\index{surface layer} whose
theory is also based on matching two metrics across a surface. A
surface layer differs from a shock-wave in that a shock surface is
characterized by a jump in the density of the fluid across the
surface, in contrast, the density becomes infinite in a surface
layer. Mathematically, the Einstein tensor $G$, which comprises
the right hand side of equation~\eqref{eq:einsteinintro}, contains
no delta-function singularities across a shock surface, but across
a surface layer $G$ does contain a delta-function singularity. We
also note here that Barrab\`{e}s and Israel use a scalar version
of the second fundamental form where we use the form given in
equation~\eqref{eq:2formintro} which is a tangent vector in the
surface. Regarding the lightlike case, Barrab\`{e}s and Israel
state that because of the breakdown in the second fundamental
form, the lightlike case is a relatively neglected area of surface
dynamics which remains imperfectly understood~\cite{isr91}.

\section{How This Document is Organized}
In chapter~\ref{chp:2} we give an introduction to general
relativity, and a derivation of the second fundamental form. Then
in chapter~\ref{chp:3} we give a generalization of the second
fundamental form based on the idea of using a transverse vector in
place of the normal vector. Continuing on, we state and prove the
main result of this dissertation giving a set of equivalent
conditions which yield the existence of shock solutions in the
Einstein equations when the metric is Lipschitz continuous across
a hypersurface, and we finish the chapter by proving a similar
result involving spherically symmetric metrics. In
chapter~\ref{chp:4} we give an exact, spherically symmetric,
lightlike shock-wave solution of the Einstein equations based on
Smoller and Temple's work in~\cite{st95}. Chapter~\ref{chp:5}
concludes the results.

%% file: back.tex
% ---Background---------------------------------------------------
% General Relativity and the Einstein Equations
\label{chp:2}

%\section{Introduction}
The purpose of this chapter is to make this dissertation
self-contained for anyone with a good understanding of
differential geometry. We construct the geometry of spacetime from
physical assumptions, derive the stress-energy tensor for a
perfect fluid, define the Einstein equations\index{Einstein
equations}, and conclude the chapter with an introduction of the
second fundamental form.

The underlying theme of general relativity is the idea that the
physics of our universe and its geometry are interdependent in
such a way that they cannot be separated from each other. The
Einstein field equations,
$$G=\kappa T,$$
is a mathematical formulation of this idea with respect to the
physics of gravity. The Einstein equations match $G$, a second
order differential operator on the metric that is related to the
Riemannian curvature tensor, to $T$, the stress-energy tensor of
the fluid or matter in the region of spacetime that is being
considered.

\section{The Geometry of Spacetime}
\label{s:spacetime} A theory of gravitation begins with a notion
of space and time. \emph{Spacetime}\index{spacetime} is a
continuum of \emph{events}\index{event} each of which has three
spacial components, and one time component. If we give each event
a name in terms of these components, say $x=(x^0,x^1,x^2,x^3)$
where $x^0=ct$ with $c$ the speed of light in a vacuum, $t$
denoting time, and $x^a$ for $a=1,2,3$ representing the spacial
components, then we can represent each event in as a coordinate
$x$. The geometry of spacetime is the structure of how this
continuum of events fit together. We shall see, it is the
principles from physics that will determine this structure.

Mathematically spacetime is a four dimensional manifold $M$
equipped with a metric\index{metric} $g$. It is the metric that
carries the information about the rates at which clocks run, and
the distances between points. It is also the metric that describes
the geometry. The question now is how do we determine a metric
which models spacetime or subset of it?

\subsection{Special Relativity}
In order to proceed we need the notion of an \emph{inertial
coordinate frame} or sometimes called an \emph{inertial observer}.
A coordinate system is inertial if any two events in the
coordinate system share the same clock or time coordinate, and the
geometry of space at any constant time is Euclidean~\cite{schutz}.
Consider a spacetime without gravity, that is, a spacetime
containing negligible mass. This is the realm of special
relativity which can be deduced by the two following postulates:

\begin{postulate}[Principle of Relativity]
\index{Principle of Relativity} \label{pt1} Let $x$ be an inertial
coordinate system. Then any other coordinate system $\bar{x}$
which moves uniformly, and is non-rotating relative to $x$ is also
an inertial coordinate system. Furthermore, "The laws of nature
are in concordance for all inertial systems"\cite{einst22}.
\end{postulate}

\begin{postulate}[Propagation of light]
\index{Propagation of light} \label{pt2} The speed of light in a
vacuum $c$ is observed to be the same from any inertial coordinate
system.
\end{postulate}

The principle of relativity is equivalent to the statement that
any body in uniform motion remains in that state unless acted upon
by an external force~\cite{schutz}. The universality of the speed
of light means that two observers, one moving with non-zero
uniform velocity with respect to the other, will each observe a
single light ray moving through vacuum to be travelling at same
speed $c$.

\subsubsection{Measuring Distance in a Special Relativistic Spacetime}

\begin{figure}[hbtp]
  \centering
  \includegraphics[width=.65\textwidth]{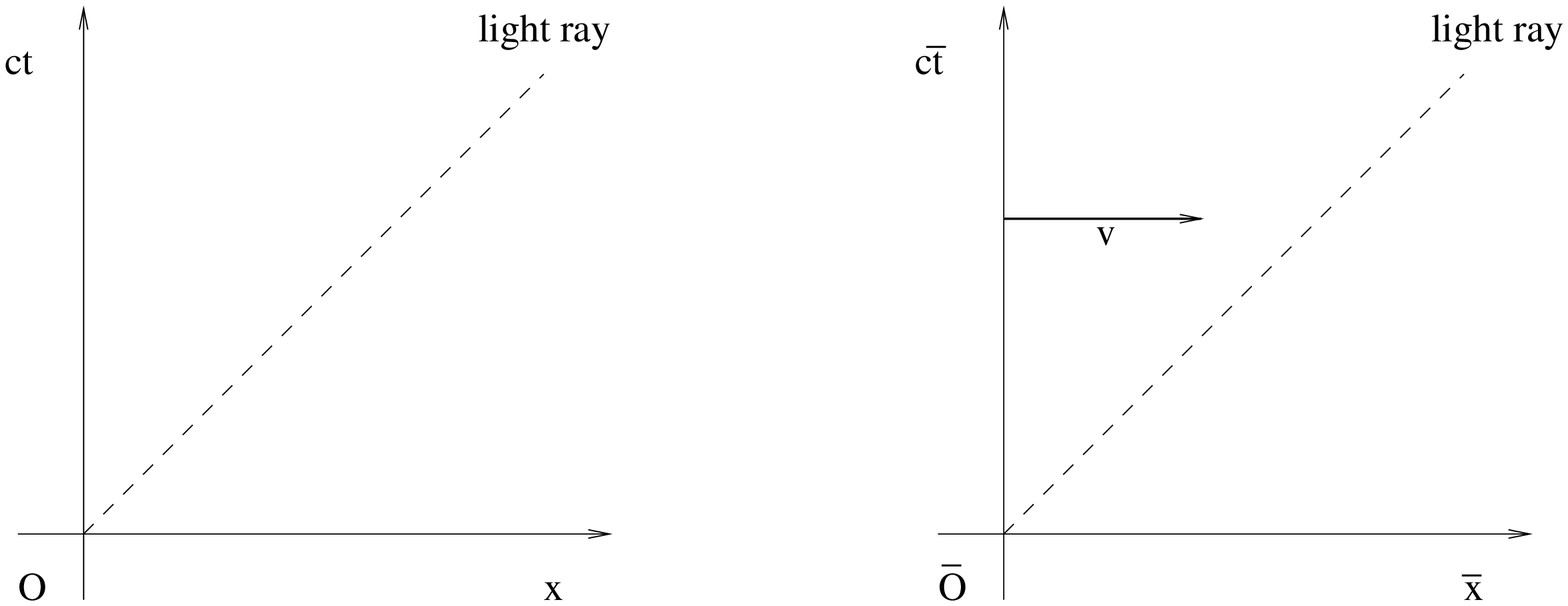}
  \caption{}\label{inertfram}
\end{figure}

Consider two inertial frames $O$ and $\bar{O}$, each having only
one spatial dimension for simplicity, with coordinates $(x,ct)$,
and $(\bar{x},c\bar{t})$ respectively. Assume $\bar{O}$ is moving
along the $x$ axis with relative uniform velocity $v$ see
figure~\ref{inertfram}.
% ---Figure for inertial frames-----------------------------------
% ---End Figure for inertial frames-------------------------------
Let two events be connected by a light ray. Then in the $O$ frame
the squared distance between the events is given by
\begin{equation}\label{eq:sqr1}
\Delta s^2=-(\Delta ct)^2 +(\Delta x)^2=0.
\end{equation}
By the propagation of light law the same two events, in the
$\bar{O}$ frame , also satisfy
\begin{equation}\label{eq:sqr2}
\Delta \bar{s}^2=-(\Delta c\bar{t})^2 +(\Delta \bar{x})^2=0.
\end{equation}
Now, assume that coordinate transformation from $O$ to $\bar{O}$
is linear, and that their origins coincide. Then if follows that
$$\Delta s^2 = \Delta \bar{s}^2$$ for any two events in spacetime
for any inertial coordinate frames $O$ and $\bar{O}$ whose origins
coincide, and are related by a linear transformation.

Recall that the metric $g$ is required to give data on the rates
at which clocks run, and the distances between spatial points.
Since any data determined by the metric also needs to satisfy
postulate~\ref{pt1}, we require that $g$ be a coordinate
independent quantity. Therefore, $g$ should incorporate the
coordinate invariant quantity $\Delta s^2,$ which composed of
squared displacements between events.

\begin{figure}[hbtp]
  \centering
  \includegraphics[width=.6\textwidth]{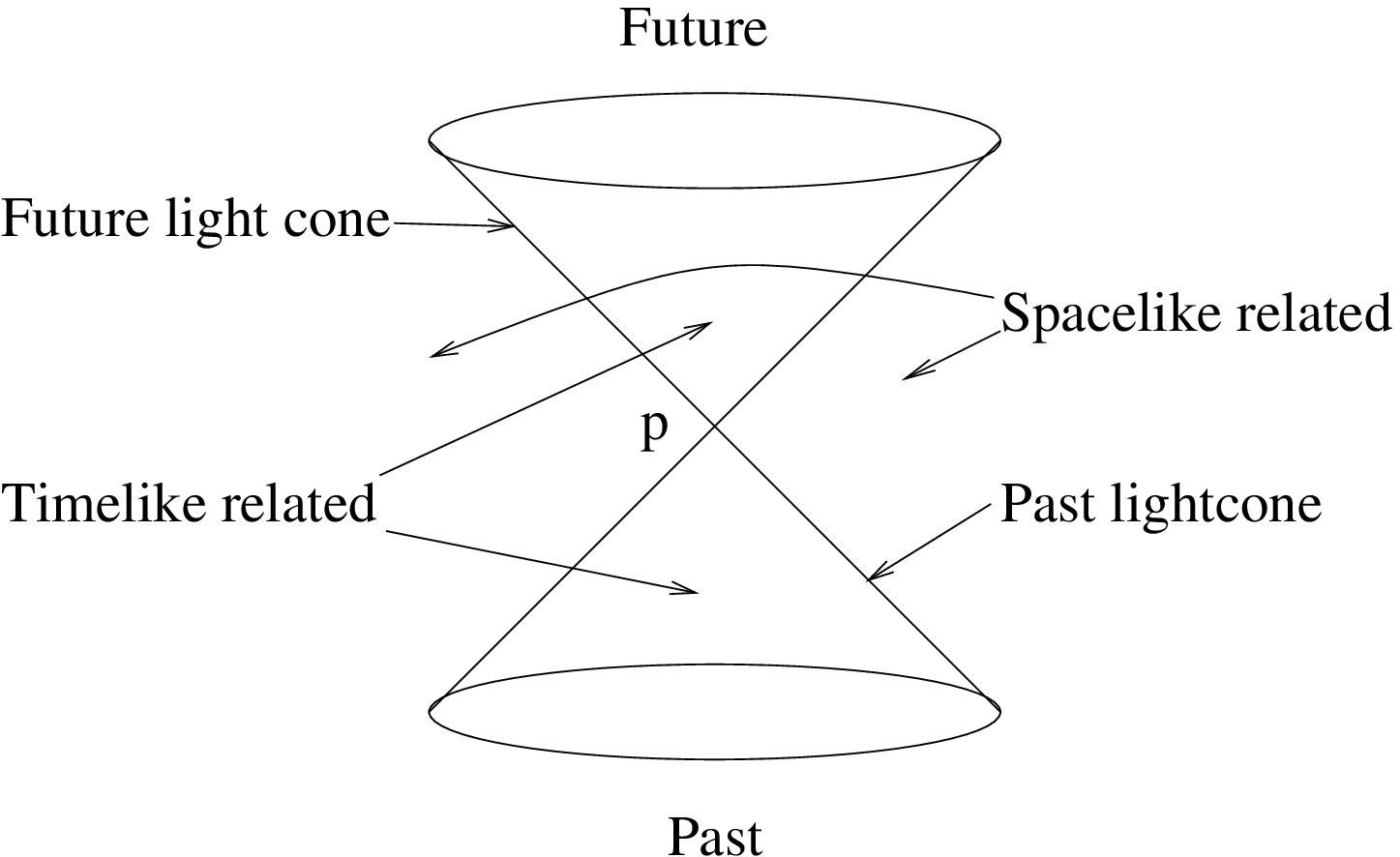}
  \caption{}\label{fig:causal}
\end{figure}

What is the formula for $g$ in a given inertial frame? Let
$\vec{A}$ denote a displacement vector between two events in some
inertial coordinate system. Then the squared distance between the
two events in the coordinates is given by
$$g(\vec{A},\vec{A})=-(A^0)^2+(A^1)^2+(A^2)^2+(A^3)^2,$$
which consistent with the squared distance in~\eqref{eq:sqr2}.
This leads to the inner product between any two vectors $\vec{A}$,
and $\vec{B}$ given by
$$g(\vec{A},\vec{B})=-A^0B^0+A^1B^1+A^2B^2+A^3B^3.$$
Notice this is very similar to the standard dot product with the
only difference being the negative sign on the product of the
first components. We can also write this metric in the following
form:
\begin{equation}\label{eq:mink}
g(\vec{A},\vec{B})=
\begin{pmatrix}A^0 \\ A^1 \\ A^2 \\ A^3\end{pmatrix}^T
\begin{pmatrix}
-1 & 0 & 0 & 0 \\
 0 & 1 & 0 & 0 \\
 0 & 0 & 1 & 0 \\
 0 & 0 & 0 & 1
\end{pmatrix}
\begin{pmatrix}B^0 \\ B^1 \\ B^2 \\ B^3\end{pmatrix}.
\end{equation}
In this form we can see that $g$ has a signature of $+2$. A metric
with this signature is said to have a \emph{Lorentzian
signature}\index{Lorentzian signature}. We will only deal with
metrics of Lorentzian signature. Furthermore, the metric given
in~\eqref{eq:mink} is said to be
\emph{Minkowskian}\index{Minkowski metric}, and is denoted by
$\eta_{\alpha\beta}=\mbox{diag}(-1,1,1,1)$. The inner product $g$
is not positive definite. To see this consider
$$g(\vec{A},\vec{A})=\langle(1,1,0,0),(1,1,0,0)\rangle=-1+1=0.$$
However, we can always find a coordinate basis so that $g$ is of
the form $\eta_{\alpha\beta}$. Therefore $g$ is
\emph{non-degenerate}, that is, $g(\vec{A},\vec{B})=0$ for all
vectors $\vec{B}$ if and only if $\vec{A}=0$.
\begin{definition}\label{def:vec}
A vector $\vec{A}$is said to be \emph{spacelike}\index{spacelike
vector} if $g(\vec{A},\vec{A})>0,$ is
\emph{lightlike}\index{lightlike vector} if
$g(\vec{A},\vec{A})=0$, is \emph{timelike}\index{timelike vector}
if $g(\vec{A},\vec{A})<0$.
\end{definition}
From this we can give the \emph{causal structure}\index{causal
structure} of spacetime, that is the causal relationship of an
event to other events. The future or past of an event $p$ lies
inside the light cone of $p$, see figure~\ref{fig:causal}. If an
event $p$ lies at the origin of an inertial frame, then the future
and past of $p$ is given by all events $q$ such that the
displacement vectors between $p$ and $q$ are timelike.
% ---Figure for Causal Structure----------------------------------
% ---End Figure for Causal Structure------------------------------

\begin{summary}[Special Relativity]
From postulates~\ref{pt1} and~\ref{pt2} we were able to describe
how the squared distance $\Delta s^2$ should be measured, and that
it is invariant under coordinate transformations between inertial
frames. From the idea of $\Delta s^2$ we defined the inner product
$g$ given in~\eqref{eq:mink}, and designated it as our metric for
spacetime with no gravity. From this we ascertained that $g$ is a
metric with Lorentzian signature.
\end{summary}

\subsection{General Relativity}
We now consider spacetime with gravity. We are motivated by the
fact that objects under the influence of gravity move along free
fall paths. To see this imagine an observer in a spaceship with
the engines turned off, and orbiting the earth under the influence
of only earth's gravity at a distance $r_0$ above the earth. We
are assume the earth is a perfect uniform sphere, and ignoring the
fact that the earth is rotating. Then an observer inside the
spacecraft now feels weightless, and can position themselves in
such way that they remain stationary in the middle of the ship so
that they are not moving relative to the walls of the ship. The
ship, and the observer are moving along a free fall path of
earth's gravity. Such a coordinate frame is said to be
\emph{comoving}\index{comoving} with the metric.

\begin{figure}[hbtp]
  \centering
  \includegraphics[width=.6\textwidth]{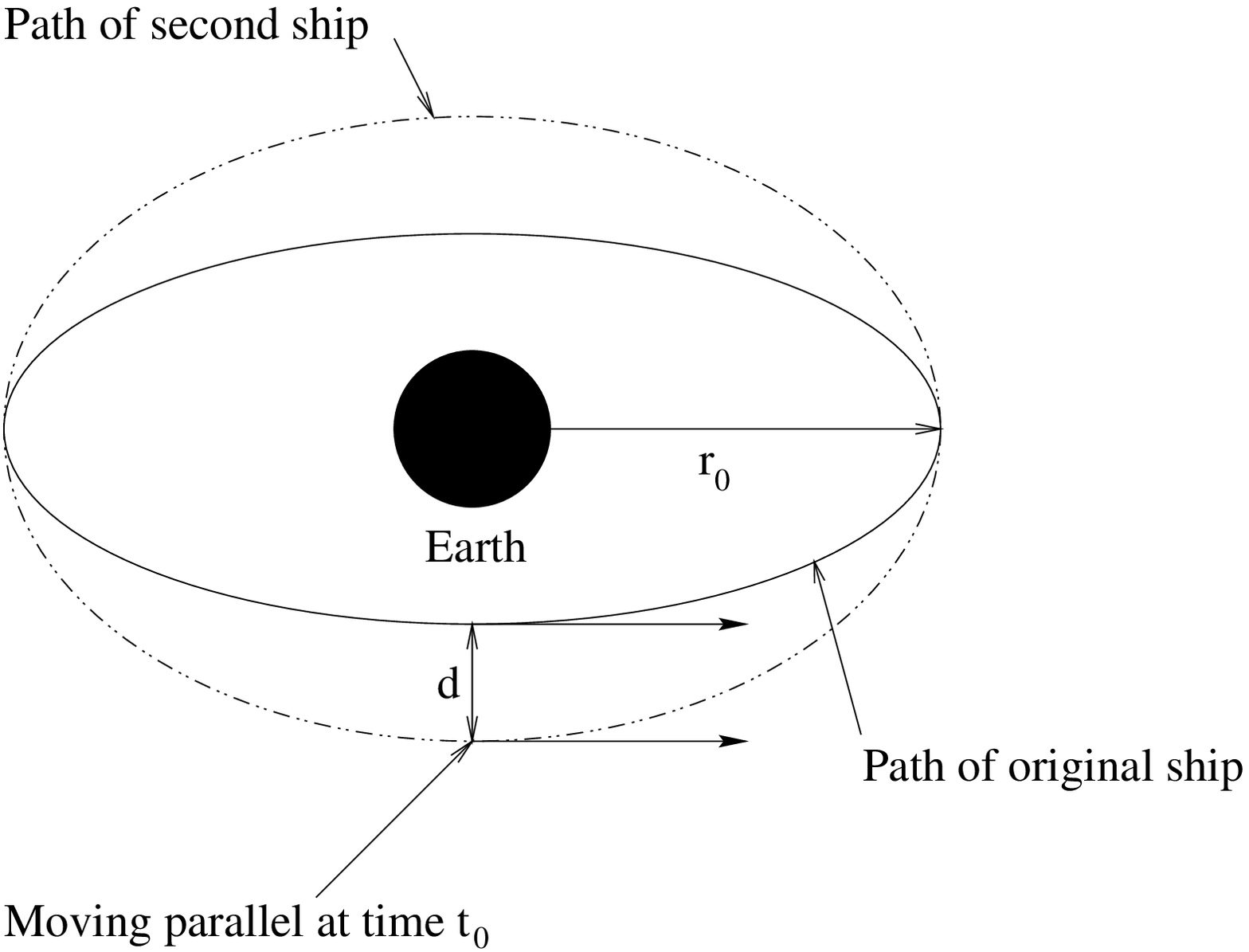}
  \caption{}\label{fig:orbit}
\end{figure}

How does the observer know that they are under the influence of a
gravitational force? Consider a second ship in orbit the same
distance $r_0$ above earth, but positioned a distance $d$ to the
right. Also, assume the second ship is moving with the same
velocity parallel to the original ship at some instant of time say
$t_0$, see figure~\ref{fig:orbit}.
% ---Figure for Orbit---------------------------------------------
% ---End Figure for Orbit-----------------------------------------
As seen in figure~\ref{fig:orbit}, the distances between the two
ships is not constant, and their paths will even cross. In this
space Euclid's parallel axiom does hold, and spacetime is curved.

Mathematically, for the local coordinate system, (inside the
original ship) the Minkowskian metric,
$g_{\alpha\beta}=\eta_{\alpha\beta}$, is appropriate just as it
was in special relativity. Also, at the point of spacetime where
the ship is located, we have $g_{\alpha\beta,\gamma}=0$. This
defines a \emph{locally inertial}\index{locally inertial} frame or
also called \emph{locally Lorentzian}\index{locally Lorentzian}.

A logical question to ask at this point is, does there exist a
global inertial frame when gravity is present? The answer is no.
We have already seen from our example that the two ships cannot be
described by the same inertial frame because they lie in a curved
space which is not Euclidean in it's spatial components violating
the definition of an inertial frame. This violation is due to the
non-uniformity of the earth's gravitational field. As we expect
with a curved space, we can only consider local frames
diffeomorphic to $\mathbb{R}^4$.

The idea that objects moving in a gravitational field are
described by comoving local inertial systems takes into account
the \emph{Equivalence Principle}\index{Equivalence Principle}. The
Equivalence Principle says that all bodies move the same way in a
gravitational field~\cite{wald}. It is the equivalence principle
that motivated Einstein to formulate the theory of general
relativity.

The implication that under the influence of gravity particles move
in a curved space, which can only be described in the context of
locally inertial frames, is that we have to use the mathematical
machinery of differential geometry to characterize the dynamics of
particles under the influence of gravity. This is in contrast to
special relativity where the curvature is zero, since their is no
acceleration of particles due to gravity, and one coordinate frame
is sufficient to describe spacetime. Therefore, we will present
the concepts from differential geometry that are needed to
describe a curved spacetime.

\subsection{The Metric in a Curved Spacetime}
The metric completely defines the curvature of the spacetime
manifold, and so we begin by constructing the spacetime metric $g$
for a possibly curved spacetime. Earlier we used displacement
vectors to construct the Minkowski metric in the flat spacetime of
special relativity. If the displacements are infinitesimally
small, then we can associate them with tangent vectors, which are
defined by directional derivatives. Let $p$ be in a point in the
spacetime manifold $M$, and suppose
$$x=(x^0,x^1,x^3,x^3):M\To \mathbb{R}^4,$$
are coordinates on $M$ whose origin is at $p$, where $x^0=ct$,
$c=1$, and $x^i$ for $i=1,2,3$ denote the spatial coordinates.
Also let $\mathcal{F}$ denote the set of all $C^\infty$ functions
that take smooth curves on $M$ into $\mathbb{R}$.
\begin{equation}\label{eq:diagram}
\begin{CD}
M@<x^{-1}<<\mathbb{R}^4\\
 @V{f}VV \\
\mathbb{R}
\end{CD}
\end{equation}
Then a tangent vector $X$ at $p$ in $M$ is defined to be a map
$X:\mathcal{F}\to \mathbb{R}$. For any coordinate system $x$ on a
neighborhood of $p$, there exist a coordinate basis given by
\begin{equation}\label{eq:tanvec}
X_\alpha\left(f(p)\right)=\left.\frac{\partial}{\partial
x^\alpha}(f\circ x^{-1})\right|_{x(p)}.
\end{equation}
See diagram in~\eqref{eq:diagram}. That is, $\{\partial/\partial
x^\alpha\}$ gives a basis for $\tanm$, the tangent space of $M$ at
$p$.

We now introduce the \emph{Einstein summation
convention}\index{Einstein summation convention} which says that
for any expression,  equal up and down indices are summed over all
possible values the index can take. For example,
using~\eqref{eq:tanvec}, if $X\in\tanm$, then
$$X=\sum_{\alpha=0}^3 X^\alpha \frac{\partial}{\partial
x^\alpha} = X^\alpha \frac{\partial}{\partial x^\alpha}.$$ Notice
that the upper index in the denominator is considered a down
index.

With tangent vectors defined as directional derivatives, they are
in one-to-one correspondence with displacement vectors. Since the
metric is quadratic in displacements we define $g$ as a map where
$$g:\tanm \times \tanm \To \mathbb{R}.$$
Also, we assume that $g$ is symmetric, and non-degenerate as was
the case in special relativity.

The metric $g$ being symmetric guarantees, that in coordinates,
the components of the metric $g_{\alpha\beta}$ make up a symmetric
matrix, and therefore there will always exist a linear
transformation to another coordinate system that will take
$g_{\alpha\beta}$ to $\eta_{\alpha\beta}=\mbox{diag}(-1,1,1,1),$
the Minkowski metric\index{Minkowski metric}. Furthermore, it can
also be shown that the metric $g$ is locally
Lorentzian\index{locally Lorentzian} or locally
inertial\index{locally inertial}, that is, given $p$ in $M$, there
exist coordinates $x$ whose origin is at $p$ such that
$$g_{\alpha\beta}(p)=\eta_{\alpha\beta} \quad \mbox{for} \quad
\alpha,\beta=0,1,2,3,$$
$$g_{\alpha\beta,\mu}(p)=0\quad \mbox{for} \quad
\alpha,\beta,\gamma=0,1,2,3,$$ and
$$g_{\alpha\beta,\mu\nu}(p)\neq0,$$
at the point $p$ for at least some values of $\alpha,$ $\beta$,
$\mu$ and $\nu$ if spacetime is not flat, that is, there exists a
gravitational field. This corresponds exactly with the notion of
locally inertial\index{locally inertial} frames given in the
previous section.

\subsection{Tensors in Spacetime}
An alternative version of the Equivalence
Principle\index{Equivalence Principle} is the \emph{General
Covariance Principle}\index{General Covariance Principle} which
motivates the use of tensors to measure physical quantities that
depend linearly on displacements. The General Covariance
Principle\index{General Covariance Principle} states that an
equation holds in a gravitational field if the two following
conditions are met~\cite{weinberg}:
\begin{enumerate}
\item The equation holds in the absence of a gravitational field;
that is, it is consistent with postulates~\ref{pt1} and~\ref{pt2}.

\item The equation is covariant\index{covariant}, which means that the equation
holds under any coordinate transformation $x\to y$.
\end{enumerate}
We take the following definition of a tensor from~\cite{wald}.
Given a finite vector space $V$, and denoting its dual space by
$V^*$, a \emph{tensor, T, of type} $(k,l)$\index{tensor} on the
space $V$ is a multilinear map
\begin{equation}\label{eq:tensorkl}
T:(V^*)^k \times V^l \To \mathbb{R}.
\end{equation}
Note that our definition of tensor makes no mention vector or dual
vector components. A tensor gives the same real number for a
particular set of vectors and dual vectors independent of the
coordinates the components are computed in.

Our first example of a tensor is the metric $g$ which is a $(0,2)$
tensor, that can be written in terms of coordinate basis one-forms
as
\begin{equation}\label{eq:metcomp}
ds^2=g_{\alpha\beta}\,dx^\alpha \otimes dx^\beta
=g_{\alpha\beta}\,dx^\alpha dx^\beta.
\end{equation}
It is customary to drop the exterior product sign between the
one-forms $dx^\alpha$. Furthermore, the metric provides a mapping
between vectors and one-forms (dual vectors) at every point. Thus,
given a vector field $X(p)$, there is a unique one-form field
given by $\tilde{X}(p)=g(X(p),\quad )$.

Using this map between vectors and one-forms we can construct a
unique and useful one-form basis. Suppose $\{\partial/\partial
x^\alpha\}$ is a coordinate basis for $\tanm$. Then
$\{dx^\alpha\equiv g(\partial/\partial x^\alpha,\quad)\}$ gives a
one-form basis.

Since $g$ is non-degenerate, there exists an inverse which takes
one-forms to vectors. If the indices of the basis one-forms are
``up" and ``down" for basis vectors, then we can use the Einstein
summation convection\index{Einstein summation convention} to keep
track of these mappings. Notice that for the Einstein summation
convention to work we have to have the indices of vector
components ``up," and the indices on components of the one-forms
``down." For example, see equation~\eqref{eq:metcomp}. The
components for the inverse of $g$ are denoted by
$g^{\alpha\beta}$. Notice that
$g_{\alpha\sigma}g^{\sigma\beta}=\delta_\alpha^\beta$ as we would
expect from matrix multiplication. Using the metric $g$ as an
invertible map from vectors to one-forms can be executed by
raising and lowering an index of a tensor by contracting the index
with metric. For example, we can map the vector $X$ to a one-form
as
$$g_{\alpha\sigma}X^\sigma=X_\alpha.$$
Another example mapping  a $(3,1)$ tensor to a $(2,2)$ tensor:
$$g^{\beta\sigma}R^{\alpha}_{\sigma\mu\nu}=R^{\alpha\beta}_{\mu\nu}.$$

The components of a $(k,l)$ tensor given in one coordinate system
can be written in terms of another set of coordinate system via
the \emph{tensor transformation law}\index{tensor transformation
law}:
\begin{equation}\label{eq:tensortran}
T^{\alpha_1,...,\alpha_k}_{\beta_1,...,\beta_l}
=T^{\mu_1,...,\mu_k}_{\nu_1,...,\nu_l}\frac{\partial
y^{\alpha_1}}{\partial x^{\mu_1}}\cdots\frac{\partial
y^{\alpha_k}}{\partial x^{\mu_k}}\frac{\partial
x^{\nu_1}}{\partial y^{\beta_1}}\cdots\frac{\partial
x^{\nu_l}}{\partial y^{\beta_l}}.
\end{equation}
Here we have written a $y$-coordinate tensor component
$T^{\alpha_1,...,\alpha_k}_{\beta_1,...,\beta_l}$ in terms of
$x$-coordinates. Also, note that the Jacobian matrix satisfies
$\frac{\partial x}{\partial y}=\left(\frac{\partial y}{\partial
x}\right)^{-1}.$

\subsection{Parallel Transport and the Derivative Operator}
Motivated by the idea of describing the curvature of spacetime
intrinsically, as opposed to describing spacetime as an embedding
in some other space, we will define curvature in terms of
\emph{parallel transport}\index{parallel transport}. Intuitively,
a vector field $Y$ defined on every point along a curve is said to
be parallel transported along the curve if the vectors of $Y$ are
parallel, and are of equal length at infinitesimally close points.
Mathematically, in a locally Lorentzian frame at a point $p$, the
components of the $Y$ vectors must stay constant along the curve
near $p$. If we let $x(\xi)$ be a parameterization of the said
curve, and denote the its tangent by $X=dx/d\xi,$ then
$dY^\alpha/d\xi=0$ at $p$. However,
\begin{equation}\label{eq:paratrans}
\frac{d Y^\alpha}{d\xi} =\frac{dx}{d\xi}\frac{dY^\alpha}{dx}
=X^\beta \frac{\partial Y^\alpha}{\partial x^\beta}=0.
\end{equation}
Equation~\eqref{eq:paratrans} leads to the idea that in order to
define parallel transport for any coordinate frame, not just an
locally inertial system, we require a notion of how to take
derivatives of vector fields.

Consider the vector field $Y$ in an arbitrary coordinate frame,
not necessarily Lorentzian. We differentiate $Y$ as follows:
\begin{equation}\label{eq:dY}
\nabla_\beta Y \equiv\frac{\partial}{\partial
x^\beta}\left(Y^\alpha \frac{\partial }{\partial x^\alpha}\right)
=\frac{\partial Y^\alpha}{\partial
x^\beta}\frac{\partial}{\partial x^\alpha} +Y^\alpha
\frac{\partial }{\partial x^\beta}\frac{\partial }{\partial
x^\alpha}.
\end{equation}
Now, $$\frac{\partial }{\partial x^\beta}\frac{\partial }{\partial
x^\alpha}$$ is a vector and can be written in terms of the
coordinate basis, that is,
\begin{equation}\label{eq:chrisof1}
\frac{\partial }{\partial x^\beta}\frac{\partial }{\partial
x^\alpha}=\G^\mu_{\alpha\beta}\frac{\partial }{\partial x^\mu},
\end{equation}
where $\G^\mu_{\alpha\beta}$ is called a \emph{Christoffel
symbol}\index{Christoffel symbol} which is yet to be determined.
Therefore, we can write equation~\eqref{eq:dY} as
\begin{equation}\label{eq:covarY}
\nabla_\beta Y \equiv Y_{;\beta} =\frac{\partial
Y^\alpha}{\partial x^\beta} \frac{\partial}{\partial x^\alpha}
+Y^\alpha \G^\mu_{\alpha\beta} \frac{\partial }{\partial x^\mu}
\equiv Y^\alpha_{\hspace{5pt},\beta}\frac{\partial}{\partial
x^\alpha} +Y^\alpha \G^\mu_{\alpha\beta} \frac{\partial }{\partial
x^\mu}.
\end{equation}
The differential operator $\nabla$ is called the \emph{covariant
derivative}\index{covariant derivative}. Given arbitrary vector
fields $X$ and $Y$ the covariant derivative is defined as
\begin{equation}\label{eq:covarXY}
\nabla_X Y=X^\alpha \nabla_\alpha Y = X^\alpha
Y^\beta_{\hspace{5pt};\alpha}\, \frac{\partial}{\partial x^\beta}.
\end{equation}

Notice that our definition of covariant derivatives did not
involve the metric. However, if we recall that the metric $g$ maps
vectors into one-forms, then it would seem that the metric would
have something to do with how their derivatives are related. In a
locally Lorentzian coordinate frame we know that for any vector
field $Y,$
$$Y_{,\alpha} = Y_{;\alpha},$$
since the derivatives of the basis vectors are zero just as in
special relativity. The same relation holds true for any tensor in
a Lorentzian frame including the metric. Therefore, in a locally
Lorentzian frame, the covariant derivative of metric components
are given by
\begin{equation}\label{eq:metflat}
g_{\alpha\beta;\gamma}=g_{\alpha\beta,\gamma}=0.
\end{equation}
Invoking the General Covariance Principle, we find that
\begin{equation}\label{eq:met0}
g_{\alpha\beta;\gamma}=0
\end{equation}
holds in any coordinate frame. Using this result any Christoffel
symbol can be written in terms of the metric. To do this we must
first show that $\G^\sigma_{\mu\nu}=\G^\sigma_{\nu\mu}.$

Consider a scalar field $\phi$ in a Lorentzian coordinate frame.
Then $\nabla \phi$ is a one-form with components $\phi_{,\beta}$.
In the Lorentzian frame
\begin{equation}\label{eq:chrissymm}
\phi_{,\alpha;\beta} =\phi_{,\alpha,\beta} =\phi_{,\beta,\alpha}
=\phi_{,\beta;\alpha},
\end{equation}
since partial derivatives commute. Once more, invoking the General
Covariance Principle, we have that the symmetry in
equation~\eqref{eq:chrissymm} holds in any coordinate system.
Therefore, in any coordinate frame
$$\phi_{,\alpha;\beta}-\phi_{,\beta;\alpha}
=\phi_{,\alpha,\beta}-\G^\sigma_{\alpha\beta}\phi_{,\sigma}
-\phi_{,\beta,\alpha}+\G^\sigma_{\beta\alpha}\phi_{,\sigma}=0,$$
and we have that
\begin{equation}\label{eq:Gsym}
\G^\sigma_{\alpha\beta}=\G^\sigma_{\beta\alpha}.
\end{equation}
Then, since $g_{\alpha\beta;\gamma}=0$ in any coordinate system,
we can write
\begin{equation}\label{eq:Gformu}
\begin{aligned}
g_{\alpha\beta,\gamma} &=\G^\sigma_{\alpha\gamma}g_{\sigma\beta}
+\G^\sigma_{\beta\gamma}g_{\sigma\alpha},\\
g_{\gamma\alpha,\beta} &=\G^\sigma_{\gamma\beta}g_{\sigma\alpha}
+\G^\sigma_{\alpha\beta}g_{\sigma\gamma},\\
g_{\beta\gamma,\alpha} &=\G^\sigma_{\beta\alpha}g_{\sigma\gamma}
+\G^\sigma_{\gamma\alpha}g_{\sigma\beta}.
\end{aligned}
\end{equation}
Thus
\begin{multline}\label{eq:Gcalc}
-g_{\alpha\beta,\gamma} +g_{\gamma\alpha,\beta}
+g_{\beta\gamma,\alpha}
\\=-\G^\sigma_{\alpha\gamma}g_{\sigma\beta}
-\G^\sigma_{\beta\gamma}g_{\sigma\alpha}
+\G^\sigma_{\gamma\beta}g_{\sigma\alpha}
+\G^\sigma_{\alpha\beta}g_{\sigma\gamma}
+\G^\sigma_{\beta\alpha}g_{\sigma\gamma}
+\G^\sigma_{\gamma\alpha}g_{\sigma\beta}.
\end{multline}
Using the symmetry property~\eqref{eq:Gsym},
equation~\eqref{eq:Gcalc} the Christoffel symbol\index{Christoffel
symbol} can be given in terms of the metric as
\begin{equation}\label{eq:cs}
  \Gamma^{\sigma}_{\alpha\beta}=
  \frac{1}{2} g^{\sigma\gamma}\{-g_{\alpha\beta,\gamma}+
  g_{\gamma\alpha,\beta}+g_{\beta\gamma,\alpha}\}.
\end{equation}

Now that we have a differential operator in hand we can define
parallel transport without reference to a coordinate system. This
is accomplished by expressing the parallel transport condition
given in equation~\eqref{eq:paratrans} in a covariant way.
Rewriting~\eqref{eq:paratrans} as
\begin{equation}\label{eq:paratranscov}
X^\beta Y^\alpha_{\hspace{5pt},\beta} =X^\beta
Y^\alpha_{\hspace{5pt};\beta}=0,
\end{equation}
which holds, if and only if,
\begin{equation}\label{eq:ptrans}
\nabla_X Y=0,
\end{equation}
see equation~\eqref{eq:covarXY}. Since equation~\eqref{eq:ptrans}
is a tensor equation, it must hold for any coordinate system.
Equation~\eqref{eq:ptrans} defines parallel transport of a vector
$Y$ along a curve with tangent vector $X$. Furthermore, given any
initial value of $Y$, equation~\eqref{eq:ptrans} has a unique
solution, and so given any initial vector $Y$ at a point on a
curve this defines a unique vector at any other point on the
curve. Then we can use the structure of parallel transport to map
$\tanm$ to $T_qM$ along a given curve from point $p$ to point $q$.
The covariant derivative, which gives this notion of identifying
tangent spaces from different points along a curve, is also
referred to as a \emph{connection}\index{connection}, and the
Christoffel symbols $\G^\sigma_{\alpha\beta}$ are referred to as
\emph{connection coefficients}\index{connection coefficients}.

With a definition of parallel transport we are in a position to
describe, mathematically, the ``free fall" paths in a
gravitational field. These paths correspond to the
\emph{geodesics}\index{geodesic} in the spacetime manifold. The
geodesics of the flat spacetime in special relativity are straight
lines. Qualitatively, geodesics can be thought of as the
straightest lines possible in a curved manifold, and they also
correspond to the ``straight" lines at the origin in a Lorentzian
coordinate frame. In flat manifold the straight lines are the only
curve in which the tangent vector parallel transports itself.
Generally, a geodesic is given by the equation
\begin{equation}\label{eq:geodesic}
\nabla_X X=0,
\end{equation}
which finds the curves where the tangent vectors parallel
transport themselves. If we let $\xi$ be a parameter which gets
mapped to curve, then
$$X^\alpha_{\hspace{5pt},\beta} =\frac{dx^\beta}{d\xi}\frac{\partial \;}{\partial
x^\beta}=\frac{d}{d\xi},$$ and we can write the geodesic
equation~\eqref{eq:geodesic}\index{geodesic equation} as
\begin{equation}\label{eq:geo2}
\frac{d}{d\xi}\left(\frac{dx^\alpha}{d\xi}\right)
+\G^\alpha_{\sigma\beta}\frac{dx^\sigma}{d\xi}\frac{dx^\beta}{d\xi}=0,
\end{equation}
which is a second order quasi-linear differential equation for
$x^\alpha(\xi)$. Equation~\eqref{eq:geo2} has a unique solution
given an initial position $x^\alpha_0=x(\xi_0)$, and an initial
direction $X^\alpha_0=X^\alpha(\xi_0)$ at $\xi_0$~\cite{schutz}.
This idea will be used in chapter~\ref{chp:3}, that is, we will be
given a point in spacetime and vector at that the point, and be
able to define a unique geodesic through the point.

\subsection{The Curvature Tensor} Having defined parallel
transport we are now able to define curvature. The \emph{Riemann
Curvature Tensor}\index{Riemann curvature tensor} describes the
failure of initially parallel geodesics to remain
parallel~\cite{wald}. This corresponds to the situation involving
the spacecraft orbiting earth in figure~\ref{fig:orbit}.

\begin{figure}[hbtp]
  \centering
  \includegraphics[width=.6\textwidth]{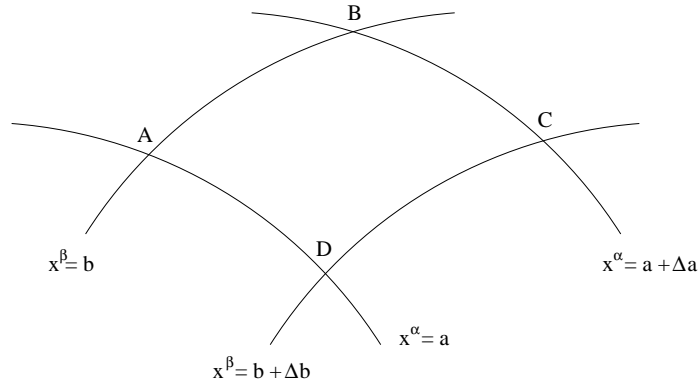}
  \caption{A piece of the coordinate grid.}\label{fig:coordgrid}
\end{figure}

If we parallel transport a vector $X$ around a coordinate grid
loop starting and ending at point $A$ in
figure~\eqref{fig:coordgrid},
% ---Figure for Coordinate Grid Loop-----------------------------
% ---End Figure Coordinate Grid----------------------------------
then it can be shown, using the general coordinates $x^\alpha$ and
$x^\beta$, the change of $X^\gamma$ by parallel transport along
$\Delta a\,\partial/\partial x^\alpha$, then $\Delta b\,
\partial/\partial x^\beta$, then $-\Delta a\, \partial /\partial
x^\alpha,$ and back to $A$ along $\Delta b\, \partial/\partial
x^\beta$ is given by
$$\Delta X^\gamma =\Delta a \Delta b \left(\G^\gamma_{\mu\alpha,\beta}
-\G^\gamma_{\mu\beta,\alpha}
+\G^\gamma_{\sigma\beta}\G^\sigma_{\mu\alpha}
-\G^\gamma_{\sigma\alpha}\G^\sigma_{\mu\beta} \right)X^\mu,$$
see~\cite[section 6.5]{schutz}. From this comes the definition of
the Riemann curvature tensor\index{Riemann curvature tensor}:
\begin{equation}\label{eq:riemcurv}
  R^{\mu}_{\hspace{4pt}\alpha\nu\beta}
  =\Gamma^{\mu}_{\alpha\beta, \nu} -\Gamma^{\mu}_{\alpha\nu,\beta}
  +\Gamma^{\mu}_{\sigma\nu} \Gamma^{\sigma}_{\alpha\beta} -
  \Gamma^{\mu}_{\sigma\beta} \Gamma^{\sigma}_{\alpha\nu}.
\end{equation}
The Riemann curvature tensor
$R^{\mu}_{\hspace{4pt}\alpha\nu\beta}$ transforms as tensor, via
the tensor transformation law~\eqref{eq:tensortran}, giving a
covariant measure of the second derivative of the metric $g$.

\begin{proposition}\label{riemp1}
The Riemann curvature tensor satisfies the following identities:
\begin{align}
R^{\mu}_{\hspace{4pt}\alpha\nu\beta}
 &=-R^{\mu}_{\hspace{4pt}\alpha\beta\nu}\label{eq:riemp1} \\
 R^{\mu}_{\hspace{4pt}[\alpha\nu\beta]}
&=R^{\mu}_{\hspace{4pt}\alpha\nu\beta}
+R^{\mu}_{\hspace{4pt}\beta\alpha\nu}
+ R^{\mu}_{\hspace{4pt}\nu\beta\alpha}=0 \\
 R_{\mu\alpha\nu\beta}
 &=g_{\mu\sigma}R^{\sigma}_{\hspace{4pt}\alpha\nu\beta}\\
 R_{\mu\alpha\nu\beta} &=R_{\nu\beta\mu\alpha}\label{eq:riemp1sym}
\end{align}
\end{proposition}
The \emph{Ricci tensor}\index{Ricci tensor} is a contraction of
the Riemann tensor, and is given by
\begin{equation}\label{eq:riccit}
  R_{\alpha\beta}=R^{\nu}_{\hspace{4pt}\alpha\nu\beta}.
\end{equation}
Due to the symmetries given in proposition~\ref{riemp1} any
contraction of the Riemann tensor reduces to $\pm
R_{\alpha\beta}$, and $R_{\alpha\beta}=R_{\beta\alpha}.$
Similarly, the \emph{Ricci scalar}\index{Ricci scalar} is given by
\begin{equation}\label{eq:riccis}
  R=g^{\alpha\beta}R_{\alpha\beta}.
\end{equation}
\subsubsection{Einstein Tensor} The \emph{Einstein
tensor}\index{Einstein tensor}, which comprises the left hand side
of the Einstein equations $G=\kappa T$, is the simplest $(0,2)$
tensor constructed from the Riemann tensor
$R^{\mu}_{\hspace{4pt}\alpha\nu\beta}$ and the metric
$g_{\alpha\beta}$ such that $$\nabla G \equiv \mbox{div}\,G=0.$$
We shall see in section~\ref{sec:matter} that div$T=0$, so that
when $G=\kappa T$ we must have div$G=0.$ To derive the Einstein
tensor we consider the Riemann curvature components
$R^{\mu}_{\hspace{4pt}\alpha\nu\beta}$ in a Lorentzian coordinate
frame. In a Lorentzian frame the Riemann tensor, given in
equation~\eqref{eq:riemcurv}, can be written as
\begin{equation}\label{eq:riemcurvlor}
R^{\mu}_{\hspace{4pt}\alpha\nu\beta}
  =\frac{1}{2}g^{\mu\sigma}(g_{\sigma\beta,\alpha\nu}
  -g_{\sigma\nu,\alpha\beta}+g_{\alpha\nu,\sigma\beta}
  -g_{\alpha\beta,\sigma\nu}).
\end{equation}
Lowering the index $\mu$ in~\eqref{eq:riemcurv}, and taking the
partial derivative with respect to $x^\lambda$ yields
\begin{equation}\label{eq:riemcurvder}
R_{\mu\alpha\nu\beta,\lambda}
  =\frac{1}{2}(g_{\mu\beta,\alpha\nu\lambda}
  -g_{\mu\nu,\alpha\beta\lambda}+g_{\alpha\nu,\mu\beta\lambda}
  -g_{\alpha\beta,\mu\nu\lambda})
\end{equation}
Exploiting the symmetries of the metric
$g_{\alpha\beta}=g_{\beta\alpha}$, and the fact that partial
derivatives commute, equation~\eqref{eq:riemcurvder} can be
written as
\begin{equation}\label{eq:bianchi0}
R_{\mu\alpha\nu\beta,\lambda} +R_{\mu\alpha\lambda\nu,\beta}
+R_{\mu\alpha\beta\lambda,\nu}=0.
\end{equation}
Now, since our coordinate frame is locally Lorentzian we have
$\G^\sigma_{\alpha\beta}=0$, and so equation~\eqref{eq:bianchi0}
is equivalent to
\begin{equation}\label{eq:bianchi1}
R_{\mu\alpha\nu\beta;\lambda} +R_{\mu\alpha\lambda\nu;\beta}
+R_{\mu\alpha\beta\lambda;\nu}=0.
\end{equation}
The relations of the components given in
equation~\eqref{eq:bianchi1} are called the \emph{Bianchi
Identities}\index{Bianchi Identities}, and
since~\eqref{eq:bianchi1} is a tensor equation it is valid in any
coordinate system.

The covariant derivative is commutative with respect to
contraction, that is,
$$\nabla_\mu\left(
T^{\alpha_1\cdots\gamma
\cdots\alpha_k}_{\beta_1\cdots\gamma\cdots\beta_l} \right)
=\nabla_\mu\, T^{\alpha_1\cdots\gamma
\cdots\alpha_k}_{\beta_1\cdots\gamma\cdots\beta_l}
$$
for any tensor $T$. Also, by equation~\eqref{eq:met0},
$g_{\alpha\beta;\gamma}=0$, and since the inverse components
$g^{\alpha\beta}$ a functions of the metric components
$g_{\alpha\beta}$ it follows that
$$g^{\alpha\beta}_{\hspace{10pt};\gamma}=0.$$
Therefore, if we apply the Ricci contraction~\ref{eq:riccit} to
the Bianchi identities~\eqref{eq:bianchi1}, then
\begin{equation}\label{eq:bcont}
g^{\mu\nu}\left(R_{\mu\alpha\nu\beta;\lambda}
+R_{\mu\alpha\lambda\nu;\beta}
+R_{\mu\alpha\beta\lambda;\nu}\right)
=R_{\alpha\beta;\lambda}-R_{\alpha\lambda;\beta}
+R^{\nu}_{\alpha\beta\lambda;\nu}=0,
\end{equation}
where we have used the antisymmetric property of the Riemann
curvature tensor given in equation~\eqref{eq:riemp1}. Now, we
contract again on the indices $\alpha$ and $\beta$, yielding
\begin{equation}\label{eq:bcont2}
g^{\alpha\beta}\left(
R_{\alpha\beta;\lambda}-R_{\alpha\lambda;\beta}
+R^{\nu}_{\alpha\beta\lambda;\nu}\right)
=R_{;\lambda}-R^\beta_{\lambda;\beta} -R^{\nu}_{\lambda;\nu}=0,
\end{equation}
or
\begin{equation}\label{eq:bcont2a}
(2R^\nu_\lambda-\delta^\nu_\lambda R)_{;\nu}=0.
\end{equation}
Then we define the Einstein curvature tensor as
\begin{equation}\label{eq:einsteintens}
  G^{\alpha\beta}
  =R^{\alpha\beta}-\frac{1}{2}g^{\alpha\beta}R =G^{\beta\alpha},
\end{equation}
and by equation~\eqref{eq:bcont2a} we have
$$G^{\alpha\beta}_{\hspace{10pt};\beta}=0.$$

\section{Describing Matter in Spacetime}\label{sec:matter}
In this section we will derive the stress-energy tensor. It is
this tensor, $T$, which describes matter, that will be inexorably
linked to the geometry of spacetime via the Einstein field
equations $G=\kappa T$. We assume that physical objects in the
region of spacetime under consideration are regarded as a fluid,
that is a continuous distribution of matter. In particular, we
assume the distribution of matter is approximated by a perfect
fluid\index{perfect fluid} which is described as a fluid with no
heat conduction or viscosity. No heat conduction implies that
energy can flow only if particles can flow, and no viscosity means
that all forces are perpendicular to the interface between
particles~\cite{schutz}. There are other stress-energy tensors for
fluids not incorporating the assumptions of a perfect fluid, but
the tensor derived below will be that for a perfect fluid.

\subsection{The Stress-Energy Tensor}
We begin with the notion of the four-momentum. A curve passing
through an event is classified as timelike\index{timelike curve},
lightlike\index{lightlike curve}, or spacelike\index{spacelike
curve} according to whether the inner product of its tangent
vectors $g_{\alpha\beta}X^{\alpha}X^{\beta}$ is timelike,
lightlike, or spacelike, see definition~\ref{def:vec}. The path of
any material particle passing through an event $p$ in spacetime
must lie inside the lightcone of $p$, see figure~\ref{fig:causal},
otherwise an observer at $p$ would see the particle moving faster
than the speed of light. Thus any moving particle must lie on a
timelike curve. The timelike path of a moving particle may be
parameterized by the \emph{proper time}\index{proper time} $\tau$
which is defined by
\begin{equation}\label{eq:propertime}
  \tau=\int (-g_{\alpha\beta}X^{\alpha}X^{\beta})^{1/2}\,d\xi,
\end{equation}
where $\xi$ is any arbitrary parameterization of the path, and
$X^{\alpha}$ is the tangent vector to the curve in this
parameterization~\cite{wald}.

The \emph{four-velocity}\index{four-velocity} is defined by the
tangent vector $u^{\alpha}$ of a timelike curve parameterized by
proper time $\tau.$ In the presence of gravity any material
particle has an associated \emph{rest mass}\index{rest mass} $m$,
which shows up in the equations of motion~\cite{wald}. The
\emph{four-momentum}\index{four-momentum} vector, $p^{\alpha}$ of
a particle with rest mass $m$ is given by
\begin{equation}\label{eq:4moment}
 p^{\alpha}=m\frac{dx^\alpha}{d\tau}=m u^{\alpha}.
\end{equation}
The \emph{Energy}\index{energy} of a particle with four-velocity
$v^\beta$ is given by
\begin{equation}\label{eq:energy}
  E=-g_{\alpha\beta}p^{\alpha}v^{\beta},
\end{equation}
taken from~\cite{wald}. In a Lorentzian frame moving with the
particle, called a \emph{rest frame}\index{rest frame}, we have
$d\tau=dx^0$, and $u^0=1,$ and $u^a=0$ for $a=1,2,3$, and so
$E=-g_{00}p^0 u^0=p^0=m$. We note that this is the well known
equation $E=mc^2$ in our units with the speed of light $c=1$. This
We identify the energy of a particle with the 0th component of its
four-momentum. Furthermore, this also shows that, in general, we
cannot differentiate between mass and energy.

Now, we define the general stress-energy tensor, in terms of its
components, as
\begin{equation}\label{eq:strendef}
T^{\alpha\beta}\equiv\left\{
\begin{array}{l}
\mbox{The rate at which the $\alpha$th component of} \\
\mbox{momentum crosses a surface of constant}\; x^{\beta}.
\end{array}
\right\}.
\end{equation}
From this definition, along with our assumption of no heat
conduction, and no viscosity we can derive the components of the
stress-energy tensor in a Lorentzian frame moving with the same
velocity as the fluid.
\begin{itemize}
\item $T^{00}$ is the rate at which energy crosses the surface of
constant time, and is referred to as the mass/energy density or
simply the energy density, and is denoted by Greek letter $\rho$.

\item $T^{a0}=T^{0a}=0$ is a consequence of no heat conduction.
The rate of momentum in the direction $x^a$ across a surface
constant in time is zero because without heat flow the momentum of
a particle cannot change without moving, which means $T^{a0}=0$.
Since the coordinate frame is moving with zero velocity with
respect to the fluid, with no heat conduction the energy is
constant across a surface of constant $x^a$ which means
$T^{0a}=0$.

\item $T^{ab}$ is the rate of the $a$th momentum component across
a surface of constant $x^a$. The absence of viscosity, a force
parallel to the interface between particles, implies that the
force should be perpendicular to the interface. Consequently, a
surface of constant $x^a$ will only be the force per unit area
against it in the $x^a$ direction which is equal for each $a$.
Mathematically, this means $T^{ab}=p\delta^{ab}$ where $p$ denotes
the pressure.
\end{itemize}
Therefore, in a Lorentzian frame\index{Lorentzian frame}
comoving\index{comoving} with a perfect fluid the stress-tensor is
given by
\begin{equation}\label{eq:stren}
T^{\alpha\beta}=pg^{\alpha\beta}+(p+\rho)u^\alpha u^\beta,
\end{equation}
and written in covariant form
\begin{equation}\label{eq:strencov}
T=(\rho+p)u\otimes u+pg.
\end{equation}

\subsubsection{Conservation of Energy and Momentum}
The stress-energy tensor $T$ is a representation of the energy and
momentum in a fluid, therefore it seems reasonable to expect there
to be a way to express the conservation of energy and momentum
using $T$. Consider a cubical fluid element with each side of
length $l$. Then amount of energy-momentum coming in must equal
the amount of energy-momentum going out of the cube. That is,
\begin{equation}\label{eq:flow}
\frac{\partial}{\partial x^0}l^3 T^{\alpha 0} =\sum_{\beta=1}^3
l^2\left\{T^{\alpha\beta}(x^\beta=0)
-T^{\alpha\beta}(x^\beta=l)\right\},
\end{equation}
where $l^2T^{\alpha\sigma}(x^\sigma=0)$ is the rate of flow in the
$x^\sigma$ direction across the surface of constant $x^\alpha$ at
$x^\sigma=0$, and $-l^2T^{\alpha\sigma}(x^\sigma=l)$ the rate at
$x^\sigma=l$ with the minus sign indicating that the flow is
coming in from the opposite direction of that at $x^\sigma=0$. Now
we divide each side of equation~\eqref{eq:flow} by $l^3$, and take
the limit as $l\to 0$ to find
\begin{equation}\label{eq:conslaw}
T^{\alpha\beta}_{\hspace{10pt},\beta}=0,
\end{equation}
where we have used the definition of the derivative
\begin{equation}\label{eq:deriv}
\lim_{l\to 0}\frac{T^{\alpha\beta}(x^\beta=0)
-T^{\alpha\beta}(x^\beta=l)}{l}=-\frac{\partial}{\partial x^\beta}
T^{\alpha\beta}.
\end{equation}
In our Lorentzian frame,
\begin{equation}\label{eq:covconslaw0}
T^{\alpha\beta}_{\hspace{10pt};\beta}=T^{\alpha\beta}_{\hspace{10pt},\beta}=0,
\end{equation}
and since this is a tensor equation we can say
\begin{equation}\label{eq:covconslaw}
\nabla_\beta T^{\alpha\beta}=0
\end{equation}
is the covariant expression for conservation of energy and
momentum.

\section{The Einstein Equations}
The idea behind the Einstein equations\index{Einstein equations}
is that the sources of the gravitational field determine the
metric\cite{schutz}. Classically, this is given by Poisson's
equation
\begin{equation}\label{eq:newton}
-\Delta \Phi=4\pi \mathcal{G} \rho,
\end{equation}
whose solution is given
$$\Phi(\x)=\int_{\mathbb{R}^3}
\frac{\mathcal{G}\rho(\mathbf{y})}{|\x-\mathbf{y}|}d\mathbf{y},$$
where $\mathcal{G}$ is Newton's gravitational constant. The
function $\Phi$ is the Newtonian gravitational source due the mass
density $\rho.$ This leads to the question, what quantity acts as
the source for the gravitational field? Since the stress-energy
tensor $T$ contains the mass-energy density $\rho$, the component
$T^{00},$ and is covariant we postulate $T$ as the source for the
gravitational field. Then to derive the relativistic
generalization of equation~\eqref{eq:newton}, $T$ should be
coupled to a symmetric $(0,2)$ tensor whose covariant derivative
vanishes, and created from the geometry of spacetime.  Finally,
the coupling should reduce to the Newtonian equivalent given in
equation~\eqref{eq:newton} when velocities and the gravitational
field are sufficiently small. The Einstein tensor\index{Einstein
tensor} $G$ was specifically designed to satisfy these criteria.
Therefore,
\begin{equation}\label{einstein}
G=\kappa T,
\end{equation}
with $\kappa=8\pi \mathcal{G}/c^4,$ which reduces to
equation~\eqref{eq:newton} in the Newtonian limit of low
velocities, and weak gravitational fields. For a perfect fluid the
Einstein equations are written as
\begin{equation}\label{eq:einsteinperf}
G^{\alpha\beta}= R^{\alpha\beta}-\frac{1}{2}g^{\alpha\beta}R =8\pi
\mathcal{G} \left(pg^{\alpha\beta}+(p+\rho)u^\alpha u^\beta
\right) = 8\pi \mathcal{G} T^{\alpha\beta}.
\end{equation}

\section{The Second Fundamental Form}
The shock-wave solutions of the Einstein equations\index{Einstein
equations} considered in chapters~\ref{chp:3} and~\ref{chp:4} are
hypersurfaces of the spacetime manifold. A
\emph{hypersurface}\index{hypersurface} is a $n-1$ dimensional
submanifold of an $n$ dimensional manifold. The \emph{second
fundamental form}\index{second fundamental form} describes how the
hypersurface is embedded in the spacetime manifold $M$ by
recording how the tangent spaces of the hypersurface change over
the hypersurface. The \emph{first fundamental form}\index{first
fundamental form} is the metric tensor $g$. In this section we
will derive the second fundamental form.

Let $\Si$ be a hypersurface of $M$, and denote the induced metric
on $\Si$ by $\tilde{g}$, which is assumed to be non-degenerate.
For every point $p$ in $\Si$, there exist \emph{slice
coordinates}\index{slice coordinates} $x=(x^1,\ldots,x^{n})$ on a
neighborhood $U$ of $p$ in $M$ such that $U \cap \Si$ is given by
$x^n=0$, and $(x^1,\ldots,x^{n-1})$ form local coordinates for
$\Si$~\cite{lee}.

Since the second fundamental form\index{second fundamental form}
is a measure of how a tangent space changes on a
hypersurface\index{hypersurface} as is moves from point to point,
it should involve the the connection\index{connection} $\nabla$ on
the surface. Therefore, we need the idea of an \emph{induced
connection}\index{connection!induced} on the hypersurface, and how
to describe vector fields on the surface.

Given vector fields $X$ and $Y$ on $M$, we already have already
have a way to measure the vector rate of change of $Y$ in the
$X_p$ direction by the vector field $\nabla_X Y$ on $M$. We should
also note that the connection $\nabla$ is unique, since one can
show that the condition, equation~\eqref{eq:met0}, $\nabla_\alpha
g_{\beta\gamma}=0$ implies that $\nabla$ is unique~\cite[theorem
3.1.1]{wald}. We want to somehow induce the differential structure
on to our hypersurface.

Any smooth vector vector field $\tilde{X}$ on a submanifold always
has a \emph{local extension}\index{local extension} into its
ambient manifold~\cite{warner}. In our case, where $\Si$ is a
hypersurface of $M$, this means given a smooth vector field
$\tilde{X}$ on $\Si$, for each point $p\in\Si$, there exists a
neighborhood $U\subset\Si$ of $p$ and a neighborhood $V\subset M$
of $p$ such that $U\subset V$, and there exists a smooth vector
field $X$ on $V$ such that
$$X|_U=\tilde{X}|_U.$$
For each $p$ in $\Si$ the metric $g$ splits $T_p M$ into the
direct sum
$$T_p M =\tans\oplus (\tans)^{\perp}.$$
Therefore, if $X\in\tanm$, then $$X=X^{tan}+X^{nor},$$ where
$X^{tan}\in \tans$ and $X^{nor}\in(\tans)^{\perp}$.
\subsubsection{The Induced Connection}\index{connection!induced}
Suppose $\Si$ is a hypersurface\index{hypersurface} of the
spacetime manifold $M$, and $\widetilde{X}, \widetilde{Y}$ are
vector fields on $\Si$. If $X,Y$ are local extensions of
$\widetilde{X}, \widetilde{Y}$ to $M,$ then we define
\begin{equation}\label{eq:indcon}
  \widetilde{\nabla}_{\widetilde{X}} \wt{Y}=\left(\nabla_X Y\right)^{tan},
\end{equation}
where $\wt{\nabla}$ denotes the connection on $\Si$.

\begin{definition}
If $\wt{X},\wt{Y}$ are local vector fields on $\Si$, then
$$I\!I(\wt{X},\wt{Y})=\nabla_X Y-\wt{\nabla}_{\wt{X}}\wt{Y}
=(\nabla_X Y)_\perp,$$ is said to be the \emph{second fundamental
form tensor}~\cite{oneill}\index{second fundamental form!tensor}.
\end{definition}
The tensor $I\!I$ is a local vector field on $M$ normal to $\Si$

\begin{definition}
Let $\n$ be a unit normal vector field on a hypersurface $\Si$ of
$M$. The $(1,1)$ tensor field $K$ on $\Si$ such that
$$\langle K(X), Y\rangle=\langle I\!I(X,Y),\n\rangle$$
for all vector fields $X,Y$ on $\Si$ is called the \emph{second
fundamental form}\index{second fundamental form!definition}. $K$
is also sometimes called the \emph{shape operator}\index{shape
operator}~\cite{oneill}.
\end{definition}
Here, $K$ determines a linear operator $K:T_p\Si \to T_p\Si$ at
each point $p\in \Si$~\cite{oneill}

\begin{lemma}
The second fundamental form,
\begin{equation}\label{eq:2form}
K=-\nabla_{X} \n
\end{equation}
at each point $p\in\Si$, and
$K:T_p\Si \to T_p\Si.$
\end{lemma}

\begin{proof}
Since $\langle \n,\n\rangle$ is constant, $\langle\;\ ,\;\rangle$
satisfies the Leibnitz rule, and is symmetric we have that
$\langle \nabla_X \n,\n\rangle=0$ for $X\in T_p\Si.$ Thus,
$\nabla_X \n$ is tangent to $\Si$ for all $X\in T_p\Si.$ Now, for
all vector fields $Y$ on $\Si$,
$$\langle K(X), Y\rangle=\langle I\!I(X,Y),\n\rangle
=\langle\nabla_X Y -(\nabla_X Y)^{tan},\n\rangle =\langle\nabla_X
Y,\n\rangle,$$ and since
$$\langle\nabla_X
Y,\n\rangle+\langle\nabla_X\n,Y\rangle=\nabla_X \langle Y, \n
\rangle  =0.,$$ we have that
$$\langle K(X), Y\rangle=-\langle\nabla_X\n,Y\rangle.$$
\end{proof}

% ---End Background Chapter---------------------------------------

%% file: main.tex
% ---Main---------------------------------------------------------
%Shock-Wave Solutions of Einstein Equations in the Lightlike Limit
\label{chp:3}

%\section{introduction}
The central results of this dissertation are contained in this
chapter. The main result, Theorem~\ref{thm:main}, is a set of
equivalent conditions that imply conservation of energy across a
surface in which the metric is only Lipschitz continuous. The
other pertinent result, Theorem~\ref{thm:sphere}, is a set of
conditions which are equivalent to conservation of energy across a
spherically symmetric hypersurface. In the case when the
hypersurface under consideration is non-null, both
Theorems~\ref{thm:main} and~\ref{thm:sphere} were first proved by
Smoller and Temple in~\cite{st94}. The same results can also be
found in~\cite{stg01}, and~\cite{stgerm99}. Our contribution here
will be to generalize Smoller and Temple's results to include the
lightlike (or null) case.

There are two main difficulties in the lightlike case. The first
is that the metric restricted to the lightlike surface is
degenerate. We will deal with this issue by considering the
problem in the context of the whole spacetime manifold where the
metric is not degenerate. The other, more problematic difficulty,
is that the second fundamental form, defined by $K=-\nabla_X \n$
where $X$ is a vector tangent to the surface and $\n$ is the
normal to the surface, cannot be used to describe the dynamics of
the surface in the ambient spacetime. This is because $K$ measures
the change in the normal vector $\n$ in the direction tangent to
the surface, but in a lightlike surface\index{lightlike
hypersurface} $\n$ lies in the tangent space of the surface; thus
$K$ can no longer give geometrical information about how the
surface relates to the ambient manifold. We rectify this, using an
idea of Barrab\`{e}s and Israel given in~\cite{isr91}, by defining
a generalized second fundamental form\index{second fundamental
form!generalized} $\K=-\nabla_X \N$ where $\N$ is a vector
transverse to the surface.

In the first section we define the generalized the second
fundamental form $\K$, and construct a coordinate system
analogous to Gaussian normal coordinates so that in this
coordinate system the components of $\K$ can be written in terms
of the derivative of the corresponding metric components in the
direction of $\N$. In the second section we state the main
theorem of this dissertation. Then in the following section we
give the supporting lemmas and theorems which will yield the
proof of the main theorem. We finish off the chapter by stating
and proving the theorem involving matching spherically symmetric
metrics across a hypersurface.

\section{Generalizing the Second Fundamental Form}
Let $M$ denote a manifold equipped with a metric $g$ with fixed
Lorentzian signature $(-+\cdots+)$, and $\Si$ denote a
hypersurface\index{hypersurface}, which may or may not be null,
that divides $M$ into two regions $M^L$ and $M^R.$ Locally, we
define $\Si$ by $\psi(y)=0$, where $\psi$ is a smooth function
such that
  \begin{equation}\label{eq:norm}
    n_i\,dy^i=\frac{\partial \psi}{\partial y^i} dy^i \neq 0
  \end{equation}
for any coordinate system. Let $g^L$ and $g^R$ denote the metrics
on $M^L$ and $M^R$ respectively, and the metric on $M$ is given
by $g=g^L\cup g^R$.  Assume $g^L,$ and $g^R$ are smooth, that is
at least $C^2$, on each side of $\Si$ with their derivatives
uniformly bounded across $\Si.$

\subsection{The failure of the Second Fundamental Form in the Lightlike
Case} Given a point, $p$ in $\Si$ let $\tans$ denote the tangent
space of $\Si$ at $p$. A vector $X$ is in $\tans$ if
\begin{equation}\label{eq:tang}
\langle\n, X\rangle=g_{\alpha\beta}n^{\alpha} X^{\beta}=0.
\end{equation}
If $\Si$ is a lightlike hypersurface\index{lightlike
hypersurface}, then $\langle\n,\n\rangle=0$, which implies that
$\n$ is also in $\tans.$ An unfortunate consequence of this is
that $[K]=0$ is always true, where $[ \, \cdot \, ]\equiv
\K^L-\K^R$ denotes the jump in $\K$ across $\Si$, and hence does
not yield any information about how $\Si$ is embedded in $M$.
Therefore, we cannot use the standard second fundamental form to
determine the necessary and sufficient conditions for a null
hypersurface to be a shock surface\index{shock surface} or a
surface layer\index{surface layer}. This failure of the second
fundamental form to yield conditions on the existence of lightlike
shocks is a consequence of $\n$ not being transverse to $\Si$;
therefore any alternative to $K$, designed to extract this
information, should involve a vector transverse to $\Si$.

\subsection{Generalized Second Fundamental Form}
We will generalize the second fundamental\index{second fundamental
form!generalized} form given by equation~\eqref{eq:2form} so that
it will still be defined in the case of lightlike hypersurfaces.
We will use an idea introduced by Barrab\`{e}s and Israel
in~\cite{isr91}. Even though they use a scalar version of the
second fundamental form and we are not, the idea is still the
same.

We begin by choosing a vector $\N$ transverse to $\Si$, that is,
a vector not in $\tans.$ Then we replace the second fundamental
form $K(X)=-\nabla_X \n$ with what we will call the
\emph{generalized second fundamental form}\index{second
fundamental form!generalized}
\begin{equation}\label{eq:t2form}
\K(X)=-\nabla_X \N,
\end{equation}
where $X$ is in $\tans.$ If we compare this with the definition
of the second fundamental form in equation~\eqref{eq:2form}, then
we notice that $\K$ depends on the vector $\N$ that we choose
just as $K$ depends on the normal vector $\n$. To make $\K$
well-defined we need to place certain restrictions on which
transverse vectors $\N$ can be chosen to define $\K$.

We need $\N$ to be continuous across $\Si$, so we require for all
points $p$ in $\Si$
\begin{equation}\label{eq:njump}
 [\langle\N,X_a\rangle]=[N_a]=0,
\end{equation}
where $\{X_a\}_{a=1}^{n-1}$ is a basis of $\tans.$ Furthermore,
the length of $\N$ on each side of $\Si$ must coincide; hence
\begin{equation} \label{eq:njump1}
 \left[\langle\N,\N\rangle\right]=0.
\end{equation}
The transverse vector chosen under the conditions~\eqref{eq:njump}
and~\eqref{eq:njump1} is not unique. Any vector $\N$ satisfying
equations~\eqref{eq:njump} and~\eqref{eq:njump1} is invariant
under the transformation
\begin{equation}\label{eq:Ntrans}
  \N \to \N'=\N+\lambda^a \,X_a,
\end{equation}
where, again, $\{X_a\}_{a=1}^{n-1}$ is a basis of $\tans,$ and
$\lambda^a$ is an arbitrary function. Thus, we are free to choose
$\N$ so that
\begin{equation}\label{eq:Nneta}
\langle\N,\n \rangle=\eta\neq 0,
\end{equation}
everywhere on $\Si$, where $\eta$ is any non-zero constant we
like. Now, under the transformation~\eqref{eq:Ntrans}, $\K$
transforms as
\begin{equation}\label{eq:Ktrans}
  \K_{a}^{b}X^a \To
\K_{a}^{b}X^a-\lambda^{c}\, \G^{b}_{ac}X^a.
\end{equation}
Equation~\eqref{eq:Ktrans} comes from the computation
$$ -\left(\nabla_X \left(\N+\lambda^l X_l\right) \right)^b
=\K_{a}^{b}X^a-(\lambda^{b}_{\hspace{4pt},a}+\lambda^{c}\,
\G^{b}_{ac})X^a =\K_{a}^{b}X^a-\lambda^{c}\, \G^{b}_{ac}X^a.$$
Equation~\eqref{eq:Ktrans} tells us that $\K$ does not transform
as a tensor, since $\lambda^{k}\, \G^{j}_{ik}$ does not. Since we
are seeking a tool to describe how $\Si$ in embedded in ambient
spacetime that is coordinate independent, this presents a slight
difficulty. This is rectified by the insight of Barrab\`{e}s and
Israel in~\cite{isr91} that $[\K]$ is a tensor. In the lemma below
we will show that $[\K]$ is invariant under the
transformation~\eqref{eq:Ntrans}, thus $[\K]$ does not depend on
the particular transverse vector $\N$ that is chosen other than
that it must satisfy the the jump conditions~\eqref{eq:njump}
and~\eqref{eq:njump1}. The reason for the invariance is that the
metric $g$ and its tangential derivatives are always continuous
across $\Si$, but its transversal derivatives may not be.

This means we can now say that $\K$ is a true generalized form of
the standard second fundamental form $K$ when we compare their
jumps across $\Si$. That is, if $\Si$ is non-lightlike, then we
choose $\n$, or another suitable vector, as our transverse vector
we get $[K]=[\K]$, and in the limit, as $\Si$ goes lightlike,
$[\K]$ remains well-defined; hence $[\K]$ is the tensor we seek
to determine how $\Si$ is embedded in spacetime even when $\Si$
is lightlike.

\begin{lemma}\label{l2}
The jump $[\K]=(\K^R-\K^R)$ across $\Si$ is independent of the
choice of $\N,$ and is a tensor.
\end{lemma}
\begin{proof}
Lowering the index $j$ in equation (\ref{eq:Ktrans}) we get
\begin{equation}\label{eq:Ktransl}
  \K_{ab}X^a\To \K_{ab}X^a-\lambda^{c}g_{db}\G^{d}_{ac}X^a.
\end{equation}
Since
$$g_{db}\G^{d}_{ac}=\frac{1}{2}(-g_{ac,b}+g_{ba,c}+g_{ca,b})$$
only involves derivatives in the tangential direction, and
\begin{equation}\label{}
  [\K_{ab}-\lambda^{c}\, g_{lb}\G^{l}_{ac}]X^a
  =[\K_{ab}]X^a-\lambda^{c}\,[ g_{lb}\G^{l}_{ac}]X^a,
\end{equation}
under the transformation (\ref{eq:Ntrans}), we have
\begin{equation}\label{eq:jumptrans}
  [\K_{ab}]X^a\To [\K_{ab}]X^a-\lambda^{c}\,[ g_{db}\G^{d}_{ac}]X^a
  =[\K_{ab}]X^a
\end{equation}
\end{proof}

\subsection{Modified Gaussian Skew (MGS) Coordinates}

We construct a new set of coordinates which will make computation
of the Riemann and Einstein curvature tensors more manageable.

One important ingredient in characterizing solutions of the
Einstein equations across an interface, first formulated by
Israel~\cite{isr66}, and also used in Smoller and Temple's
work~\cite{st94,st95,stgerm99}, is that in Gaussian normal
coordinates you can write the second fundamental form in terms of
the metric as
$$K_{ab}=-\frac{1}{2}g_{ab,0},$$
where $0$ represents the $0$th coordinate, whose corresponding
basis vector is the normal to $\Si.$ In Gaussian normal
coordinates, values of $[K]$ can be used to easily compute the
corresponding jumps in the Riemann, Ricci, and Einstein tensors.
In the lightlike case, however, Gaussian normal coordinates are
undefined since the normal vector $\n$ is a tangent vector. To
rectify this we employ Modified Gaussian Skew (MGS)\index{MGS
coordinates} coordinates\index{Modified Gaussian Skew
coordinates}.

A MGS coordinate system is constructed in an manner analogous to
Gaussian normal coordinates by replacing the vector normal to
$\Si$ with a transverse vector that satisfy the jump
conditions~\eqref{eq:njump} and~\eqref{eq:njump1}. The property
that the $0$th coordinate vector is orthogonal to the other $n-1$
coordinate vectors is lost in a MGS coordinate system, but such a
coordinate system can be defined even when $\Si$ is a null
hypersurface. However, to deal with the degeneracy of the vector
subspace $\tans$ we will place certain restrictions on the basis
of $\tans$ we choose. These additional restrictions will define
MGS coordinates, and in doing so, the metric, in MGS coordinates,
will be almost diagonal.

We model the construction of MGS coordinates after the
construction of Gaussian normal coordinates by Smoller and Temple
from~\cite{st94, stgerm99}. One can construct MGS coordinate as
follows:

\begin{itemize}
  \item We first choose coordinates in the surface. Let $\tans$ denote
  the tangent space of $\Si$ at point $p\in\Si.$
  The surface coordinates $y^1,\ldots,y^{n-1}$ will be chosen to coincide with
  a particular set of basis vectors of $\tans$. First,
  choose $y^1$ so that $\frac{\partial \,}{\partial y^1} =\n$.
  The remaining $y^i$, $i=2,\ldots,n-1$, are chosen so that
  $g_{ij}=\left\langle \frac{\partial}{\partial y^i},
  \frac{\partial}{\partial y^j}\right\rangle=0$ for $i\neq j.$
  We will show below that $g_{ii}>0$. Furthermore, we will also place the
  restriction on each $y^i$ that
  \begin{equation}\label{eq:Nperp}
  \left\langle \frac{\partial}{\partial y^i},
  \N\right\rangle=0,
  \end{equation}
  for $i=2,\ldots,n-1$.

  Notice that we can always choose such a coordinate system. First,
  $\n$ is in $\tans,$ and is non-zero, so we can
  choose an orthogonal basis incorporating $\n$. Also recall that once
  $\N$ is chosen we are free to modify it via a transformation of the
  form given in equation~\eqref{eq:Ntrans}. Thus we can always transform
  \begin{equation}\label{eq:Nmgstrans}
  \N\to \N-\sum_{i=2}^{n-1}\left\langle \frac{\partial}{\partial y^i},
  \N\right\rangle \frac{\partial}{\partial y^i}\end{equation}
  in order to satisfy
  equation~\eqref{eq:Nperp}.
  \item For each $p\in\Si$, let $\gamma_p(s)$ denote the
  geodesic which satisfies
  $$\gamma_p(0)=p, \qquad \dot{\gamma}_p(0)=\N,$$
  where $\N$ is a vector transverse to $\Si$ satisfying the
  jump conditions~\eqref{eq:njump} and~\eqref{eq:njump1}.
  \item Assume $\N$ points into the right hand side of $\Si$, for
  convenience.
  \item Choose coordinate $w^0$ so that if $\gamma_p(s)=q$, then $w^0(q)=s.$
  As a consequence, $w^0<0$ on the left hand side of $\Si$, and $w^0>0$
  on the right hand side of $\Si.$
  \item The coordinates $w^1,\ldots,w^{n-1}$ are given by $w^i(p)=y^i(p)$ of
  $p\in\Si,$ and $w^i(q)=w^i(p)$ if and only if $q=\gamma_p(s),$
  for some $p$ and $s,$ where $i=1,\ldots,n-1$.
\end{itemize}
By construction, we can justify the following lemma.
\begin{lemma} In MGS coordinates,
\begin{equation}\label{mgsmetric}
ds^2=g_{00}d(w^0)^2+2\eta\,dw^0 dw^1+g_{ii}dw^i dw^i,
\end{equation}
where $i=2,\ldots,n-1,$ and $\eta=g(\N,\n).$
\end{lemma}
We now show how an MGS coordinate system relates at other
coordinate system.
\begin{lemma}\label{l:C11trans}
There exists a $C^{1,1}$ transformation from any coordinate system
to an MGS coordinate system.
\end{lemma}
\begin{proof}
The hypersurface $\Si$ is given by $\psi(y)=0$ in $y$-coordinates.
Using a smooth transformation we can use \emph{slice
coordinates}\index{slice coordinates} to define $\Si$ by $y^0=0$.
We assume that $\Si$ is lightlike, hence we can write the normal
to $\Si$ in terms of vectors in $\tans$ for any $p$ in $\Si,$ so
that
$$\n=n^a\frac{\partial}{\partial y^a},\quad a\neq0.$$
Let $\mu_p(t)$ denote the geodesic which satisfies
$$\mu_p(0)=p,\quad\dot{\mu}_p(0)=\n.$$
Then we can smoothly transform the slice coordinates $y$ to a new
set of slice coordinates $u$, where $u^1=t$ and
$u^\alpha=y^\alpha$ for $\alpha\neq 1$, so that $\Si$ is given by
$u^0=0.$ In $u$-coordinates we are free to use the
transformation~\eqref{eq:Nmgstrans} to adjust $N$ so that
$$\left\langle \frac{\partial}{\partial u^i},
\N\right\rangle=0,$$ but still have $\langle \N,\n \rangle=\eta$.
Then our MGS coordinates can be written as
$w=(w^0,\ldots,w^{n-1})=(s,t,u^2,\ldots,u^{n-1}).$ Now for $q$ not
in $\Si$, but on the geodesic $\gamma_p(s)$ and in the coordinate
neighborhood, we have that $w(q)=(s,0,u^1(p),\ldots,u^{n-1}(p))$
is a smooth function on each side of $\Si$. Therefore, it only
remain to check the continuity of the derivatives across at $s=0$.
Indeed,
$$\frac{\partial u^\alpha}{\partial w^b}=\delta^\alpha_b,$$
and
$$\frac{\partial u^\alpha}{\partial w^0}=N^\alpha,$$
where $N^a$ denotes the components of $\N$ in $u$-coordinates.
This comes from the vector transformation law:
$$\N=\frac{\partial}{\partial w^0}
=\frac{\partial u^\alpha}{\partial w^0} \frac{\partial}{\partial
u^\alpha} =N^\alpha\frac{\partial}{\partial u^\alpha}.$$
Since
$\N$ is continuous across $\Si$ it follows that the derivatives
are as well.
\end{proof}

\begin{remark}\label{g00}
Notice that it is possible for $\N$ to be null, which would mean
$g_{00}=0$ in MGS coordinates. We will see in the proof of the
lemma below that the span of the vectors $\N$ and $\n$ forms a
two dimensional Lorentz vector space. Since any two dimensional
or greater Lorentz vector space contains two linearly independent
null vectors, it is possible for both $\N$ and $\n$ to be null.
For more information regarding Lorentz vector spaces
see~\cite[Chapter 5]{oneill}.
\end{remark}

\begin{lemma}\label{posdef}
In MGS coordinates,the metric components $g_{ii}>0$ for
$i=2,\ldots,n-1,$ and hence the coordinate basis vectors can be
normalized so that $g_{ii}=1$ for $i=2,\ldots,n-1.$
\end{lemma}
\begin{proof}
Let $\tanm$ denote the $n$-dimensional tangent space of the
spacetime manifold $M$ at a point $p$, and let $p \in \Si.$ Now
choose MGS coordinates on a neighborhood $U$ containing $p$. Then
the coordinate basis of $\tanm$ is $(\N,\n,e_2,\ldots,e_{n-1}),$
where $e_2,\ldots,e_{n-1}$ are orthonormal with respect to each
other, and are orthogonal to both $\N$ and $\n.$ Thus, in MGS
coordinates at $p$ the metric $g$ is given by the matrix
\begin{equation}\label{eq:gmgsmatrix}
\left(
\begin{array}{ccrcl}
1 & \eta & 0 \quad &\cdots & \quad 0 \\
\eta & 0 & 0 \quad &\cdots & \quad 0 \\
0  & 0 & g_{22} &  &\mbox{\Huge 0}  \\
\vdots  &\vdots &    &\ddots &  \\
0    &0&  \mbox{\Huge 0} & & g_{n-1\,n-1}
\end{array}\right),
\end{equation}
where $\eta$ is a constant, and $g_{ii}=\pm\, 1$ for
$i=2,\ldots,n-1$. Notice that for a transverse vector satisfying
the jump conditions~\eqref{eq:njump} and~\eqref{eq:njump1}, it is
possible for $g(\N,\N)=0,\pm 1$, which would make $g_{00}=0,\pm 1$
in MGS coordinates. However, we can always ``adjust" $\N$ by a
transformation of the form~\eqref{eq:Ntrans} to make $g_{00}=1$.
If $g(\N,\N)=0$, then let $\N \to \N+\n/2\eta$, and if
$g(\N,\N)=-1$, then let $\N \to \N+\n/\eta$. Therefore, we can say
$g_{00}=1$ in any MGS coordinate system.

Now, since $\tans$ is a degenerate subspace of the Lorentzian
vector space $\tanm$, and $g_{11}=g(\n,\n)=0$, we must have
$g_{ii}> 0$ because a degenerate subspace of a Lorentzian vector
space cannot contain any timelike vectors, and also cannot
contain more than one linearly independent lightlike
vector\footnote{See~\cite[Chapter 5]{oneill} for more information
on vector subspaces of Lorentzian vector spaces.}.
\end{proof}

\begin{remark}\label{rmk:lorentzcoord}
Using the idea from the above proof we can construct a coordinate
system so that the metric $g$ is of the form
$\mbox{diag}(-1,1,\dots,1)$, which we will use later in this
section as a candidate for an inertial coordinate frame. We
construct these coordinates, denoted by $u=(u^0,\ldots,u^{n-1})$,
on the neighborhood $U$ of $p$ in $\Si$ from MGS coordinates
$w=(w^0,\ldots,w^{n-1})$ on $U$ of $p$.  From Lemma~\ref{posdef}
we can choose $w^i$ so that $g(w^i,w^i)=1$ for $i=2,\ldots,n-1.$
Now we let $u^1=w^0$,$u^i=w^i$ for $i=2,\ldots,n-1,$ and choose
$u^0$ so that
$$\frac{\partial}{\partial u^0}
=\left(\frac{1}{\eta}\,\frac{\partial}{\partial
w^1}-\frac{\partial}{\partial w^0}\right) =(\n-\eta\N)/\eta.$$
Then, in $u$-coordinates,
$g_{00}=g\left((\n-\eta\N)/\eta,(\n-\eta\N)/\eta\right)=-1$,
$g_{11}=g(\N,\N)=1$, and
$g_{01}=g\left((\n-\eta\N)/\eta,\N\right)=0$. Thus, in
$u$-coordinates, the metric $g$ has the form
$\mbox{diag}(-1,1,\dots,1)$.
\end{remark}

\subsection{The Main Result of the Generalized Second Fundamental Form}
The following lemma shows that the components of the generalized
second fundamental form\index{second fundamental
form!generalized} can be written in terms of the derivative of
the corresponding component of the metric $g$ in the $\N$
direction. As a consequence, from values of the jump in the
transverse second fundamental form, $[\K]$, we will be able to
determine the values of the jumps of the Riemann curvature tensor
$[R^\mu_{\alpha\nu\beta}]$, and hence the Ricci tensor
$[R_{\alpha\beta}]$, the Ricci scalar $[R]$, and Einstein
curvature tensor $[G_{\alpha\beta}]$.

\begin{lemma}\label{l3}
In MGS coordinates
\begin{equation}\label{eq:gaussform}
  \K_{ab}=-\frac{1}{2}g_{ab,0},
\end{equation}
where $a,b=1,\ldots,n-1.$
\end{lemma}

\begin{proof}
For any vector field $X^{\beta}$ on the surface $\Si$ we have
\begin{equation}\label{eq:31l}
\begin{aligned}
-\K_{\beta}^{\tau} X^{\beta} & =(\nabla_{X}\N)_\beta^{\tau} X^{\beta}\\
& =\left(N^{\tau}_{\;,\beta}+N^{\sigma}\,
\G^{\tau}_{\beta\sigma}\right)X^{\beta}\\
& =\G^{\tau}_{\beta 0}X^{\beta},
\end{aligned}
\end{equation}
where $\beta,\tau, \sigma=0,\ldots,n-1$. The last equality comes
from the property that in MGS coordinates the components of the
transverse vector $\N$ are $N^a=0$, $a=1,\dots,n-1$, and $N^0=1$.

Since $g_{\sigma 0}$ is constant for all $\sigma=0,\ldots,n-1$,
the Christoffel symbol\index{Christoffel symbol} in equation
(\ref{eq:31l}) can be written as
\begin{equation}\label{eq:32l}
\G^{\tau}_{\beta 0}
=\frac{1}{2}\left\{g^{\tau\alpha}\set{-g_{\beta 0,\alpha} +
g_{\alpha\beta,0}+g_{0\alpha,\beta}}\right\} =\frac{1}{2}g^{\tau
a}g_{ab,0}.
\end{equation}
where $a,b=1,\ldots,n-1.$ Combining equations (\ref{eq:31l}) and
(\ref{eq:32l}) gives
\begin{equation}\label{eq:33l}
  -\K_{b}^{\tau} X^{b}=-\frac{1}{2}g^{\tau a}g_{a b,0} X^{b}.
\end{equation}
Lowering the index $\tau$ in equation (\ref{eq:33l}), where
\begin{equation}\label{eq:34l}
  g_{\sigma\tau}\,g^{\tau a}g_{a b,0}
  =\delta^{a}_{\sigma}g_{a b,0}=g_{a b,0},
\end{equation}
we get the desired result
\begin{equation}\label{eq:36l)}
  \K_{ab}=-\frac{1}{2}g_{ab,0}.
\end{equation}
\end{proof}
This result will be the key tool we will use to prove the Main
theorem of this chapter, and of this dissertation! We end this
section with a corollary of the above lemma which will we use to
prove part of our main theorem.

\begin{corollary}\label{c32}
The metric components of $g=g^L\cup g^R$ in MGS coordinates are
$C^{1}$ functions of the coordinate variables if and only if
$[\K]=(\K^R-\K^L)=0$ at each point on the surface $\Sigma.$
\end{corollary}

\begin{proof}
The jumps in derivatives of metric components in the direction of
the surface $\Si$ are always zero, that is, $g_{\alpha\beta,c}$
for $\alpha,\beta=0,\ldots,n-1$ and $c=1,\ldots,n-1$, since
$[g_{\alpha\beta}]=0$ on $\Si$. Also, $g_{0\beta}$ is constant for
each $\beta=0,\ldots,n-1$ Then by Lemma~\ref{l3},
$[\K_{ab}]=-\frac{1}{2}[g_{ab,0}]=0$ implies that $g$ is $C^1$ for
each coordinate.
\end{proof}

\section{The Main Theorem} The chief result of this
dissertation is the theorem stated below which classifies all
Lipschitz continuous solutions of the Einstein equations across a
hypersurface. Smoller and Temple first proved this result in the
case of non-lightlike hypersurfaces in~\cite{st94}, and can also
be found in~\cite{stg01} and~\cite{stgerm99}. We extend their
results here to include the case of a lightlike hypersurface using
the generalized second fundamental form, developed in the previous
section.

\begin{theorem}\label{thm:main}
Let $\Si$ be a hypersurface\index{hypersurface} in spacetime, and
let $\K$ denote the generalized second fundamental
form\index{second fundamental form!generalized} on $\Si$. Assume
that the components $g_{\alpha\beta}$ of the gravitational
metric\index{metric} $g$ are smooth on either side of $\Si,$
(continuous up to the boundary on either side separately), and
Lipschitz continuous\index{metric!Lipschitz continuous} across
$\Si$ in some fixed coordinate system.  Then the following
statements are equivalent:

\noindent(i) $[\K]=0$ at each point of $\Si.$

\noindent(ii) The components of the curvature
tensors\index{tensor!curvature} $R^{\mu}_{\alpha\nu\beta}$ and
$G_{\alpha\beta},$ viewed as second order operators on the metric
components $g_{\alpha\beta},$  produce no delta
functions\index{delta function singularity} on $\Si.$

\noindent(iii) For each point $p$ in $\Si$ there exists a
$C^{1,1}$ coordinate transformation defined in a neighborhood of
$p,$ such that, in the new coordinates, (which can be taken to be
the MGS coordinates for the surface), the metric components are
$C^{1,1}$ functions of these coordinates.

\noindent(iv) For each $p$ in $\Si,$ there exists a coordinate
frame that is locally Lorentzian\index{locally Lorentzian} at $p,$
and can be reached within the class of $C^{1,1}$ coordinate
transformations\index{$C^{1,1}$ coordinate transformation}.

Moreover, if any one of these equivalencies hold and $\Si$ is
non-null, then the Rankine-Hugoniot jump conditions\index{jump
conditions!Rankine-Hugoniot}, $[G]_{\beta}^{\alpha} N^{\beta}=0,$
which express the weak form of conservation of energy and momentum
across $\Si$ when $G=\kappa T$, hold at each point on $\Si.$

Lastly, if we add the condition that the second derivative of the
the metric inner product on $\tans^+$ is continuous, where
$\tans^+$ denotes the spacelike subspace of $\tans$, then the
Rankine-Hugoniot jump conditions\index{jump
conditions!Rankine-Hugoniot}, $[G]_{\beta}^{\alpha} N^{\beta}=0,$
also hold for a null-surface. In MGS coordinates, the condition
that the second derivative of the the metric inner product on
$\tans^+$ is continuous is equivalent to the condition
$\left[g_{tt,00}\right]=0,$ for $t=2,\ldots,n-1.$
\end{theorem}

\section{Supporting Results}
The proof of this theorem involves series of lemmas and theorems
which are developed below.

\subsection{The Rankine-Hugoniot Jump conditions}

The objective of this section is to write down the
Rankine-Hugoniot jump conditions\index{jump conditions} in terms
of our generalized second fundamental form\index{second
fundamental form!generalized}, which will give the criteria for
when the jump conditions hold. These jump conditions,
\begin{equation}\label{eq:rhjumpm}
[G]_{\beta}^{\alpha} N^{\beta}=0,
\end{equation}
express conservation of energy and momentum on the shock surface
$\Si$ when $G=\kappa T$ in the weak sense~\cite{st94, stgerm99}.
We begin with a lemma enables us to compute the connection
coefficients in terms of the generalized second fundamental
form\index{second fundamental form!generalized}.
\begin{lemma}\label{l4}
In MGS coordinates the components of the connection coefficients,
for a metric $g$ at a point $p$ in $\Si$, can be written as
\begin{align}
\Gamma^{\gamma}_{ab} & = \widetilde{\Gamma}^{\gamma}_{ab}
+g^{\gamma 0}\K_{ab}  \label{eq:l41}\\
\Gamma^{\gamma}_{\alpha 0} & = -\K^{\gamma}_\alpha, \label{eq:l42}
\end{align}
where $\gamma=0,\ldots,n-1,$, $a,b,s=1,\ldots,n-1,$ and
$\widetilde{\Gamma}^{\gamma}_{ab}$ indicates the summation of the
indices on the metric components goes from 1 to $n-1.$
\end{lemma}

\begin{proof}
Writing the connection coefficients in terms of the metric
components we see that
\begin{equation}
\begin{split}
\Gamma^{\gamma}_{ab}
    &=\frac{1}{2}g^{\gamma\sigma}\{-g_{ab,\sigma} +g_{\sigma a, b} +g_{b\sigma, a}\}\\
    &=\frac{1}{2}g^{\gamma s}\{-g_{ab,s} +g_{s a, b} +g_{bs,a}\}
      + \frac{1}{2}g^{\gamma 0}\{-g_{ab,0} +g_{0 a, b} +g_{b0,
      a}\}\\
    &=\widetilde{\Gamma}^{\gamma}_{ab}-\frac{1}{2}g^{\gamma
    0}g_{ab,0}\\
    &=\widetilde{\Gamma}^{\gamma}_{ab}+g^{\gamma 0}\K_{ab},
\end{split}
\end{equation}
which is equation~\eqref{eq:l41}. We use the same method to get
equation~\eqref{eq:l42}, that is,
\begin{equation}
\begin{split}
\Gamma^{\gamma}_{\alpha 0}
    &=\frac{1}{2}g^{\gamma\sigma}\{-g_{\alpha 0,\sigma} +g_{\sigma \alpha,0} +g_{0\sigma, \alpha}\}\\
    &=\frac{1}{2}g^{\gamma s}\{-g_{\alpha 0,s} +g_{s \alpha,0} +g_{0s,\alpha}\}
      + \frac{1}{2}g^{\gamma 0}\{-g_{\alpha 0,0} +g_{0 \alpha, 0} +g_{00,\alpha}\}\\
    &=\frac{1}{2}g^{\gamma s}g_{s\alpha,0}=-\K^{\gamma}_\alpha.
\end{split}
\end{equation}
\end{proof}
With the use of the above lemma we can now write the components
of the Riemann curvature tensor in terms of the components of the
generalized second fundamental form.

\begin{lemma}\label{l5} In MSG coordinates the components of the
Riemann curvature tensor can be written as
\begin{multline}\label{l51}
R^{\gamma}_{ajb}=g^{\gamma 0} \left( \K_{ab;j}-\K_{aj;b} \right)
+ \left( \K^{\gamma}_{b} \K_{aj}-\K^{\gamma}_{j} \K_{ab}\right)
\\+\left( \widetilde{\Gamma}^0_{aj}\K^{\gamma}_b
-\widetilde{\Gamma}^0_{ab}\K^{\gamma}_j \right) +g^{s0}\left(
\widetilde{\Gamma}^{\gamma}_{sj}\K_{ab}
-\widetilde{\Gamma}^{\gamma}_{sb}\K_{aj}
\right)+\widetilde{R}^{\gamma}_{ajb},
\end{multline}
\begin{equation}\label{l52}
R^\gamma_{aj0} =-\K^\gamma_{a,j}-\widetilde{\G}^\gamma_{aj,0}
-g^{\gamma 0}\K_{aj,0}
-\K^s_a\left(\widetilde{\G}^\gamma_{sj}+g^{\gamma 0}\K_{sj}\right)
+\K^\gamma_s\left(\widetilde{\G}^s_{aj} +g^{s0}\K_{aj}\right)
+\K^\gamma_j\K^0_a,
\end{equation}
and
\begin{equation}\label{l53}
R^{\gamma}_{0jb} =\K^\gamma_{j;b}-\K^\gamma_{b;j}.
\end{equation}
Recall that the Greek indices go from $0,\ldots,n-1$ and the
Latin indices go from $1,\ldots,n-1.$
\end{lemma}

\begin{proof}For equation~\eqref{l51} we will first consider three
equations that will be pertinent to the computation. Using
Lemma~\ref{l4}, we can write
\begin{equation}\label{eq:l51}
\begin{split}
g^{\gamma 0}\K_{ab,j} -g^{\gamma 0}&\K_{aj,b} -g^{\gamma
0}\K_{sb}\Gamma^s_{aj} +g^{\gamma 0}\K_{sj}\Gamma^s_{ab}\\
&=g^{\gamma 0} \left\{\K_{ab,j}-\K_{sb}\G^s_{aj}
-\left(\K_{aj,b}-\K_{sj}\G^s_{ab} \right)
\right\}\\
&=g^{\gamma 0} \left\{\K_{ab,j}-\K_{\sigma b}\G^\sigma_{aj}
-\left(\K_{aj,b}-\K_{\sigma
j}\G^\sigma_{ab} \right) \right\} \\
&=g^{\gamma 0}\left(\K_{ab;j}-\K_{aj;b}\right),
\end{split}
\end{equation}

\begin{equation}\label{eq:l52}
\begin{aligned}
\widetilde{\G}^\gamma_{sb} \G^s_{aj} &=\widetilde{\G}^\gamma_{sb}
\left( \widetilde{\G}^s_{aj}-\frac{1}{2}g^{s0}g_{aj,0} \right)
&=\widetilde{\G}^\gamma_{sb} \widetilde{\G}^s_{aj}
+g^{s0}\widetilde{\G}^\gamma_{sb} \K_{aj}\\
\widetilde{\G}^\gamma_{sj} \G^s_{ab} &=\widetilde{\G}^\gamma_{sj}
\widetilde{\G}^s_{ab}+g^{s0}\widetilde{\G}^\gamma_{sj} \K_{ab}, &
\end{aligned}
\end{equation}
and
\begin{equation}\label{eq:l53}
\begin{aligned}
\G^\gamma_{0j}\G^0_{ab} &= -\K^\gamma_j
\left(\widetilde{\G}^0_{ab} +g^{00} \K_{ab}\right)
&=-\widetilde{\G}^0_{ab} \K^\gamma_j
-\K_{ab}\K^\gamma_j \\
\G^\gamma_{0b}\G^0_{aj} &= -\widetilde{\G}^0_{aj} \K^\gamma_b
-\K_{aj}\K^\gamma_b.  &
\end{aligned}
\end{equation}

Using Lemma~\ref{l4}, and equations~\eqref{eq:l51},
\eqref{eq:l52}, \eqref{eq:l53} we have
\begin{equation}\label{eq:l54}
\begin{aligned}
R^\gamma_{ajb} & = \G^\gamma_{ab,j} -\G^\gamma_{aj,b}
+\G^\gamma_{\sigma j}\G^\sigma_{ab} -\G^\gamma_{\sigma b} \G^\sigma_{aj}  \\
& = \G^\gamma_{ab,j} -\G^\gamma_{aj,b}+\G^\gamma_{s j}\G^s_{ab}
-\G^\gamma_{s b} \G^s_{aj} +\G^\gamma_{0 j}\G^0_{ab} -\G^\gamma_{0
b} \G^0_{aj}\\
& =\widetilde{\G}^\gamma_{ab,j}+g^{\gamma 0}\K_{ab,j}
-\widetilde{\G}^\gamma_{aj,b}-g^{\gamma 0}\K_{aj,b}
-\G^s_{aj}\left(\widetilde{\G}^\gamma_{sb}+g^{\gamma
0}\K_{sb}\right)\\&\qquad+\G^s_{ab}\left(\widetilde{\G}^\gamma_{sj}+g^{\gamma
0}\K_{sj}\right)+\G^\gamma_{0 j}\G^0_{ab} -\G^\gamma_{0
b} \G^0_{aj}\\
&= g^{\gamma 0} \left\{\K_{ab,j}-\K_{sb}\G^s_{aj}
-\left(\K_{aj,b}-\K_{sj}\G^s_{ab} \right)
\right\}\\
& \qquad +\widetilde{\G}^\gamma_{sj} \G^s_{ab}
-\widetilde{\G}^\gamma_{sb} \G^s_{aj}
+\widetilde{\G}^\gamma_{ab,j} -\widetilde{\G}^\gamma_{aj,b}
+\G^\gamma_{0 j}\G^0_{ab}
-\G^\gamma_{0 b} \G^0_{aj}\\
&=g^{\gamma 0} \left( \K_{ab;j}-\K_{aj;b} \right) + \left(
\K^{\gamma}_{b} \K_{aj}-\K^{\gamma}_{j} \K_{ab}\right)  +\left(
\widetilde{\Gamma}^0_{aj}\K^{\gamma}_b
-\widetilde{\Gamma}^0_{ab}\K^{\gamma}_j \right) \\ &
\qquad+g^{s0}\left( \widetilde{\Gamma}^{\gamma}_{sj}\K_{ab}
-\widetilde{\Gamma}^{\gamma}_{sb}\K_{aj}
\right)+\widetilde{R}^{\gamma}_{ajb},
\end{aligned}
\end{equation}
where $\widetilde{R}^{\gamma}_{ajb}=\widetilde{\G}^\gamma_{ab,j}
-\widetilde{\G}^\gamma_{aj,b} +\widetilde{\G}^\gamma_{s
j}\widetilde{\G}^s_{ab} -\widetilde{\G}^\gamma_{s b}
\widetilde{\G}^s_{aj}.$

For equations~\eqref{l52} and~\eqref{l53}, we also use
Lemma~\ref{l4} to write
\begin{equation}\label{eq:l74}
\begin{aligned}
R^\gamma_{aj0} &= \G^\gamma_{a0,j}-\G^\gamma_{aj,0}
+\G^\gamma_{sj}\G^s_{a0}-\G^\gamma_{s0}\G^s_{aj}
+\G^\gamma_{0j}\G^0_{a0}-\G^\gamma_{00}\G^0_{aj} \\
&=-\K^\gamma_{a,j}-\widetilde{\G}^\gamma_{aj,0} -g^{\gamma
0}\K_{aj,0} -\K^s_a\left(\widetilde{\G}^\gamma_{sj}+g^{\gamma
0}\K_{sj}\right) \\&\qquad +\K^\gamma_s\left(\widetilde{\G}^s_{aj}
+g^{s0}\K_{aj}\right) +\K^\gamma_j\K^0_a ,
\end{aligned}
\end{equation}
and
\begin{equation}
\begin{aligned}
R^{\gamma}_{0jb}
 &= \G^\gamma_{0b,j} -\G^\gamma_{0j,b}
+\G^\gamma_{\sigma j}\G^\sigma_{0b} -\G^\gamma_{\sigma b} \G^\sigma_{0j}\\
&= -\K^\gamma_{b,j} +\K^\gamma_{j,b} -\G^\gamma_{\sigma
j}\K^\sigma_{b} +\G^\gamma_{\sigma b}
\K^\sigma_{j}\\
&= \K^\gamma_{j,b} +\G^\gamma_{\sigma b}\K^\sigma_{j}
-\left(\K^\gamma_{b,j}
+\G^\gamma_{\sigma j}\K^\sigma_{b}\right) \\
 &=\K^\gamma_{j;b}-\K^\gamma_{b;j}.
\end{aligned}
\end{equation}
\end{proof}

\begin{remark}\label{xjump} For any quantities $A$, and $B$ such that both $[A]=0$
and $[B]=0$, the jump in the product of the quantities is
$$[AB]=A^LB^L-A^L B^R +A^L B^R -A^R B^R=A^L[B] +
B^R[A]=0.$$
\end{remark}
The next three lemmas are more of the technical variety which we
also need to facilitate computations to be made later on.

\begin{lemma}\label{l4a}If $\left[\K_{ab} \right]=0$ for every point
in $\Si,$ then $\left[\K_{ab;j} \right]=0$ for every point in
$\Si.$
\end{lemma}

\begin{proof} Using Lemma~\ref{l4} we have
\begin{equation}\label{eq:l4a1}
\begin{aligned}
\left[\K_{ab;j} \right] & =\left[\K_{ab,j}-\K_{\sigma
b}\G^\sigma_{aj} \right]\\ & =\left[\K_{ab,j} \right]- g^{\sigma
0} \left[\K_{\sigma b}\K_{aj}\right] +g^{\sigma
0}\left[\K_{\sigma b}\widetilde{\G}^\sigma_{aj}\right].
\end{aligned}
\end{equation}
Since $\left[\smash{\K_{ab}}\right]=0$ for every point in $\Si$
implies that $[\K_{ab,j}]=0$ for $j=1,\ldots,n-1,$ and
$[\widetilde{\G}^\sigma_{aj}]=0$, the result holds by
equation~\eqref{eq:l4a1}.
\end{proof}

\begin{lemma}\label{l6} If $[\K]=0$, then $\left[R^{\gamma a}_{jb}\right]=0.$
\end{lemma}

\begin{proof}
If $[\K]=0$, then by lemmas~\ref{l5} and~\ref{l4a},
remark~\ref{xjump}, $[\widetilde{\G}^\gamma_{ab}]=0,$ and
$[\widetilde{R}^\gamma_{ajb}]=0,$ we can conclude both
\begin{equation}
\left[R^\gamma_{ajb}\right]=0 \qquad \mbox{and} \qquad
\left[R^\gamma_{0jb}\right]=0.
\end{equation}
Therefore, if $[\K]=0$, then
\begin{equation}\label{eq:l61}
\left[R^{\gamma a}_{jb}\right] =g^{a\sigma}\left[R^\gamma_{\sigma
jb}\right] =g^{as}\left[R^\gamma_{s jb}\right]
+g^{a0}\left[R^\gamma_{0 jb}\right]=0.
\end{equation}
\end{proof}

\begin{lemma}\label{l6r} If $[\K]=0$, then in MGS coordinates
$$\left[R^{at}_{0\,t}\right]=g^{a0}[g_{tt,00}],$$
for $a=1,\ldots,n-1,$ and $t>1.$
\end{lemma}

\begin{proof} First,
\begin{equation}\label{eq:l6r1}
R^{at}_{0\,t}=g^{t\sigma}R^a_{\sigma 0t}=-g^{tt}R^a_{tt0},
\end{equation}
since for $t>1$, $g^{tt}$ are the only non-zero components. We now
compute $R^1_{tt0}$ in MGS coordinates by using lemma~\ref{l4}
which yields
\begin{equation}\label{eq:l6r2}
\begin{aligned}
R^a_{tt0} &= \G^a_{t0,t}-\G^a_{tt,0}
+\G^a_{st}\G^s_{t0}-\G^a_{s0}\G^s_{tt}
+\G^a_{0t}\G^0_{t0}-\G^a_{00}\G^0_{tt} \\
&=-\K^a_{t,t}-\widetilde{\G}^a_{tt,0} -g^{a 0}\K_{tt,0}
-\K^s_t\left(\widetilde{\G}^a_{st}+g^{a 0}\K_{st}\right)
\\&\qquad +\K^a_s\left(\widetilde{\G}^s_{tt}
+g^{s0}\K_{tt}\right) +\K^a_t\K^0_t.
\end{aligned}
\end{equation}
When $\left[\K \right]=0$ we have
$$\left[\K^a_{t,t} \right] = g^{a\sigma}\left[\K_{\sigma
t,t}\right]=0,$$
$$\left[\tilde{\G}^a_{tt,0} \right]
=\frac{a}{2\eta^2}\left[g_{tt,0} \right]_{,a}=0,$$ and the
remaining terms are zero, when $[\K]=0$, except
\begin{equation}\label{eq:nojump}
\left[g^{a 0}\K_{tt,0} \right]=g^{a0}[g_{tt,00}]
\end{equation}
when $g^{a0}\neq 0$.
\end{proof}

The next result gives us a nice way to express the components of
the Einstein curvature tensor in terms of the components of
Riemann curvature tensor. We reproduce the statement of the
following lemma, and its proof, given
in~\cite{stg01},~\cite{st94}, and~\cite{stgerm99}. These same
computations can also be found in~\cite[section 14.2]{mtw}.

\begin{lemma}\label{l1} The components of the Einstein curvature
tensor can be written as
\begin{equation}\label{eq:l11}
G^{\alpha}_{\alpha}=-\sum_{\sigma,\tau\neq
\alpha}R^{|\sigma\tau|}_{|\sigma\tau|},\quad \alpha=1,\cdots,n
\end{equation}
and
\begin{equation}\label{eq:l12}
G^\alpha_\beta=\sum_{\tau\neq
\alpha,\beta}R^{|\alpha\tau|}_{|\beta\tau|},\quad \alpha\neq
\beta.
\end{equation}
The indices inside $|\cdot|$ are always taken as an increasing
sequence.
\end{lemma}

\begin{proof} We begin by raising the index of the components of
the Einstein tensor defined in section~\ref{s:spacetime}, and
given by
$$G_{\alpha\beta}=R_{\alpha\beta}-\frac{1}{2}g_{\alpha\beta}R,$$
from which we find
\begin{equation}\label{pl11}
G^\alpha_\beta=R^\alpha_\beta-\frac{1}{2}g^\alpha_\beta R
=R^\alpha_\beta-\frac{1}{2} \delta^\alpha_\beta R.
\end{equation}
Since $R^{\sigma \tau}_{\mu \nu}$ is antisymmetric in
$\sigma\tau$ and also in $\mu \nu$, it follows that
\begin{equation}\label{pl12}
R^\alpha_\alpha=R^{\alpha \nu}_{\alpha
\nu}=\sum_{\alpha\neq\nu}R^{\alpha \nu}_{\alpha \nu}=R^{|\alpha
\nu|}_{|\alpha \nu|}.
\end{equation}
Notice that we do not sum over the index $\alpha$ since it
represents the component, and not a summation index. The Ricci
Scalar can be written as
\begin{equation}\label{pl13}
  R=R^{\sigma\tau}_{\sigma\tau}=2R^{|\sigma\tau|}_{|\sigma\tau|}.
\end{equation}
Combining equations~\eqref{pl12} and~\eqref{pl13} we have
equation~\eqref{eq:l11},
$$G^{\alpha}_{\alpha}
=2R^{|\alpha \nu|}_{|\alpha \nu|}-R^{|\sigma
\tau|}_{|\sigma\tau|} =-\sum_{\sigma,\tau\neq
\alpha}R^{|\sigma\tau|}_{|\sigma\tau|}.$$ Furthermore, we also
find equation~\eqref{eq:l12} by the computation
$$G^\alpha_\beta=R^\alpha_\beta=R^{\alpha\tau}_{\beta\tau}
=\sum_{\tau\neq \alpha,\beta}R^{|\alpha\tau|}_{|\beta\tau|}. $$
\end{proof}

Now we are finally able to directly prove a part of our main
result, Theorem~\ref{thm:main}. The following theorem gives the
criteria for when the Rankine-Hugoniot jump conditions hold across
$\Si$.

\begin{theorem}\label{lrh}
If both $[g_{tt,00}]=0$ for $t>1$ in MGS coordinates, and $[\K]=0$
at each point of $\Si$, then $[G_{\alpha\beta}]N^{\alpha}=0$. The
condition $[g_{tt,00}]=0$ for $t>1$ in MGS coordinates means the
metric is $C^2$ in coordinates not associated with the transverse
vector $\N$ and the $\n$ to the surface.
\end{theorem}

\begin{proof}
In MGS coordinates $N^0=1$ and $N^a=0$ for $a=1,\ldots n-1,$
therefore
\begin{equation}\label{eq:lrh1}
[G_{\alpha_\beta}]N^\alpha=[G_{0\beta}]
=g_{00}[G^0_\beta]+g_{01}[G^1_\beta].
\end{equation}
When $\beta=0,$ we have
\begin{equation}\label{eq:lrh2a}
[G_{0\beta}] =[G_{00}]=g_{00}[G^0_0]+g_{01}[G^1_0],
\end{equation}
$\beta=1$,
\begin{equation}\label{eq:lrh2b}
[G_{0\beta}] =[G_{01}]=g_{10}[G^0_0]+g_{11}[G^1_0]=g_{10}[G^0_0],
\end{equation}
and when $\beta=t>1$ we have
\begin{equation}\label{eq:lrh2c}
[G_{0\beta}] =[G_{0t}]=g_{i0}[G^t_0].
\end{equation}
See equation~\eqref{eq:gmgsmatrix} for the components of $g$ in
MGS coordinates. By lemma~\ref{l1} we have
$$\left[G^0_0\right]=\left[-\sum_{s,t\neq 0}R^{|st|}_{|st|}\right]
=-\frac{1}{2}\sum_{s,t\neq 0} [R^{st}_{st}],$$ and, by
lemma~\ref{l6}, $[\K]=0$ implies that $\left[R^{st}_{st}\right]=0$
for all $s,t=1,\ldots,n-1,$ hence $[G^0_0]=0.$ Also, by
lemma~\ref{l1},
$$\left[G_0^a\right]=\left[\,\sum_{t\neq 0,a}R^{|at|}_{|0t|}\right]
=\sum_{t\neq 0,a}\left[R^{|at|}_{|0\,t|}\right].$$ Now,
lemma~\ref{l6r}, $[\K]=0$ implies that
$$[R^{at}_{0\,t}]=g^{a0}[g_{tt,00}],$$
for all $t=2,\ldots,n-1.$ Then $[g_{tt,00}]=0$ for $t>1$ implies
that $[G^a_0]=0.$ Therefore, in MGS coordinates,
$$\left[G_{\alpha\beta}\right]N^\alpha =0.$$
\end{proof}

\subsection{Curvature Tensor Components as Delta Functions}
\label{delta function} The physical condition for the hypersurface
$\Si$  to be considered a shock surface, when $G=\kappa T$, is
that there exist no delta function singularities\index{delta
function singularity} at a point $p$ in $\Si$ in the components of
the Einstein curvature tensor $G$.

What does this mean? When considering weak solutions of the
Einstein equations $G=\kappa T$ across a discontinuity if there
exists a delta function singularity in a component of the Einstein
tensor, then this singularity manifests itself in the coupled
stress energy tensor $T$ as an infinite spike in the mass/energy
density of matter at the discontinuity~\cite{isr66}. In this case
the hypersurface of discontinuity is referred to as a
\emph{surface layer}\index{surface layer}. On the other hand, if
there exist no delta function singularities in the Einstein
tensor, then there is merely a discontinuity in the mass/energy
density, and the surface is considered a shock.

We will show below that there exist no delta function
singularities in the Riemann curvature tensor if and only if
$[\K]=0$ for every point on $\Si,$ and hence will give a necessary
and sufficient condition for a hypersurface to be considered a
shock surface.

To see this, notice that each component of the Einstein tensor
can be written in terms of the metric components
$g_{\alpha\beta}$, its first and second weak derivatives, and the
components of the inverse of the metric $g^{\alpha\beta}$. Since
$g_{\alpha\beta}$ and $g^{\alpha\beta}$ are Lipschitz continuous
across $\Si$ it follows that the only possible delta function
singularities contained in the the Einstein tensor are the second
order derivatives, $g_{\alpha\beta,\sigma\tau}$, of the metric
components. Also, notice that derivatives in the direction
tangent to $\Si$ always cancel out because $[g_{\alpha\beta}]=0$
at all points of $\Si.$ This means that, in MGS coordinates,
$[g_{\alpha\beta,s}]=0$ for $s=1,\ldots,n-1,$ and the only
possible delta functions are of the form $g_{ab,00}$. Therefore
in MGS coordinates, when $[\K_{ab}]=[g_{ab,0}]=0$, the metric $g$
is $C^1$ just as we saw in corollary~\ref{c32}, and $g_{ab,00}$
is not a delta function. However, when
$[\K_{ab}]=[g_{ab,0}]\neq0$, then $g_{ab,00}$ is a delta
function.

In the following lemmas we will show that components of the
Riemann curvature tensor\index{Riemann curvature tensor}, and
hence the components of the Ricci tensor\index{Ricci tensor}, the
Ricci scalar\index{Ricci scalar}, and the Einstein
tensor\index{Einstein tensor} contain no delta functions if and
only if $[\K]= 0.$

\begin{lemma}\label{l7}
In a MGS coordinate frame, with respect to the surface $\Si$, the
components of the Ricci tensor can be written as follows:
\begin{align}
R_{00} &=-\frac{1}{2}g^{ab}g_{ab,00}+\,\mbox{lower order 0-derivatives},\label{l73}\\
R_{a0} &=-\frac{1}{2}g^{i0}g_{ai,00}+\,\mbox{lower order
0-derivatives},\label{l72}\\
 R_{ab} &=\frac{1}{2}g_{ab,00}+\,\mbox{lower order
0-derivatives}.\label{l71}
\end{align}
Furthermore, the Ricci scalar can be written as
\begin{equation}\label{l7b1}
R=-g^{s0}g^{t0}g_{st,00}+\,\mbox{lower order 0-derivatives}.
\end{equation}
\end{lemma}

\begin{proof}
From equations~\eqref{l51}~\eqref{l52} and~\eqref{l53} we have
\begin{align}\label{eq:l51b}
R^{\gamma}_{ajb}&=\;\mbox{lower order \emph{0}-derivatives},\\
\label{eq:l51c} R^\gamma_{aj0} &= -\frac{1}{2}g^{\gamma
0}g_{aj,00} +\,\mbox{lower order \emph{0}-derivatives},\\
\label{eq:152c} R^{\gamma}_{0jb}&=\;\mbox{lower order
\emph{0}-derivatives}.
\end{align}
Thus, using Lemma~\ref{l4} we can write equation~\eqref{l73},
\begin{equation}\label{eq:175}
\begin{aligned}
R_{00}&=R^a_{0a0}+R^0_{000}= R^a_{0a0}\\
&=\G^a_{00,a}-\G^a_{0a,0} +\G^a_{\sigma a }\G^\sigma_{00}
-\G^a_{\sigma 0}\G^\sigma_{0a} \\
&=\K^a_{a,0}-\G^a_{s 0}\G^s_{0a}\\
&=-\frac{1}{2}g^{as}g_{as,00}-\K^a_s\K^s_a \\
&=-\frac{1}{2}g^{as}g_{as,00}+\;\mbox{lower order
\emph{0}-derivatives},
\end{aligned}
\end{equation}
where we have used $\G^\gamma_{00}=0.$ Also, from
equation~\eqref{eq:l51c} we have
\begin{equation}
R_{a0}=R^s_{as0}+R^0_{a00}=R^s_{as0} =-\frac{1}{2}g^{s
0}g_{as,00} +\,\mbox{lower order \emph{0}-derivatives},
\end{equation}
and from Lemma~\ref{l5} and equation~\eqref{eq:l74} we have
\begin{equation}
R_{ab}= R^s_{asb}+R^0_{a0b}=\frac{1}{2}g_{ab,00}+\,\mbox{lower
order 0-derivatives}.
\end{equation}
Using the results above we can write
$$\begin{aligned}
R &= g^{\alpha\beta}R_{\alpha\beta}
=g^{ab}R_{ab}+g^{a0}R_{a0}+g^{0b}R_{0b}+g^{00}R_{00} \\
&=g^{ab}R_{ab}+2\,g^{a0}R_{a0}+g^{00}R_{00}\\
&=\frac{1}{2}g^{ab}g_{ab,00}
+2\,g^{a0}\left(-\frac{1}{2}g^{t0}g_{at,00}\right)
-\frac{1}{2}g^{ab}g_{ab,00} +\,\mbox{lower order 0-derivatives}\\
&=- g^{a0}g^{t0}g_{at,00}+\,\mbox{lower order 0-derivatives}.
\end{aligned}$$
\end{proof}

\begin{corollary}\label{c7}
Assume $g=g^L\cup g^R$ is smooth on either side of $\Si$, and
Lipschitz continuous across $\Si.$ Then in MGS coordinates the
jump condition $[\K]=0$ exists at a point $p$ in $\Si$, if and
only if, the curvature tensors\index{tensor!curvature}
$R^{\alpha}_{\beta\gamma\delta}$, $R_{\alpha\beta}$,
$G_{\alpha\beta},$ and the Ricci scalar $R$, viewed as second
order differential operators in the weak sense on the metric
components $g_{\alpha\beta},$  produce no delta function
singularities\index{delta function singularity} at the point $p$
in $\Si.$
\end{corollary}

\begin{proof} By applying lemmas~\ref{l5} and~\ref{l7}, we see that,
in MGS coordinates, the components
$R^{\alpha}_{\beta\gamma\delta}$, $R_{\alpha\beta}$,
$G_{\alpha\beta},$ and the Ricci scalar $R$ can be written in the
form
$$ A(\x)g_{ab,00}+\,\mbox{lower order \emph{0}-derivatives},$$ where
$A(\x)$ is some function of the coordinates which is as least
Lipschitz continuous. From Lemma~\ref{l3} we know, in MGS
coordinates, $[\K]=0$ at a point $p$ in $\Si$, if and only if,
$[g_{ab,0}]=0$ for all $a,b=1,\ldots,n-1$ at the point $p$ in
$\Si.$ We also know $g_{ab,00}$ is not a delta function, if and
only if, $[g_{ab,0}]=0$. This completes the proof.
\end{proof}

Now our goal is to generalize corollary~\ref{c7} for any
coordinate system, not just for a MGS coordinate system. To do
this we will need the following lemma which considers how the
Riemann curvature tensor, defined in terms of second order weak
derivatives of the metric components $g_{\alpha\beta},$
transforms from one set of coordinates to another by a $C^{1,1}$
transformation.  Recall by a $C^{1,1}$ transformation we mean a
function whose first derivatives are Lipschitz continuous. The
lemma was originally proved by Temple and Smoller in~\cite{st94},
and can also be found in~\cite{stg01} and in~\cite{stgerm99}. We
will reproduce their result here for convenience.

\begin{lemma}\label{l8}
Let $\mathcal{R}=R^\lambda_{\mu\nu\xi}=L[g]$ denote the
components of the Riemann curvature tensor in $x$-coordinates
where $L$ is the second order linear operator on the metric
components $g_{\alpha\beta}$ which defines
$R^\lambda_{\mu\nu\xi}$. Similarly, let
$\bar{\mathcal{R}}=\bar{R}^\alpha_{\beta\gamma\delta}=L[\bar{g}]$
denote the components in $y$-coordinates which are related to
$x$-coordinates by a $C^{1,1}$ transformation.

If $\mathcal{R}$ is a weak solution of $\mathcal{R}=L[g]$ in
$x$-coordinates, then $\mathcal{R}\frac{\partial x}{\partial y}$
is a weak solution of $\bar{\mathcal{R}}=L[\bar{g}]$ for any
coordinate system $y$ related to $x$ by a $C^{1,1}$
transformation. Here we have used the short-hand notation
$$
\mathcal{R}\frac{\partial x}{\partial y} =R^\lambda_{\mu\nu\xi}
\frac{\partial x^\mu}{\partial y^\beta} \frac{\partial
x^\nu}{\partial y^\gamma}\frac{\partial x^\xi}{\partial
y^\delta}\frac{\partial y^\alpha}{\partial x^\lambda},
$$
and multiplication by a function is taken in the weak sense.
\end{lemma}

\begin{proof}
Let $g$ be smooth and $\varphi$ be an arbitrary smooth test
function with compact support. Furthermore, let
$$\int_{\Real^4} L[g]\varphi = \int_{\Real^4} L^*[g,\varphi],$$
where $L^*[g,\varphi]$ is defined as the expression obtained from
$L[g]$ by integrating the second order derivatives in $g$ once by
parts. Now, in any coordinate system, $L$ is given by
\begin{equation}\label{eq:l81}
\begin{aligned}
R^\lambda_{\mu\nu\xi}=L[g] &=
\G^\lambda_{\mu\xi,\nu}-\G^\lambda_{\mu\nu,\xi} +\;\mbox{lower
order \emph{0}-derivatives}\\
&= \left(g^{\lambda\sigma}\{-g_{\mu\xi,\sigma}+g_{\sigma\mu,\xi}
+g_{\xi\sigma,\mu}\}\right)_{,\nu}
-\left(g^{\lambda\sigma}\{-g_{\mu\nu,\sigma}+g_{\sigma\mu,\nu}
+g_{\nu\sigma,\mu}\}\right)_{,\xi} \\ & \hspace{1.0in}
+\;\mbox{lower order \emph{0}-derivatives} \\
&=g^{\lambda\sigma}\left( -g_{\mu\xi,\sigma\nu}
+g_{\xi\sigma,\mu\nu} +g_{\mu\nu,\sigma\xi}
-g_{\nu\sigma,\mu\xi}\right) +\;\mbox{lower order
\emph{n}-derivatives}.
\end{aligned}
\end{equation}
Thus, $L^*[g,\varphi]$ is composed of the metric components
$g_{\mu\nu}$, the test function $\varphi$, their first
derivatives, and the inverse metric entries $g^{\mu\nu}$.
Therefore, $L^*[g,\varphi]$ is integrable over $\Real^4$ for any
Lipschitz continuous metric and any Lipschitz continuous test
function $\varphi$ with compact support.

Now, suppose
$$\langle\mathcal{R},\varphi\rangle\equiv \int_{\Real^4}\mathcal{R}\varphi
=\int_{\Real^4} L^*[g,\varphi],$$ for all Lipschitz continuous
test functions $\varphi$ with compact support, that is,
$\mathcal{R}$ is a weak solution of $\mathcal{R}=L[g].$ Let
$$\bar{g}\equiv\bar{g}_{\alpha\beta}
=g_{\mu\nu}\frac{\partial x^\mu}{\partial y^\alpha}
\frac{\partial x^\nu}{\partial y^\beta}=g\frac{\partial
x}{\partial y}.$$ If $\frac{\partial x}{\partial y}$ is Lipschitz
continuous, then $L^*\left[g(\partial x/\partial
y),\varphi\right]$ is  bounded for any Lipschitz continuous test
function.

Suppose $g$ is an arbitrary, non-degenerate, Lipschitz continuous
metric, and $\varphi$ an arbitrary Lipschitz continuous test
function. Furthermore, suppose the coordinates $x$ and $y$ are
related by a $C^{1,1}$ transformation, that is, $\frac{\partial
x}{\partial y},\,\frac{\partial y}{\partial x}\in C^{0,1}$ Let
$\bar{g}^\epsilon_{\alpha\beta}$ denote the approximation of
$\bar{g}_{\alpha\beta}$ by a smooth function, and also let
$x^\epsilon(y)$ denote the approximation of the coordinate map
$x(y)$ where $x^\epsilon(y)$ is smooth, and has a smooth
inverse\footnote{In our case here, it is sufficient that if $u$
is a continuous function on a bounded set, then there exists a
smooth function $u^\epsilon$ such that
$\|u^\epsilon-u\|_{\infty}\To 0$ as $\epsilon\To 0.$ For more on
approximations by smooth functions see~\cite[Section 5.3]{evans}
or~\cite[page 62]{smol}.}. These approximations can be chosen so
that
\begin{align*}
&\bar{g}^\epsilon_{\alpha\beta} \To \bar{g}_{\alpha\beta}
\;\mbox{in}\;C^{0,1}, \quad x^\epsilon(y) \To
x(y)\;\mbox{in}\;C^{1,1},
\\ &\frac{\partial x^\epsilon}{\partial y}(y) \To \frac{\partial
x}{\partial y}(y) \;\mbox{in}\;C^{0,1},\; \mbox{and} \quad
\frac{\partial y}{\partial x^\epsilon}(x^\epsilon) \To
\frac{\partial y}{\partial x}(x)\;\mbox{in}\;C^{0,1}.
\end{align*}
Consequently,
$$g^\epsilon\equiv\bar{g}^\epsilon
\frac{\partial y}{\partial x^\epsilon} \To g,\quad
\mbox{and}\quad \bar{g}^\epsilon\To \bar{g}\; \mbox{in}\;
C^{0,1}.$$  Now, if we define
\begin{equation}\label{eq:l83}
\bar{\mathcal{R}}^\epsilon=L[g^\epsilon],
\end{equation}
and
\begin{equation}\label{eq:l84}
 \mathcal{R}^\epsilon=\bar{\mathcal{R}}^\epsilon
 \frac{\partial y}{\partial x^\epsilon},
\end{equation}
then, by definition,
\begin{equation}\label{eq:l85}
\left\langle \mathcal{R}^\epsilon \frac{\partial
x^\epsilon}{\partial y}, \varphi \right\rangle =\left\langle
\bar{\mathcal{R}}^\epsilon,\varphi \right\rangle = \int_{\Real^4}
L^*\left[g^\epsilon \frac{\partial x^\epsilon}{\partial
y},\varphi \right],
\end{equation}
for all test functions $\varphi.$ From this equation we can say
$\mathcal{R}^\epsilon \frac{\partial x^\epsilon}{\partial y}$ is
the Riemann curvature tensor for $g^\epsilon \frac{\partial
x^\epsilon}{\partial y}.$ Now, because every function in
equation~\eqref{eq:l85} is sufficiently smooth, the equation
holds, if and only if,
\begin{equation}\label{eq:l86}
\left\langle \mathcal{R}^\epsilon, \varphi \right\rangle =
\int_{\Real^4} L^*\left[g^\epsilon, \varphi \right]
\end{equation}
Since $g^\epsilon \rightarrow g$ as $\epsilon \rightarrow 0$ in
$C^{0,1},$ if follows that
$$\int_{\Real^4} L^*\left[g^\epsilon, \varphi \right] \To
\int_{\Real^4} L^*\left[g, \varphi \right] \;\mbox{as}\; \epsilon
\rightarrow 0. $$ Therefore, from equation~\eqref{eq:l86},
$\mathcal{R}^\epsilon \rightarrow \mathcal{R}$ as $\epsilon
\rightarrow 0$ in the weak sense. In the same way,
$\bar{\mathcal{R}}^\epsilon \rightarrow \bar{\mathcal{R}},$ and
thus, by equation~\eqref{eq:l84}, we conclude $\bar{\mathcal{R}}
=\mathcal{R}\frac{\partial x}{\partial y}$ in the weak sense,
which completes the proof.
\end{proof}

With Lemma~\ref{l8}, we are now able to prove the result in
corollary~\ref{c7} for any coordinate system with our Lipschitz
continuous metric $g=g^L\cup g^R$.

\begin{theorem}\label{lcurvtens} Assume $g=g^L\cup g^R$ is smooth
on either side of $\Si$, and Lipschitz continuous across $\Si.$
Then the jump condition $[\K]=0$ exists at a point $p$ in $\Si$,
if and only if, the curvature tensors\index{tensor!curvature}
$R^{\alpha}_{\beta\gamma\delta}$, and $G_{\alpha\beta}$ viewed as
second order differential operators in the weak sense on the
metric components $g_{\alpha\beta},$  produce no delta function
singularities\index{delta function singularity} at $p$ in $\Si.$
\end{theorem}

\begin{proof}
By corollary~\ref{c7}, the theorem holds in MGS coordinates. Now,
for any metric $g=g^L\cup g^R$ which is smooth on either side of
$\Si$, and Lipschitz continuous across $\Si,$ a transformation
from our MGS coordinate system to another coordinate system must
be least $C^{1,1}$ by lemma~\ref{l:C11trans}. Since the
transformation is an invertible $C^{1,1}$ function, by
Lemma~\ref{l8}, the theorem holds for any coordinate system.
\end{proof}

\subsection{Relating to a $C^{1,1}$ metric by a $C^{1,1}$ transformation}
Here we prove a result that relates values of $[\K]$ to the
smoothness of the metric. In particular, the condition that
$[\K]=0$ for each point $p$ in $\Si$ is equivalent to the
existence of a coordinate system at each point in $\Si$ such that
the metric components are $C^1$ functions in these coordinates.

\begin{theorem}\label{c11}
Assume $g=g^L\cup g^R$ is smooth on either side of $\Si$, and
Lipschitz continuous across $\Si$ in some fixed coordinate system
at each point in $\Si$. Then for each point $p$ in $\Si$,
$[\K]=0$, if and only if, there exists another coordinate system,
defined in a neighborhood of $p,$ such that the metric components
are $C^{1,1}$ functions of these coordinates, and are related to
the original coordinates by a $C^{1,1}$ coordinate
transformation.
\end{theorem}

\begin{proof}
The initial coordinate system is related to an MGS coordinate
system by a $C^{1,1}$ transformation. By corollary~\ref{c32}, if
$[\K]=0$, then the metric components in the MGS coordinate system
are $C^1.$ Furthermore, since $g$ is $C^2$ on each side of $\Si$
it follows that $g_{\alpha\beta,\mu\nu}$ is bounded for each
$\alpha,\beta,\gamma=0,\ldots,n-1$, therefore each components is
$C^{1,1}.$

Conversely, if the metric components $g_{\alpha\beta}$ in the
original coordinates are equivalent to a $C^{1,1}$ metric in
another coordinate system by a $C^{1,1}$ transformation, then in
an MGS coordinate system, which is related to the latter
coordinates by a $C^{1,1}$ transformation, the metric components
will be $C^{1,1}$ functions of these latter coordinates.
Therefore, by corollary~\ref{c32}, we have $[\K]=0$, and the
theorem is proved.
\end{proof}

\subsection{Existence of Local Lorentzian Frame Via a $C^{1,1}$ Transformation}

\begin{theorem}\label{lorentz}
Assume $g=g^L\cup g^R$ is smooth on either side of $\Si$, and
Lipschitz continuous across $\Si$ in some fixed coordinate system
at each point in $\Si$. The $[\K]=0$ at each point of $\Si$, if
and only if, for each $p$ in $\Si,$ there exists a coordinate
system that is locally Lorentzian\index{locally Lorentzian} at
$p,$ and is related to the original coordinates by a $C^{1,1}$
coordinate transformation\index{$C^{1,1}$ coordinate
transformation}.
\end{theorem}

\begin{proof}
If $[\K]\neq 0$ at some point $p$ in $\Si$, then, by
Theorem~\ref{c11}, there exists no coordinate system defined in a
neighborhood of $p$, and related to the original coordinates by a
$C^{1,1}$ coordinate transformation,  such that the metric
components are $C^{1}$ functions of these coordinates. As a
consequence, there cannot exist a Lorentzian coordinate frame
containing $p$ and related to the original coordinates by a
$C^{1,1}$ coordinate transformation.

Conversely, suppose $[\K]=0$ at every point in $\Si.$ Denote the
original $n$ dimensional coordinate system, at a point $p$ in
$\Si$, by $x=(x^0,\dots,x^{n-1}).$ Then, by
lemma~\ref{l:C11trans}, these coordinates are related to a MGS
coordinate system by a $C^{1,1}$ transformation. By
Lemma~\ref{posdef} we know that $g_{ii}>0,$ for $i=2,\ldots,n-1,$
in any MGS coordinate system. Thus, we can choose $g_{ii}=1$ and
$g_{ii,j}=0$ at $p$ for $i,j=2,\ldots,n-1$. Now, we construct a
coordinate system, denoted by $u=(u^0,\ldots,u^{n-1}),$ just as we
did in Remark~\ref{rmk:lorentzcoord}, so that the metric is of the
form $g=\mbox{diag}(-1,1,\ldots,1)$. By the construction of this
coordinates system (see the proof of Lemma~\ref{posdef} and
Remark~\ref{rmk:lorentzcoord}), the transformation from our MGS
coordinate system to $u$ coordinates is in the class of $C^{1,1}$
functions. Therefore, it remains to show that
$g_{\alpha\alpha,\beta}=0$ for all $\alpha,\beta=0,\ldots,n-1$.
For $i,j=2,\ldots,n-1$, we already have $g_{ij,a}=0$ where
$a=1,\ldots,n-1$, and $\K_{ij}=g_{ij,0}=0$. Also, we have
$g_{00,\beta}=0$ for $\beta=0,\ldots,n-1.$

To show $g_{11,\beta}=0$ we first consider the metric components
in MGS coordinates denoted by $g^{MGS}_{\alpha\beta}$.  In MGS
coordinates, $g^{MGS}_{11}=0$ throughout $\Si$, so
$g^{MGS}_{11,a}=0$ for $a=1,\ldots,n-1,$ and
$\K_{11}=g^{MGS}_{11,0}=0$. The metric component $g_{11}$ in
$u$-coordinates, since $x$ and $u$ coordinates describe the same
neighborhood of $p$, can be written in terms of metric components
in MGS coordinates as
$$g_{11}
=g^{MGS}_{00}-\frac{2}{\eta}\,g^{MGS}_{01}
+\frac{1}{\eta^2}\,g^{MGS}_{11}.$$ Therefore, since $g^{MGS}_{00,
\beta}=0$ and $g^{MGS}_{01,\beta}=0$ at $p$, then $g_{11,\beta}=0$
for all $\beta=0,\ldots,n-1$.
\end{proof}

\subsection{Justification of the Main Theorem, Theorem~\ref{thm:main}}
The justification of Theorem~\ref{thm:main} comes from the
statement of Theorems~\ref{lrh},~\ref{lcurvtens},~\ref{c11},
and~\ref{lorentz}.

\section{The Spherically Symmetric Case in Four Dimensional Spacetime}
% ---Spherical Symmetric Theorem---------------------------------
We conclude this chapter with another central result which
considers the case of matching two spherically symmetric metrics
across a hypersurface in four dimensional spacetime. The theorem
below will show that the weak form of conservation across the
shock surface, $\left[G^{\alpha\beta}\right]N_\alpha N_\beta=0$,
is implied by only a single condition when the areas of the
spheres of symmetry match smoothly at the shock surface and change
monotonically as the shock moves transversely to the areas of the
spheres of symmetry. The argument in the proof is formulated so
that this implication holds even when the shock surface is null.
Smoller and Temple have already proved this theorem for the
non-null case in~\cite{st94}, and  the proof, of the non-null
case, can also be found in~\cite{stg01} and~\cite{stgerm99}.
Therefore, we will prove the theorem here only in the case when
$\n$ is null. Our argument differs only slightly from Smoller and
Temple's in that we do not make use of the Israel condition
$$\left[G_{00}\right]=\left[\mbox{trace}(K^2)-(\mbox{trace}\,K)^2 \right],$$
but instead compute  $\left[G_{00}\right]$ directly. However, we
could prove the non-null case using the same argument below, and
replacing $\N$ with $\n$ and MGS coordinates with Gaussian normal
coordinates.

\begin{theorem}\label{thm:sphere}
Assume that $g$ and $\bar{g}$ are two spherically symmetric
metrics that match Lipschitz continuously across a
three-dimensional shock surface $\Sigma$ to form the matched
metric $g \cup \bar{g}$.  That is, assume that $g$ and $\bar{g}$
are Lorentzian metrics given by
\begin{equation}\label{metric}
  ds^2=-a(t,r)dt^2+b(t,r)dr^2+c(t,r)d\Omega^2
\end{equation}
 and
 \begin{equation}\label{metricbar}
  d\bar{s}^2 =-\bar{a}(\bar{t},\bar{r})d\bar{t}^2
+\bar{b}(\bar{t},\bar{r})d\bar{r}^2
+\bar{c}(\bar{t},\bar{r})d\Omega^2,
\end{equation}
where $d\Omega^2=d\theta^2 + \sin^2 \theta d\varphi^2$ is the
standard metric on the unit 2 sphere. Assume that there exists a
smooth coordinate transformation $\Psi : (t,r)\rightarrow
(\bar{t},\bar{r})$, defined in a neighborhood of the shock
surface $\Sigma$ given by $r=r(t)$, such that the metrics agree
on $\Sigma.$ (We implicitly assume that $\theta$ and $\varphi$
are continuous across the surface.) Assume that
\begin{equation}\label{eq:c}
c(t,r)=\bar{c}(\bar{t},\bar{r}),
\end{equation}
in an open neighborhood of the shock surface $\Sigma$, so that, in
particular, the area of the two-spheres of symmetry in the barred
and unbarred metrics agree on the shock surface.  Assume also that
the shock surface $r=r(t)$ in unbarred coordinates is mapped to
the surface $\bar{r}=\bar{r}(\bar{t})$  by
$(\bar{t},\bar{r}(\bar{t}))=\Psi(t,r(t))$. Let $\N(c)$ denote the
derivative of the function $c$ in the direction of the vector
$\N$, and assume that $\N(c)\neq 0$.

Then the following are equivalent to the statement that the
components of the metric $g \cup \bar{g}$ in any MGS coordinate
system are $C^{1,1}$ functions of these coordinates across the
surface $\Sigma:$

\begin{equation}\label{eq:s1}
  [G_{\alpha\beta}]N^\alpha=0
\end{equation}

\begin{equation}\label{eq:s2}
  [G_{\alpha\beta}]N^\alpha N^\beta=0,
\end{equation}
and
\begin{equation}\label{eq:s3}
  [\K]=0.
\end{equation}

\end{theorem}

\begin{proof} From Theorem~\ref{thm:main} we have that $[\K]=0$
at each point of $\Si$, if and only if, for each point $p$ in
$\Sigma$ there exists a $C^{1,1}$ coordinate transformation
defined in a neighborhood of $p,$ such that, in the new
coordinates the metric components are $C^{1,1}$ functions of these
coordinates.

The functions $c(t,r)$ and $\bar{c}(\bar{t},\bar{r})$, which are
defined in equations~\eqref{metric} and~\eqref{metricbar},
transform as functions under $(t,r)$-transformations. By
equation~\eqref{eq:c} we have $g_{22}=c=\bar{c}=\bar{g}_{22}$ on
$\Si$, and therefore
$g_{22,0}=\N(c)=\N(\bar{c})=\bar{g}_{22,0}\neq0$ and
$g_{33,0}=\N(c)\sin^2\theta
=\N(\bar{c})\sin^2\theta=\bar{g}_{33,0}\neq 0$ on $\Si.$
Furthermore, $g_{22,00}=\N'(c)=\N'(\bar{c})=\bar{g}_{22,00}$ and
$g_{33,00}=\N'(c)\sin^2\theta
=\N'(\bar{c})\sin^2\theta=\bar{g}_{33,00}$ on $\Si.$ Then we also
have from Theorem~\ref{thm:main} that since $[\K]=0$ and
$[\K_{tt,0}]=[g_{tt,00}]=0$ for $t=2,3$, it follows that
$[G_{\alpha\beta} ]N^\alpha=0$. Furthermore, if $[G_{\alpha\beta}
]N^\alpha=0,$ then $[G_{\alpha\beta}]N^\alpha N^\beta=0.$
Therefore, to complete the proof, it remains to show that
$[G_{\alpha\beta}]N^\alpha N^\beta=0$ implies that $[\K]=0$.

To this end, we choose a smooth coordinate system
 $(w^1,w^2,w^3)=(z^1,\theta,\varphi)$ where $\partial/\partial w^1
=\n$. Now extend the above coordinates to a MGS coordinate system
$w=(w^0,w^1,w^2,w^3).$ Notice that by the construction of a MGS
coordinate system, $\N$ is perpendicular to the coinciding
2-spheres of symmetry, that is, $\N$ depends only on the time and
radial components of any coordinate system. Furthermore, in MGS
coordinates, the metric $g\cup\bar{g}$ is diagonal except for
$(g\cup\bar{g})_{01}=\eta$ which is a constant. Since, by
Lemma~\ref{l3}, $\K_{ab}=-(1/2)(g\cup\bar{g})_{ab,0}$, we have
that $\K_{ab}$ is diagonal. Therefore, the only non-zero
components of $\K$ are
\begin{equation}
\K_{11}=-\frac{1}{2}(g\cup\bar{g})_{11,0},
\end{equation}
\begin{equation}
\K_{22}=-\frac{1}{2}(g\cup\bar{g})_{22,0},
\end{equation}
and
\begin{equation}
\K_{33}=-\frac{1}{2}(g\cup\bar{g})_{33,0}.
\end{equation}
Thus, we have
\begin{equation}\label{eq:k22}
\left[\K_{22}\right]=0,
\end{equation}
and
\begin{equation}\label{eq:k33}
\left[\K_{33}\right]=0,
\end{equation}
across $\Si$. Thus it only remains to show that
$\left[\K_{11}\right]=0$. To do this we must use our condition
$[G_{\alpha\beta}]N^\alpha N^\beta=0$

In MGS coordinates this is
\begin{equation}\label{eq:Gab}
\left[G_{\alpha\beta}\right]N^\alpha N^\beta
=\left[G_{00}\right]=g_{0\sigma}\left[G^\sigma_0\right]
=g_{00}\left[G^0_0\right]+g_{01}\left[G^1_0\right]
=\left[G^0_0\right]+\eta\left[G^1_0\right]=0.
\end{equation}
By equation~\eqref{eq:l11},
\begin{equation}\label{eq:G00a}\begin{aligned}
\left[G^0_0\right]=-\sum_{s,t\neq0}\left[R^{|st|}_{|st|}\right]
&=-\left[R^{12}_{12}\right]
-\left[R^{13}_{13}\right]-\left[R^{23}_{23}\right] \\
&=-g^{22}\left[R^1_{212}\right] -g^{33}\left[R^1_{313}\right]
-g^{33}\left[R^2_{323}\right]
\end{aligned}
\end{equation}
and by equation~\eqref{l51} we find that
$$\left[R^1_{t1t}\right]
=-\left[\K^1_1 \K_{tt}\right]
-\left[\widetilde{\G}^0_{tt}\K^1_{1}\right] =-g^{11}\left(\K_{tt}
+\widetilde{\G}^0_{tt}\right)\left[\K_{11}\right],$$ for $t=2,3$.
Also, from equation~\eqref{l51} it follows
$\left[R^2_{323}\right]=0$. Now since, for $t=2,3$,
$$\widetilde{\G}^0_{tt}=\frac{1}{2}g^{0 s}\left\{-g_{tt,s}
+g_{st,t} +g_{ts,t}\right\} =-\frac{1}{2}g^{01}g_{tt,1}
=-\frac{1}{2\eta}g_{tt,1},$$ where we have used that
$g^{01}=1/\eta$ and $g_{01}=\eta$ are the only non-zero off
diagonal components of $g$ and $g^{-1}$, and so we can say
\begin{equation}\label{eq:Ri1i}
\left[R^1_{t1t}\right]=-g^{1 1}\left(\K_{tt}
-\frac{1}{2\eta}g_{tt,1}\right)\left[\K_{11}\right].
\end{equation}
Therefore, by equations~\eqref{eq:G00a} and~\eqref{eq:Ri1i}, and
also using $g^{11}=-1/\eta^2$, $g_{22,1}=\n(c)$, and
$g_{33,1}=\n(c)\sin^2\theta$, we can write
\begin{equation}\label{eq:G00}
\begin{aligned}
\left[G^0_0\right] &=-g^{22}\left[R^1_{212}\right]
-g^{33}\left[R^1_{313}\right] \\
&=\frac{g^{22}}{\eta^2}\left(\K_{22}
-\frac{1}{2\eta}g_{22,1}\right)\left[\K_{11}\right]
+\frac{g^{33}}{\eta^2}\left(\K_{33}
-\frac{1}{2\eta}g_{33,1}\right)\left[\K_{11}\right]\\
&=\frac{1}{c\eta^2}\left(\N(c)
-\frac{1}{2\eta}\n(c)\right)\left[\K_{11}\right]
+\frac{1}{c\sin^2\theta\eta^2}\left(\N(c)\sin^2\theta
-\frac{1}{2\eta}\n(c)\sin^2\theta\right)\left[\K_{11}\right]\\
&=\frac{2}{c\eta^2}\left(\N(c)
-\frac{1}{2\eta}\n(c)\right)\left[\K_{11}\right].
\end{aligned}
\end{equation}
Using equation~\eqref{eq:l12}
\begin{equation}\label{eq:G10a}
\left[G^1_0\right]=\sum_{t\neq0,1}\left[R^{|1t|}_{|0t|}\right]
=\left[R^{12}_{02}\right]+\left[R^{13}_{03}\right]
=g^{22}\left[R^1_{202}\right] +g^{33}\left[R^1_{303}\right],
\end{equation}
where we have used $g_{ij}=0$ for $i\neq j$ and $i,j=2,3.$ By
equation~\eqref{l52}, for $t=2,3$,
\begin{equation}\label{eq:R101}
\left[R^1_{t0t} \right]=-\left[R^1_{tt0} \right]
=\frac{1}{\eta^2}\left[K_{11}\right]\left( \frac{1}{\eta}K_{tt}
+\frac{1}{2\eta^2}g_{tt,1}\right).
\end{equation}
Therefore,
\begin{equation}\label{eq:G10}
\eta\left[G^1_0\right] =\frac{2}{c\eta^2}\left(\N(c)
+\frac{1}{2\eta}\n(c)\right)\left[\K_{11}\right]
\end{equation}
Using equations~\eqref{eq:G00} and~\eqref{eq:G10} we can write
equation~\eqref{eq:Gab} as
\begin{equation}\label{eq:Gabf}
\left[G_{\alpha\beta}\right]N^\alpha N^\beta
=\frac{4\N(c)}{c\eta^2}\left[\K_{11}\right]=0.
\end{equation}
Since we assumed that $\N(c)\neq 0$, we must have that
$\left[\K_{11}\right]=0$, and the theorem is proved.
\end{proof}

% ---End Chapter 3------------------------------------------------

%% file: exam.tex
%---Example-------------------------------------------------------
%An Astrophysical Shock-Wave in the Lightlike Limit
\label{chp:4}

%\section{Introduction}
The goal of this chapter is to construct exact, spherically
symmetric, lightlike shock-wave solutions of the Einstein
equations. We will follow Smoller and Temple's work
from~\cite{st95} in which they matched a
Friedman-Robertson-Walker (FRW) metric to a
Tolman-Oppenheimer-Volkoff (TOV) metric across a timelike
hypersurface. We will use their same matching technique here to
show the existence of lightlike shock-wave solutions of the
Einstein equations.

We are modelling spherically symmetric lightlike shock-wave
expanding into an static spacetime. The geometry of the region
behind the shock will be that of our universe given by the FRW
metric, and the region in front of the shock will be static
spacetime whose geometry is given by the TOV metric.

In this chapter we first derive the exact FRW and TOV type
solutions of the Einstein equations, and then match these
solutions along the 2-spheres of symmetry to get an exact,
lightlike, shock-wave solution. As stated in the introduction we
will follow~\cite{st95}, but will use the prescription given in
the previous chapter to obtain our lightlike shock solution.

\section{The Metrics}
As we stated in chapter~\ref{background}, physical principles
determine our spacetime metric. The case is no different here.
Both the FRW and TOV metrics are derived from assumptions which
seem to concur with observational data.

\subsection{The Friedman-Roberson-Walker (FRW) Metric}
The FRW metric is generally regarded as the metric which models
our universe on a large scale, around $10^8$ to $10^9$ light
years~\cite{weinberg}. It is derived from the \emph{cosmological
principle}\index{cosmological principle} which are assumptions
that our universe is both \emph{homogeneous}\index{homogeneous}
and \emph{isotropic}\index{isotropic} at each point.
Homogeneous\index{homogeneous} means that, on a large enough
scale, in any given instant of time, each point of space looks
like any other.  Isotropic\index{isotropic} means that there are
no preferred directions in space; observations do not depend upon
which direction we look when we consider these large scales. We
will not derive metric here, but the derivation can be found in
most texts on the subject of general relativity.
See~\cite{hawkellis}, \cite{mtw}, \cite{schutz}, \cite{wald},
and~\cite{weinberg} just to name a few. The FRW metric is defined
as
\begin{equation}\label{frw}
ds^2=-dt^2+\frac{R^2(t)}{1-kr^2}dr^2+r^2R^2(t)d\Omega^2,
\end{equation}
where $d\Omega^2=d\theta^2+\sin^2\theta\,d\phi^2$ denotes the
standard metric on the unit 2-sphere. The constant $k$ can be
chosen to be $+1,$ $-1,$ or $0$ each giving the spatial geometry
of the 3-sphere, flat space, and hyperbolic space respectively.
The function $R(t)$ is sometimes referred to as the "cosmological
scale factor."

This FRW metric was first derived by H. P. Robertson and A. G.
Walker in the 1930's, but it was incomplete in that it did not
give a prediction for the function $R(t).$ In 1922 Alexandre
Friedmann had made some assumptions about the material content of
the universe, and if one derives FRW metric as a solution of the
Einstein equations with these assumptions, then $R(t)$ can be
computed~\cite{weinberg}.

\subsection{The Tolman-Oppenheimer-Volkov (TOV) Metric} The TOV
metric describes the geometry in front of our shock surface which
models a static and isotropic gravitational field. An especially
good derivation of this metric can be found in~\cite[chapter
11]{weinberg}. The TOV metric is given by
\begin{equation}\label{tov}
d\bar{s}^2=-B(\bar{r})d\bar{t}^2+A(\bar{r})^{-1}d\bar{r}^2+\bar{r}^2d\Omega^2.
\end{equation}
The TOV metric given in barred coordinates to distinguish it from
the unbarred coordinates of the FRW metric.

\section{The FRW and TOV Solutions}

In this section we derive the exact FRW and TOV type solutions of
the Einstein equations which will encompass the regions of
spacetime behind and in front of the shock surface.

We assume matter in spacetime is modelled by a perfect fluid that
is comoving relative to the coordinates, that is, in free fall
relative to the coordinates~\cite[section 11.8]{weinberg}. Then
the four-velocity of the fluid is given by
\begin{equation}\label{eq:comov}
u^0=\sqrt{-g_{00}}, \quad \mbox{and} \quad u^i=0, \; i=1,2,3.
\end{equation}
Recall that (see equation~\eqref{eq:stren}) for a perfect fluid
the stress-energy tensor becomes
\begin{equation}\label{eq:stren1}
T^{\alpha\beta}=pg^{\alpha\beta}+(p+\rho)u^\alpha u^\beta,
\hspace{0.25in} \alpha,\beta=0,\ldots, 3.
\end{equation}

\subsection{The FRW Solution}
The FRW solution of the Einstein equations will be the solution
behind the shock. Substituting the FRW metric~\eqref{frw} into
the Einstein field equations
$$G=\kappa T,$$
where $T$ is given by~\eqref{eq:stren1}, and assuming an equation
of state of the form $p=p(\rho)$ yields the following pair of
differential equations
\begin{equation}\label{eq:ode1}
p=-\rho-\frac{R\dot{\rho}}{3\dot{R}},
\end{equation}
and
\begin{equation}\label{eq:ode2}
\dot{R}^2+k=\frac{8\pi\mathcal{G}}{3}\rho R^2,
\end{equation}
from which we can solve for the unknown functions $R(t)$, and
$\rho(t)$~\cite{st94, st95, stgerm99}.

As a simplification we restrict ourselves to the case $k=0$ so
that the FRW metric in~\eqref{frw} becomes
$$ds^2=-dt^2+R^2(t)dr^2+r^2R^2(t)d\Omega^2,$$
and the differential equations~\eqref{eq:ode1}
and~\eqref{eq:ode2} become
\begin{equation}\label{eq:ode1s}
p=-\rho-\frac{R\dot{\rho}}{3\dot{R}},
\end{equation}
and
\begin{equation}\label{eq:ode2s}
\dot{R}^2=\frac{8\pi\mathcal{G}}{3}\rho R^2.
\end{equation}
Rewriting the equation~\eqref{eq:ode2s} as
\begin{equation}\label{eq:ode2ss}
  \dot{R}=\pm\left(\frac{8\pi\mathcal{G}\rho}{3}\right)^{1/2} R,
\end{equation}
and then substituting into equation~\eqref{eq:ode1s} we get
\begin{equation}\label{eq:ode1ss}
  p=-\rho\mp\frac{\dot{\rho}}{\sqrt{24\pi\mathcal{G}\rho}}.
\end{equation}
Notice that the $\pm$ signs in equation~\eqref{eq:ode2ss}
correspond directly to the $\mp$ signs in
equation~\eqref{eq:ode1ss}. We can solve~\eqref{eq:ode1ss}
explicitly when the equation of state $p=p(\rho)$ is given.
Indeed, solving~\eqref{eq:ode1ss} for $dt$ yields
\begin{equation}\label{eq:frwdt}
dt=\mp\frac{d\rho}{(\rho + p)\sqrt{24\pi \mathcal{G}\rho}},
\end{equation}
which we integrate to get
\begin{equation}\label{eq:frw_t}
t-t_0 =\mp\int_{\rho_0}^{\rho}\frac{d\xi}{[\xi +
p(\xi)]\sqrt{24\pi \mathcal{G}\xi}}.
\end{equation}
Then using by equation~\eqref{eq:ode2ss}, and
\begin{equation}\label{eq:Rdot}
\dot{R}=\frac{d\rho}{dt}\frac{dR}{d\rho} =\mp(\rho +
p)\sqrt{24\pi \mathcal{G}\rho}\frac{dR}{d\rho}
\end{equation}
we have
\begin{equation}\label{eq:dRR}
  \frac{dR}{R}=\frac{-d\rho}{3(\rho+p)}.
\end{equation}
It follows that equation~\eqref{eq:dRR} has solution
\begin{equation}\label{frwsln}
R=R_0\exp\left(\int_{\rho_0}^{\rho}\frac{-1}{3\left(\xi +
p(\xi)\right)}\,d\xi\right).
\end{equation}

\subsection{The TOV Solution}
Now we give the the TOV solution which encompasses the region in
front of the shock, and represents a general relativistic version
of static, singular isothermal sphere~\cite{st95}.

We proceed in a similar manner at the FRW metric by substituting
the TOV metric~\eqref{tov} into the Einstein
equations\index{Einstein equations}
$$G=\kappa T,$$ where $T$ is given by~\eqref{eq:stren1}.
In this case, the substitution (see~\cite{weinberg}) gives the
following

\begin{equation}\label{eq:A}
  A(\bar{r})=\left(1-\frac{2\mathcal{G}M}{\bar{r}}\right),
\end{equation}
\begin{equation}\label{eq:dmdr}
  \frac{dM}{d\bar{r}}=4\pi \bar{r}^2\bar{\rho},
\end{equation}
and
\begin{equation}\label{eq:ove}
  -\bar{r}^2\frac{d}{d\bar{r}}\bar{p}
  =\mathcal{G}M\bar{\rho}\left\{1+\frac{\bar{p}}{\bar{\rho}}\right\}
  \left\{1+\frac{4\pi \bar{r}^3\bar{p}}{M}\right\}
  \left\{1-\frac{2\mathcal{G}M}{\bar{r}}\right\}^{-1}.
\end{equation}
Here the unknown functions, $\bar{\rho}(\bar{r}),$
$\bar{p}(\bar{r})$, and $M(\bar{r})$, depend only on $\bar{r}$.
Analogous to the FRW solution we assume an equation of state
$p=p(\rho)$ for the TOV metric so that the differential
equations~\eqref{eq:dmdr} and~\eqref{eq:ove} give solutions for
$M(\bar{r})$ and $\bar{\rho}(\bar{r}).$ The function $M(\bar{r})$
denotes the mass inside radius $\bar{r}$, and can be written as
\begin{equation}\label{eq:M}
  M(\bar{r})=\int_0^{\bar{r}}4\pi \xi^2\bar{\rho}(\xi)\,d\xi.
\end{equation}
To find the metric component $B(\bar{r})$ we look to the equation
for hydrostatic equilibrium~\cite[equation 11.1.8]{weinberg},
\begin{equation}\label{eq:dBB}
  \frac{B'(\bar{r})}{B}=-2\frac{\bar{p}'(\bar{r})}{\bar{p}+\bar{\rho}}.
\end{equation}

Now we further restrict the equation of state for the TOV metric
to be of the form
\begin{equation}\label{eq:pbar}
 \bar{p}=\bar{\sigma}\bar{\rho},
\end{equation}
where $\bar{\sigma}$ is a constant, and also assume the energy
density is given by
\begin{equation}\label{eq:rhobar}
 \bar{\rho}=\frac{\gamma}{\bar{r}^2},
\end{equation}
with $\gamma$ a constant. With these assumptions
equation~\eqref{eq:M} becomes
\begin{equation}\label{eq:Ms}
 M(\bar{r})
 =\int_0^{\bar{r}}4\pi\xi^2\bar{\rho}(\xi)\,d\xi=4\pi \gamma\bar{r}.
\end{equation}
Combining equations~\eqref{eq:pbar} -~\eqref{eq:Ms}
with~\eqref{eq:ove} gives
\begin{equation}\label{eq:gamma}
  \gamma =\frac{1}{2\pi\mathcal{G}}
  \left(\frac{\bar{\sigma}}{1+6\bar{\sigma}+\bar{\sigma}^2}\right),
\end{equation}
and inserting equation ~\eqref{eq:Ms} into equation~\eqref{eq:A}
yields
\begin{equation}\label{eq:As}
  A=1-8\pi\mathcal{G}\gamma.
\end{equation}
Now, we find $B$ by inserting equations~\eqref{eq:pbar}
and~\eqref{eq:rhobar} into equation~\eqref{eq:dBB}, which becomes
\begin{equation}\label{eq:dBBS}
  \frac{dB}{B}=-\frac{2\bar{\sigma}}{(1+\bar{\sigma})}\frac{d\bar{\rho}}{\bar{\rho}}.
\end{equation}
Then solving we get
\begin{equation}\label{eq:B}
  B=B_0\left(\frac{\bar{r}}{\bar{r}_0}\right)^{4\bar{\sigma}/(1+\bar{\sigma})}.
\end{equation}

\section{Matching the FRW metric to the TOV metric}
In~\cite{st94} Smoller and Temple derive a coordinate
transformation that takes $(\bar{t},\bar{r}) \to (t,r)$ so that
the FRW metric~\eqref{frw} and the TOV metric~\eqref{tov} match
Lipschitz continuously across a shock surface $\Si.$ The same
derivation can also be found in~\cite{stg01} and~\cite{stgerm99}.
Here we summarize their procedure here as briefly as possible
highlighting their results which we will use to construct a
lightlike shock-wave solution.

Since our goal is to be able to apply Theorem~\ref{thm:sphere} we
begin by letting
\begin{equation}\label{eq:rRr}
  \bar{r}(t,r)=R(t)r,
\end{equation}
which ensures that
$$\bar{r}^2\,d\Omega^2=R^2r^2\,d\Omega^2.$$
A consequence of equation~\eqref{eq:rRr} is that we can write the
following
\begin{align}
d\bar{r}&=Rdr+\dot{R}rdt, \label{eq:drbar}\\
dr&=\frac{1}{R}d\bar{r}-\frac{\dot{R}}{R}rdt, \label{eq:dr}\\
\dot{r}&=\frac{\dot{\bar{r}}}{R}-\frac{\dot{R}r}{R}
\label{eq:rdot}.
\end{align}
Now, it can be shown that in $(t,\bar{r})$-coordinates the FRW
metric~\eqref{frw} is written as
\begin{equation}\label{frw_bar}
ds^2=\frac{1}{R^2-k\bar{r}^2}
\left\{-R^2(1-\frac{8\pi{\mathcal{G}}}{3}\rho R^2r^2)dt^2
+R^2d\bar{r}^2-2R\dot{R}\bar{r}\,dt\,d\bar{r}\right\}
+\bar{r}^2d\Omega^2.
\end{equation}
At this point a mapping $t=t(\bar{t},\bar{r})$ is constructed to
eliminate the cross term $dt\,d\bar{r}$ in
equation~\eqref{frw_bar}. This is done first for a general metric
given by
\begin{equation}\label{eq:gen_met}
d\tilde{s}^2=-C(t,\bar{r})dt^2+D(t,\bar{r})d\bar{r}^2+2E(t,\bar{r})dtd\bar{r},
\end{equation}
and choosing $\psi=\psi(t,\bar{r})$ such that
\begin{equation}\label{eq:psi}
\frac{\partial}{\partial \bar{r}}(\psi C)
=-\frac{\partial}{\partial t}(\psi E),
\end{equation}
so that
\begin{equation}\label{eq:dpsi}
d\bar{t}=\psi(t,\bar{r})\{C(t,\bar{r})dt-E(t,\bar{r})d\bar{r}\},
\end{equation}
is an exact differential. This means that
equation~\eqref{eq:gen_met} can be written as
\begin{equation}\label{eq:gen_mets}
d\tilde{s}^2 =-(\psi^{-2}C^{-1})d\bar{t}^2 +\left(D
+\frac{E^2}{C}\right)d\bar{r}^2.
\end{equation}
Therefore, the FRW metric in $(\bar{t},\bar{r})$-coordinates is
\begin{equation}\label{frw_bars}
ds^2=\frac{1}{R^2-k\bar{r}^2} \left\{-(\psi^2C)^{-1}d\bar{t}^2
+\left(D+\frac{E^2}{C}\right)
d\bar{r}^2\right\}+\bar{r}^2d\Omega^2,
\end{equation}
where
\begin{equation}\label{eq:CDE}
C=R^2\{1-\frac{8\pi{\mathcal{G}}}{3}\rho R^2r^2\}, \quad D=R^2,
\quad \mbox{and} \quad E=-R\dot{R}\bar{r}.
\end{equation}

Now we can finally define the shock surface at which the FRW
metric~\eqref{frw} and the TOV metric~\eqref{tov} match Lipschitz
continuously. Indeed, using equation~\eqref{eq:CDE} write
\begin{equation}\label{eq:DE2C}
D+\frac{E^2}{C} =R^2+\frac{R^2\dot{R}^2\bar{r}^2}{R^2
\left(1-\frac{8\pi{\mathcal{G}}}{3}\rho\bar{r}^2\right)}
=R^2+\frac{\dot{R}^2R^2r^2}{1-\frac{8\pi{\mathcal{G}}}{3}\rho
R^2r^2},
\end{equation}
and equate the $d\bar{r}^2$ components in the TOV
metric~\eqref{tov} and the FRW metric in
$(\bar{t},\bar{r})$-coordinates~\eqref{frw_bars}. Then it follows
that
\begin{equation}\label{eq:shcks0}
\left(R^2-k\bar{r}^2\right)
\left(1-\frac{2{\mathcal{G}}M}{\bar{r}}\right)^{-1}
=R^2+\frac{\dot{R}^2R^2r^2}{1-\frac{8\pi{\mathcal{G}}}{3}\rho
R^2r^2},
\end{equation}
and using the differential equation~\eqref{eq:ode2} we get
\begin{equation}\label{eq:shock}
M(\bar{r})=\frac{4\pi}{3}\rho(t)\bar{r}^3,
\end{equation}
which implicitly defines the shock surface $\Si$. To express
$\Si$ in $(t,r)$-coordinates we use the transformation
$\bar{r}=R(t)r.$ For the shock surface $\Si$ to remain in the
domain of definition of the FRW metric we must assume that
$1-kr^2>0$ for $k>0$~\cite{st94}. We also must obtain the
conditions under which the function $\psi$ defined in
equation~\eqref{eq:psi}, which determines $t$ from
$(\bar{t},\bar{r})$, can be found uniquely. Since the
$d\bar{t}^2$ terms must match on the shock surface, it follows
that
\begin{equation}\label{eq:initpsi}
 \frac{1}{R^2-k\bar{r}^2}\frac{1}{\psi^2 C}=B(\bar{r})
\end{equation}
must hold on the shock surface. Then $\psi$ can be determined by
the partial differential equation
\begin{equation}\label{eq:psipde}
 C\psi_{\bar{r}}+E\psi_{t}=f(t,\bar{r},\psi),
\end{equation}
where $C(t,\bar{r})$ and $E(t,\bar{r})$ are given in
equation~\eqref{eq:CDE}. As shown in~\cite{stg01},~\cite{st94},
and~\cite{stgerm99}, the partial differential
equation~\eqref{eq:psipde} with initial value given
by~\eqref{eq:initpsi} can be solved uniquely in neighborhood of a
point on the shock surface provided
\begin{equation}\label{eq:nonchar}
 \dot{\bar{r}}=\frac{d\bar{r}}{dt}\neq\frac{C}{E}.
\end{equation}
Here $\dot{\bar{r}}$ denotes the speed of the shock surface, and
equation~\eqref{eq:nonchar} is the condition that the shock
surface be non-characteristic at a point~\cite{st94}.

Now that we have the shock surface, we restate another result of
Smoller and Temple, proposition 1 in~\cite{stg01},~\cite{st94},
and~\cite{stgerm99}, which will be used for computations
involving the transformation from $(t,r)$-coordinates to
$(\bar{t},\bar{r})$-coordinates and back.
\begin{proposition}\label{p1}
On the shock surface given by equation~\eqref{eq:shock}, the
following identities hold:
\begin{align}
\frac{1}{\psi^2C^2}&=B\left(1+\frac{AE^2}{C^2}\right)=\frac{B}{A}(1-kr^2),\label{eq:p11}\\
C&=R^2A,\label{eq:p12}\\
\frac{E}{C}&=\frac{-\dot{R}r}{A},\label{eq:p13}\\
\frac{E^2}{C^2}&=\frac{-A+(1-kr^2)}{A^2}, \label{eq:p14}\\
\dot{R}^2r^2&=-A+(1-kr^2)\label{eq:p15}.
\end{align}
\end{proposition}

In addition to using proposition~\ref{p1} to transform quantities
in $(t,r)$-coordinates to and from
$(\bar{t},\bar{r})$-coordinates we will also need the following
proposition.

\begin{proposition}\label{p2}
On the shock surface given by~\eqref{eq:shock}, the following
partial derivatives can be written as
\begin{alignat}{2}
\frac{\partial t}{\partial \bar{t}}&=(\psi C)^{-1} ,&\qquad
\frac{\partial t}{\partial \bar{r}}&= \frac{E}{C},\\
\frac{\partial r}{\partial \bar{t}}&= \frac{A}{R}\frac{E}{C}(\psi
C)^{-1}, &\qquad \frac{\partial r}{\partial \bar{r}}&=
\frac{1-kr^2}{RA}.
\end{alignat}
\end{proposition}

\begin{proof}
From equation~\eqref{eq:dpsi} we can write
\begin{equation}\label{eq:dt}
dt=(\psi C)^{-1}d\bar{t}+\frac{E}{C}d\bar{r},
\end{equation}
from which, follows
\begin{equation}\label{eq:dtdtbar}
\frac{\partial t}{\partial \bar{t}}=(\psi C)^{-1},
\end{equation}
and
$$\frac{\partial t}{\partial \bar{r}}= \frac{E}{C}.$$ By
equation~\eqref{eq:dr} we have
$$dr=\frac{1}{R}d\bar{r}-\frac{\dot{R}r}{R}\frac{dt}{d\bar{t}}\,d\bar{t},$$
and so, using equation~\eqref{eq:p13},
$$\frac{\partial r}{\partial \bar{t}} =\frac{\dot{Rr}}{R}\frac{dt}{d\bar{t}}
=\frac{A}{R}\frac{E}{C}(\psi C)^{-1}.$$ Lastly, using
equations~\eqref{eq:rRr},~\eqref{eq:p13},~\eqref{eq:p14}, and
~\eqref{eq:dtdtbar} yields
\begin{multline*}
\frac{\partial r}{\partial \bar{r}} =\frac{\partial}{\partial
\bar{r}} \left(\frac{\bar{r}}{R}\right)
=\frac{1}{R}-\frac{\dot{R}\bar{r}}{R^2} \frac{\partial
t}{\partial \bar{t}} =\frac{1}{R}-\frac{\dot{R}r}{R} (\psi
C)^{-1} \\
=\frac{1}{R}+\left(A\frac{E}{C}\right)\frac{E}{C}\frac{1}{R}
=\frac{1}{R}\left(1-\frac{A+(1-kr^2)}{A}\right)
=\frac{1-kr^2}{RA}.
\end{multline*}
\end{proof}

We summarize the conditions under which we can match the FRW and
TOV metrics, given in~\eqref{frw} and~\eqref{tov} respectively,
in the following theorem, which is a restatement of theorem 6
in~\cite{stg01}, and theorem 5 in~\cite{st94} and~\cite{stgerm99}.

\begin{theorem}\label{thm:surf}
Let the point $(t_0,\bar{r}_0)$ satisfy
\begin{equation}\label{eq:surf1}
M(\bar{r})=\frac{4\pi}{3}\rho(t)\bar{r}^3,
\end{equation}
and let~\eqref{eq:surf1} define the shock surface
$\bar{r}=\bar{r}(t)$ in a neighborhood of $(t_0,\bar{r}_0)$.
Furthermore, assume
\begin{equation}\label{eq:surf2}
R(t)r=\bar{r},
\end{equation}
requiring the barred and unbarred coordinates to be equivalent on
the spheres of symmetry of the FRW and TOV metrics, given
in~\eqref{frw} and~\eqref{tov}, and the shock surface in
$(t,r)$-coordinates be given by $r=r(t)=\bar{r}(t)/R(t)$. Then,
when
$$1-kr^2>0,$$
$$A(\bar{r}_0)\neq 0,$$
and
$$\dot{\bar{r}}=\frac{d\bar{r}}{dt}\neq\frac{C}{E}=-\frac{\dot{R}r}{A},$$
hold at $t=t_0$, we can define the coordinate $\bar{t}$ so that
there exists a smooth regular transformation on a neighborhood of
$(t_0,r_0)$ that takes
$$(t,r)\To (\bar{t},\bar{r}),$$
and the FRW and TOV metrics, given in~\eqref{frw}
and~\eqref{tov}, match Lipschitz continuously across the shock
surface $r=r(t).$
\end{theorem}

\section{The Conservation Condition}

The goal of this section is write down a formula that expresses
the weak form of conservation of energy across a spherically
symmetric lightlike shock surface constructed in the manner of the
previous section. We should note that, so far, in this chapter we
have merely recapitulated results from~\cite{stg01},~\cite{st94},
and~\cite{stgerm99}, that is, we have not stated anything that
has to do with our shock surface being lightlike or not. It is in
this part of this chapter that our results begin to diverge from
that of Smoller and Temple in the sense that we are dealing with
the lightlike case. However, we are still using their work as a
model to obtain our own results.

\subsection{The Lightlike Surface}
Assume that equation~\eqref{eq:shock},
\begin{equation}\label{eq:shock1}
 M(\bar{r})=\frac{4\pi}{3}\rho(t)\bar{r}^3,
\end{equation}
defines the surface given by the level curve
\begin{equation}\label{eq:phi}
 \varphi(t,r)=r-r(t)=0
\end{equation}
in a neighborhood of $(\bar{t}_0,\bar{r}_0)$ which satisfies
equation~\eqref{eq:shock1}. We can also write $\varphi$ in
$(\bar{t},\bar{r})$-coordinates as
\begin{equation}\label{eq:phibar}
\varphi(\bar{t},\bar{r})
=\frac{\bar{r}}{R\left(t(\bar{t},\bar{r})\right)}
-r\left(t(\bar{t},\bar{r})\right).
\end{equation}
Then we can compute the normal to $\varphi=0$ by
\begin{equation}\label{eq:dphi}
d\varphi=n_0\,dt+n_1\, dr = -\dot{r}\,dt+dr
\end{equation}
to find
\begin{equation}\label{eq:n}
n_0=-\dot{r},\quad \mbox{and} \quad n_1=1.
\end{equation}
Although, the computations are a little more complicated, from
$\varphi$ in equation~\eqref{eq:phibar} we can write
$$d\varphi=\bar{n}_0\,d\bar{t}+\bar{n}_1\, d\bar{r}
=-\frac{\dot{\bar{r}}}{R(\psi C)}\,d\bar{t}
+\frac{\dot{\bar{r}}\dot{R}\bar{r}}{R^2 A}\, d\bar{r}$$ which
yields
\begin{equation}\label{eq:nbar}
\bar{n}_0=-\frac{\dot{\bar{r}}}{R(\psi C)},\quad \mbox{and} \quad
\bar{n}_1=\frac{1}{R}+\frac{\dot{\bar{r}}\dot{R}\bar{r}}{R^2 A}.
\end{equation}
We are assuming our surface is lightlike, hence
\begin{equation}\label{eq:n2}
 \langle \n,\n \rangle
 =g^{\alpha\beta}n_\alpha n_\beta
 =g^{00}(n_0)^2+g^{11}(n_1)^2 =-\dot{r}^2+\frac{1-kr^2}{R^2}=0,
\end{equation}
which implies that
\begin{equation}\label{eq:ndotn}
\dot{r}^2=\frac{1-kr^2}{R^2}.
\end{equation}
Notice that if we compute
$$\langle \bar{\n},\bar{\n }\rangle
=\bar{g}^{00}(\bar{n}_0)^2+\bar{g}^{11}(\bar{n}_1)^2
$$
using equation~\eqref{eq:nbar}, and the transform the result into
$(t,r)$-coordinates we get equation~\eqref{eq:n2}.

\subsection{The Transverse Vector}
Equation~\eqref{eq:n2} implies that $\n$ lies in that tangent
space of the surface, and therefore we must choose a transverse
vector $\N$ which satisfies equations~\eqref{eq:njump},
~\eqref{eq:njump1}, and~\eqref{eq:Nneta}. Now for any $\N$ we have
\begin{equation}\label{eq:Nchce}
 \langle\N,\n \rangle
 =g_{00}N^0n^0+g_{11}N^{1}n^{1}
 =N^0n_0+N^{1}n_1=-\dot{r}N^0+N^1=\eta\neq 0,
\end{equation}
In light of equation~\eqref{eq:Nchce}, we choose $\N$ so that
\begin{equation}\label{eq:N}
N^0=0,\quad\mbox{and}\quad N^1=\eta\neq 0,
\end{equation}
or, lowering indices,
\begin{equation}\label{eq:Nl}
N_0=0,\quad\mbox{and}\quad N_1=g_{11} N^1=\eta
\frac{R^2}{1-kr^2}.
\end{equation}
\begin{lemma}\label{l:N}
The transverse vector $\N$ defined in equation~\eqref{eq:N}
satisfies
\begin{equation}\label{eq:Njumps}
 \left[\left\langle\N, X_a\right\rangle\right]=0,
\end{equation} where $\{X_a\}_{a=1}^{3}$ is a basis
for the tangent of the surface at $(\bar{t}_0,\bar{r}_0)$, and
\begin{equation}\label{eq:Nnetas}
\langle\bar{\N},\bar{\n} \rangle
=\langle\N,\n \rangle=\eta\neq 0,
\end{equation}
where $\langle\bar{\N},\bar{\n} \rangle$ is written in
$(t,r)$-coordinates. Furthermore, we have that
\begin{equation} \label{eq:Njump1s}
 \left[\langle\N,\N\rangle\right]=0.
\end{equation}
\end{lemma}

\begin{proof}
To show equation~\eqref{eq:Nnetas} we first transform $N_0$ to
$\bar{N}_0$, and $N_1$ to $\bar{N}_1$. Using
propositions~\ref{p1} and~\ref{p2} we have
\begin{equation}\label{eq:Nbar0}
\bar{N}_0=N_0\,\frac{\partial t}{\partial \bar{t}}
+N_1\,\frac{\partial r}{\partial \bar{t}}
=\eta\frac{R^2}{1-kr^2}\frac{A}{R}\frac{E}{C}(\psi C)^{-1} =\eta
R\frac{A}{1-kr^2}\frac{E}{C}(\psi C)^{-1},
\end{equation}
and
\begin{equation}\label{eq:Nbar1}
\bar{N}_1=N_0\,\frac{\partial t}{\partial \bar{r}}
+N_1\,\frac{\partial r}{\partial \bar{r}}
=\eta\frac{R^2}{1-kr^2}\frac{1-kr^2}{RA} =\eta\frac{R}{A}.
\end{equation}
Then, using
equations~\eqref{tov},~\eqref{eq:p11},~\eqref{eq:p13},
and~\eqref{eq:nbar} we find
\begin{equation}\label{eq:Nn}
\begin{aligned}
 \left\langle\bar{\N}, \bar{\n}\right\rangle
 &=\bar{g}^{00}\bar{N}_0 \bar{n}_0 +\bar{g}^{11}\bar{N}_1 \bar{n}_1\\
 &= \left(-\frac{1}{B}\right) \eta R\frac{A}{1-kr^2}\frac{E}{C}(\psi C)^{-1}
 \left(\frac{-\dot{\bar{r}}}{R(\psi C)}\right)
+A\eta\frac{R}{A}\left(\frac{1}{R}
+\frac{\dot{\bar{r}}\dot{R}\bar{r}}{R^2 A}\right)\\
&=\eta\left(\frac{B^{-1}A\dot{\bar{r}}}{1-kr^2}\frac{E}{C}\frac{B(1-kr^2)}{A}
+1+\frac{\dot{\bar{r}}\dot{R}r}{A}\right)\\
&=\eta\left(-\frac{\dot{\bar{r}}\dot{R}\bar{r}}{A}+1
+\frac{\dot{\bar{r}}\dot{R}\bar{r}}{A} \right)=\eta.
\end{aligned}
\end{equation}
Then equations~\eqref{eq:Nchce},~\eqref{eq:N},and~\eqref{eq:Nn}
prove ~\eqref{eq:Nnetas} that holds.

For equation~\eqref{eq:Njumps} we must choose a basis for the
tangent of the surface at $(\bar{t}_0,\bar{r}_0).$ Let
$$X_1=\n, \ X_2=\frac{\partial}{\partial \theta},\ \mbox{and} \
X_3=\frac{\partial}{\partial \phi}.$$ Then by
equations~\eqref{eq:Nchce}, and~\eqref{eq:Nn} we have
$$\left[\left\langle\N, X_1\right\rangle\right]
=\left\langle\N, \n\right\rangle -\left\langle\bar{\N},
\bar{\n}\right\rangle=\eta-\eta=0.$$ Now,
$$\left[\left\langle\N, X_a\right\rangle\right]=0$$
for $a=2,3$ since
$$
\left\langle\N,\frac{\partial}{\partial \theta}\right\rangle
=\left\langle\bar{\N},\frac{\partial}{\partial
\theta}\right\rangle =\left\langle\N,\frac{\partial}{\partial
\phi}\right\rangle =\left\langle\bar{\N},\frac{\partial}{\partial
\phi}\right\rangle =0.$$

To show equation~\eqref{eq:Njump1s} holds we compute
$\langle\bar{\N},\bar{\N}\rangle$ in $(t,r)$-coordinates. From
equations~\eqref{eq:Nbar0} and~\eqref{eq:Nbar1}, and the
relations in proposition~\ref{p1} we see that
\begin{equation}\label{eq:NNbar}
\begin{aligned}
\langle\bar{\N},\bar{\N}\rangle &=\bar{g}^{00}\bar{N}_0 \bar{N}_0
+\bar{g}^{11}\bar{N}_1 \bar{N}_1 \\ &=-\frac{1}{B}\left(\eta
R\frac{A}{1-kr^2}\frac{E}{C}(\psi C)^{-1}\right)^2
+A\left(\eta\frac{R}{A}\right)^2\\
&=\eta^2\left(-\frac{1}{B}\frac{R^2
A^2}{(1-kr^2)^2}\frac{(-A+1-kr^2)}{A^2}\frac{B(1-kr^2)}{A}
+\frac{R^2}{A}\right)\\
&=\eta^2\left(\frac{R^2}{1-kr^2}-\frac{R^2}{A}+\frac{R^2}{A}\right)\\
&=\eta^2\frac{R^2}{1-kr^2}.
\end{aligned}
\end{equation}
We also have
\begin{equation}\label{eq:NN}
\langle\N,\N\rangle =g^{00}N_0 N_0 +g^{11}N_1 N_1
=\frac{1-kr^2}{R^2}\left(\eta\frac{R^2}{1-kr^2}\right)^2
=\eta^2\frac{R^2}{1-kr^2},
\end{equation}
and thus
$$\left[\langle\N,\N\rangle\right]
=\langle\N,\N\rangle-\langle\bar{\N},\bar{\N}\rangle
=\eta^2\frac{R^2}{1-kr^2}-\eta^2\frac{R^2}{1-kr^2}=0.$$
\end{proof}

\subsection{The Conservation Condition}

The conservation of energy across the surface defined by
equation~\eqref{eq:shock} is given by
\begin{equation}\label{eq:RH}
\left[T^{\alpha\beta}N_\alpha N_\beta\right]=0,
\end{equation}
where $T^{\alpha\beta}$ is the stress energy tensor for a perfect
fluid given in equation~\eqref{eq:stren1}. For comoving
coordinates, see equation~\eqref{eq:comov}, the jump
condition~\eqref{eq:RH}, using equations~\eqref{eq:Nl}
and~\eqref{eq:Nbar0}, becomes
\begin{equation}\label{eq:strenj}
\begin{aligned}
\left[T^{\alpha\beta}N_\alpha N_\beta\right]
&=(p-\bar{p})|\N|^2+(\rho+p)N_0^2-(\bar{\rho}+\bar{p})\frac{\bar{N}^2_0}{B}
\\ &=(p-\bar{p})\eta^2\frac{R^2}{1-kr^2}-(\bar{\rho}+\bar{p})
\eta^2 R^2\left(\frac{1}{A}-\frac{1}{1-kr^2}\right)\\
&= \eta^2R^2\left(\frac{p+\bar{\rho}}{1-kr^2}-\frac{\bar{p}
+\bar{\rho}}{A} \right)=0,
\end{aligned}
\end{equation}
where we have used
$$\begin{aligned}
\bar{\N}_0^2&=\left(\eta R\frac{A}{1-kr^2}\frac{E}{C}(\psi
C)^{-1}\right)^2\\
&=\eta^2 \frac{R^2A^2}{(1-kr^2)^2}
\frac{(-A+1-kr^2)}{A^2}\frac{B(1-kr^2)}{A} \\
&=\eta^2 B R^2 \left(\frac{1}{1-kr^2}-\frac{1}{A}\right).
\end{aligned}$$

When the jump condition~\eqref{eq:strenj} holds, it follows that
on solutions of the Einstein equations $G=\kappa T$, this is
equivalent to the jump condition on the Einstein tensor
\begin{equation}\label{eq:einj}
\left[G^{\alpha\beta}N_\alpha N_\beta\right]=0.
\end{equation}
Therefore, all the equivalencies in Theorem~\ref{thm:sphere} hold
which implies all of the equivalencies hold in
Theorem~\ref{thm:main}. Thus, there exists a lightlike shock-wave
solution of the Einstein equations when Theorem~\ref{thm:surf} and
equation~\eqref{eq:strenj} holds. In the next section we give an
exact such solution.

\section{An Exact Lightlike Shock-Wave Solution of the Einstein Equations}

\subsection{The jump Conditions}
\begin{figure}
  \centering
  \includegraphics[width=.6\textwidth]{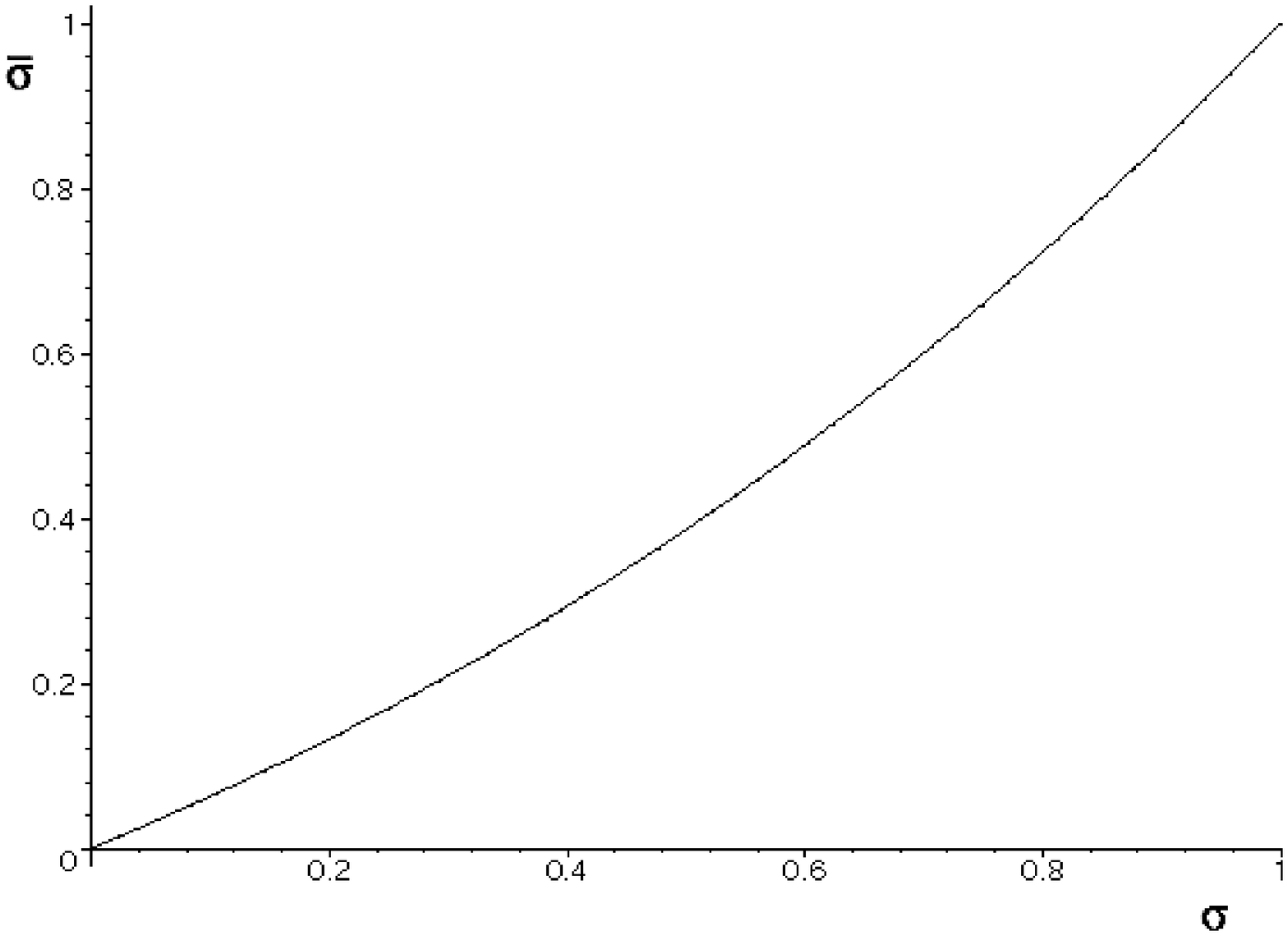}
  \caption{A graph of $\bar{\sigma}=H(\sigma)$
  defined in equation~\eqref{eq:H}. Note that
  $H(1)=1$ and $H(0)=0$.}
  \label{fig:sigvssigbar}
\end{figure}

We begin by assuming the suppositions of Theorem~\ref{thm:surf}
are satisfied, the jump condition~\eqref{eq:strenj} holds, $k=0$,
and the equation of state for the TOV metric is given by
$$\bar{p}=\bar{\sigma}\bar{\rho},$$
for some constant $\bar{\sigma}$. Furthermore, also suppose that
the TOV solution is given by~\eqref{eq:rhobar} -~\eqref{eq:As},
and~\eqref{eq:B}. On the shock surface given by
$$M(\bar{r})=\frac{4\pi}{3}\rho(t)\bar{r}^3$$
we solve for $\rho$ via the coordinate transformation
$\bar{r}(t)=r(t)R(t)$ to find
\begin{equation}\label{eq:rhorel}
\rho=\frac{3}{4\pi}\frac{M}{\bar{r}(t)^3}=\frac{3\gamma}{\bar{r}(t)^2}=3\bar{\rho}.
\end{equation}
Now, we are able to compute the pressure $p$ of the FRW metric.
Indeed, substituting $k=0$, $\bar{p}=\bar{\sigma}\bar{\rho}$ and
$\rho=3\bar{\rho}$ into the jump condition~\eqref{eq:strenj}
yields
\begin{equation}\label{eq:psrho}
p=\sigma\rho,
\end{equation}
where
$$\sigma=\frac{1}{3}\left(\frac{\bar{\sigma}+1-A}{A} \right)
=\frac{\bar{\sigma}(\bar{\sigma}+5)}{3(\bar{\sigma}+1)}.$$
Analogous to~\cite{st94}  we can state this relation as
\begin{equation}\label{eq:H}
\bar{\sigma}=\frac{1}{2}\sqrt{9\sigma^2
-18\sigma+25}+\frac{3}{2}\sigma-\frac{5}{2}\equiv H(\sigma).
\end{equation}
As in non-lightlike case, see~\cite{st95}, within the region
$0\leq \sigma, \bar{\sigma}\leq 1$, we have $H(0)=0$,
$H'(\sigma)>0$, and $\bar{\sigma}<\sigma$; see
figure~\ref{fig:sigvssigbar}.
%---Graph of H----------------------------------------------------
%---End Figure----------------------------------------------------
However, in our case $H(1)=1$, which is not true in~\cite{st95}.
As we shall see below, having to satisfy the lightlike
condition~\eqref{eq:ndotn} will restrict $\sigma$ to a single
value. So far, we have that when $H(\sigma)=\bar{\sigma}$, the
jump condition~\eqref{eq:strenj} holds, and so we can apply
Theorem~\ref{thm:sphere} which means that we also get the
equivalencies of Theorem~\ref{thm:main}. This means we can say
that the matched FRW and TOV metrics form a lightlike solution of
the Einstein equations if the shock speed in $(t,r)$-coordinates
satisfy the lightlike condition~\eqref{eq:ndotn}. Therefore, we
must compute the formulas for $\rho(t)$, $R(t)$, the shock
positions $r(t)$ and $\bar{r(t)},$ the shock speeds
$\dot{\bar{r}}(t)$ and $\dot{r(t)},$ and then find the appropriate
values of $\sigma$ so that equation~\eqref{eq:ndotn} holds.

\subsection{The Shock Solution} We begin by substituting $p=\sigma \rho$ into the FRW
solutions~\eqref{eq:frwdt}, and~\eqref{eq:dRR} to get
\begin{equation}\label{eq:frwdt2}
dt=\mp\frac{1}{\sqrt{24\pi\mathcal{G}}(1+\sigma)}\rho^{-3/2}d\rho,
\end{equation}
and
\begin{equation}\label{eq:dRR2}
\frac{dR}{R}=-\frac{1}{3(1+\sigma)}\frac{d\rho}{\rho}.
\end{equation}
From~\eqref{eq:rhorel} we can write $\bar{r}$ in terms of $\rho$,
\begin{equation}\label{eq:rhorel2}
 \bar{r}=\sqrt{3\gamma}\rho^{-1/2},
\end{equation}
and then differentiating to get
\begin{equation}\label{eq:drhorel}
\rho^{-3/2}d\rho=-\frac{2}{\sqrt{3\gamma}}d\bar{r}.
\end{equation}
Notice that we can directly substitute this expression for
$\rho^{-3/2}d\rho$ into equation~\eqref{eq:frwdt2} which yields
\begin{equation}\label{eq:frwdt3}
dt=\mp\frac{1}{(1+\sigma)}\frac{1}{\sqrt{18\pi
\mathcal{G}\gamma}}d\bar{r}.
\end{equation}
Now integrating gives the shock position
\begin{equation}\label{eq:shockpos}
\bar{r}(t)=\pm\sqrt{18\pi\mathcal{G}\gamma}(1+\sigma)(t-t_0)+\bar{r}_0,
\end{equation}
and in conjunction with \eqref{eq:rhorel} we can compute the FRW
energy density as a function of $t$,
\begin{equation}\label{eq:rhosoln1}
\rho(t) =\frac{3\gamma}{\bar{r}(t)^2}
=\frac{3\gamma}{[\pm\sqrt{18\pi\mathcal{G}\gamma}(1+\sigma)(t-t_0)
+\bar{r}_0]^2}.
\end{equation}
Now, we can solve the differential equation~\eqref{eq:dRR2} for
$R(t)$ to find
\begin{equation}\label{eq:Rsoln}
R(t)=R_0\left(\frac{\rho}{\rho_0}\right)^{-1/3(1+\sigma)}
    =R_0\left(\frac{\bar{r}(t)}{\bar{r}_0}\right)^{2/3(1+\sigma)},
\end{equation}
from which we obtain
\begin{equation}\label{eq:shockpos2}
\begin{aligned}
r(t)=\bar{r}(t)R(t)^{-1} & =\bar{r}(t)R_0{-1}\left(\frac{\bar{r}(t)}{\bar{r}_0}\right)^{-2/3(1+\sigma)}\\
             & =\bar{r}_0R_0{-1}\left(\frac{\bar{r}(t)}{\bar{r}_0}\right)^{(1+3\sigma)/(3+3\sigma)}.
\end{aligned}
\end{equation}

\subsection{The Shock Speeds} To compute the shock speeds in
$(t,\bar{r})$- and $(t,r)$-coordinates we differentiate
equations~\eqref{eq:shockpos} and~\eqref{eq:shockpos2} with
respect to $t$, and get
\begin{equation}\label{eq:shckspdbar}
\dot{\bar{r}} =3(1+\sigma)
\sqrt{\frac{\bar{\sigma}}{1+6\bar{\sigma}+\bar{\sigma}^2}},
\end{equation}
and
\begin{equation}\label{eq:shckspd}
\dot{r} =\frac{1+3\sigma}{R(t)}
\sqrt{\frac{\bar{\sigma}}{1+6\bar{\sigma}+\bar{\sigma}^2}}.
\end{equation}

\subsection{The Lightlike Condition}
In~\cite{st95}, Smoller and Temple show that the shock speed,
relative to the FRW particles, is given by
\begin{equation}\label{eq:shckspdform}
s(\sigma)=(1+3\sigma) \sqrt{\frac{\bar{\sigma}}{1
+6\bar{\sigma}+\bar{\sigma}^2}},
\end{equation}
where we use $\bar{\sigma}=H(\sigma)$, which is defined in
equation~\eqref{eq:H}, to make this expression a function of only
$\sigma$. The function $s(\sigma)$ is valid in our case as well,
where the only deviation from Smoller and Temple's argument,
given in~\cite[lemma 1]{st95}, is that we must use our version of
$H(\sigma)$ which is slightly different than the one they
derived.

\begin{figure}
  \centering
  \includegraphics[width=.7\textwidth]{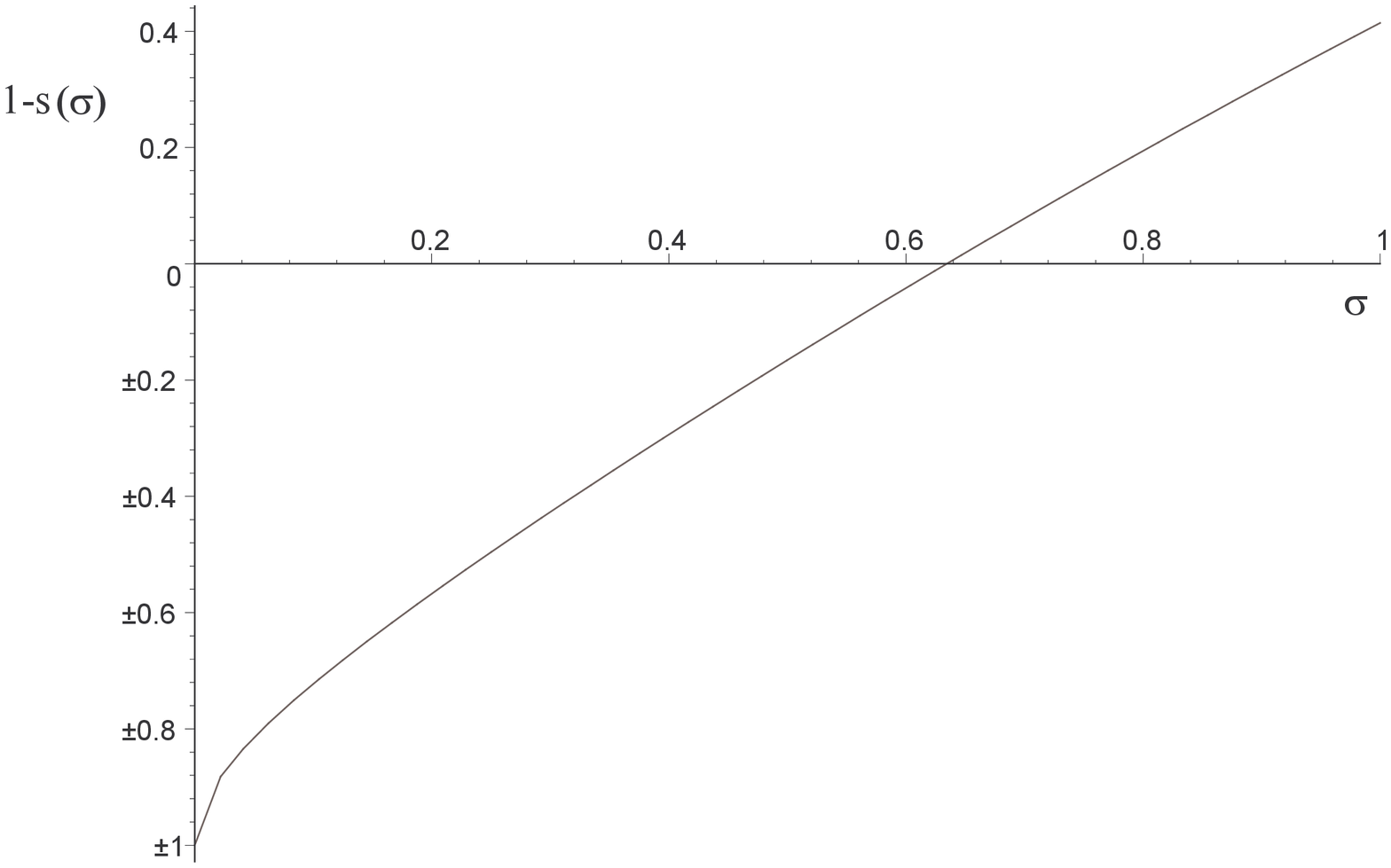}
  \caption{A graph of $1-s(\sigma)$
  defined in equation~\eqref{eq:sigma2}. The solution of
  $1-s(\sigma)=0$ gives a value of
  $\sigma\equiv\sigma_2\approx 0.63442$.}
  \label{fig:sigma2}
\end{figure}

For our lightlike solution to be valid we must satisfy the extra
condition which is the lightlike condition given
in~\eqref{eq:ndotn} where $k=0$.  Therefore, we must have that
\begin{equation}\label{eq:llcond}
s(\sigma)=1.
\end{equation}
Solving $s(\sigma)=1$ for $\sigma$, see figure~\ref{fig:sigma2},
gives
\begin{equation}\label{eq:sigma2}
\sigma\equiv\sigma_2\approx 0.63442.
\end{equation}
%---Graph of L=1--------------------------------------------------
%---End Figure----------------------------------------------------

We can now characterize our lightlike shock-wave solution which
we do in the following theorem.
\begin{theorem}\label{thm:shcksln}
Assume an equation of state of the form
$\bar{p}=\bar{p}(\bar{\rho})$ for the TOV metric, and $p=\sigma
\rho$ for the FRW metric. Also, assume that relation between
$\sigma$ and $\bar{\sigma}$, given by
$$\bar{\sigma}
=\frac{1}{2}\sqrt{9\sigma^2-18\sigma+25}
+\frac{3}{2}\sigma-\frac{5}{2}=H(\sigma)$$ holds, where we have
taken $k=0$. Then the TOV solution given by
$$\bar{\rho}=\gamma/\bar{r}^2,$$
$$M(\bar{r})=4\pi\gamma\bar{r},$$
$$A=1-8\pi\mathcal{G}\gamma,$$
and
$$B=B_0\left(\frac{\bar{r}}{\bar{r}_0}\right)^{4\bar{\sigma}_2/(1
+\bar{\sigma}_2)}$$ will match the FRW solution given by
$$\rho(t)=\frac{3\gamma}{[\pm\sqrt{18\pi\mathcal{G}\gamma}(1+
\sigma_2)(t-t_0)+\bar{r}_0]^2},$$ and
$$R(t)=R_0\left(\frac{\bar{r}(t)}{\bar{r}_0}\right)^{2/3(1+\sigma_2)}$$
across the shock surface
$$\bar{r}(t)=\pm\sqrt{18\pi\mathcal{G}\gamma}(1+\sigma_2)(t-t_0)+\bar{r}_0,$$
such that conservation of energy and momentum hold across the
surface. We have used $\sigma_2$, and
$H(\sigma_2)=\bar{\sigma}_2$, where $\sigma_2$ is the solution of
$s(\sigma)=1$, and $s(\sigma)$ denotes the shock speed given
in~\eqref{eq:shckspdform}.
\end{theorem}

\subsection{The Lax Shock Conditions}
The final question we would like to answer is how do we classify
the shock solution given in Theorem~\ref{thm:shcksln} with respect
to the shock conditions first given by Smoller and Temple
in~\cite{st95}? As in~\cite{st95} we only consider the case when
the pressure and density are greater behind the shock wave. Since
$\rho=3\bar{\rho}$ this means that the FRW region is behind the
shock, the TOV region is in front of the shock. Therefore, for our
solution we take the plus sign for expression of $\dot{R}$
in~\eqref{eq:ode2ss}, and the corresponding signs in
equations~\eqref{eq:ode1ss},~\eqref{eq:frwdt},
and~\eqref{eq:frw_t}.

The shock speed given in~\eqref{eq:shckspdform} is coordinate
dependent, and is computed in a \emph{locally Minkowskian}
frame~\cite{stg01,st95,stgerm99}. A coordinate system is called
locally Minkowskian at a point $p$ if
$g_{\alpha\beta}(p)=\mbox{diag}(-1,1,1,1)$, but not necessarily
Lorentzian where it is also required that
$g_{\alpha\beta,\gamma}(p)=0.$ Also, note that since we are
working with the radial component exclusively we only need to
consider a locally Minkowskian frame in the $(t,r)$-coordinates.
We denote our locally Minkowskian coordinates by $(t,\tilde{r}),$
and they can be obtained by letting $r=\varphi(\tilde{r})$ and
choosing $\varphi$ so that $\varphi'=1/R^2,$ which implies
$$ds^2=-dt^2+R^2dr^2 \To d\tilde{s}^2=-dt^2+d\tilde{r}^2.$$
The shock speed $s(\sigma)$ is determined in a locally
Minkowskian frame comoving with the FRW metric.

Recall that we have chosen $\sigma=\sigma_2\approx 0.63442$ to
satisfy $\dot{r}=1/R$, and since
$$\dot{r}=\frac{dr}{dt}=\frac{1}{R}\frac{d\tilde{r}}{dt},$$ it follows
that
$$s(\sigma_2)=1,$$
in $(t,\tilde{r})$-coordinates. The characteristic speeds behind
the shock are equal to the sound speed $\pm\sqrt{\sigma}$ in the
$(t,\tilde{r})$-coordinate system, since the FRW fluid is
comoving with respect to $(t,\tilde{r})$-coordinates~\cite{st95}.
That is, the speeds of the characteristics relative to the FRW
fluid are given by
$$\tilde{\lambda}^{\pm}_{FRW}\equiv\pm\frac{d\tilde{r}}{dt}
=\pm\sqrt{\sigma}.$$ In our case, with the shock moving outward
with respect to $r$ and $\bar{r}$, the Lax characteristic
condition, which says that the characteristic curves of the shock
family impinge on the shock from both sides and all other
characteristics cross the shock, hold if and only if
\begin{equation}\label{eq:lax1}
 \tilde{\lambda}^{+}_{TOV}<s<\tilde{\lambda}^{+}_{FRW},
\end{equation}
where $\tilde{\lambda}^{+}_{TOV}$ denotes the corresponding
characteristic speed on the TOV side, or the front side, of the
shock. Then in our lightlike case we have,
$$\tilde{\lambda}^{+}_{FRW}
=\sqrt{\sigma_2}\approx \sqrt{0.63442}=.79650<1=s(\sigma_2),$$
which means that our shock solution does not satisfy the Lax
characteristic condition. We can also conclude that
\begin{equation}\label{eq:shckcon1}
 \tilde{\lambda}^{-}_{FRW} <\tilde{\lambda}^{+}_{FRW} <s(\sigma_2).
\end{equation}

\begin{figure}
  \centering
  \includegraphics[width=.7\textwidth]{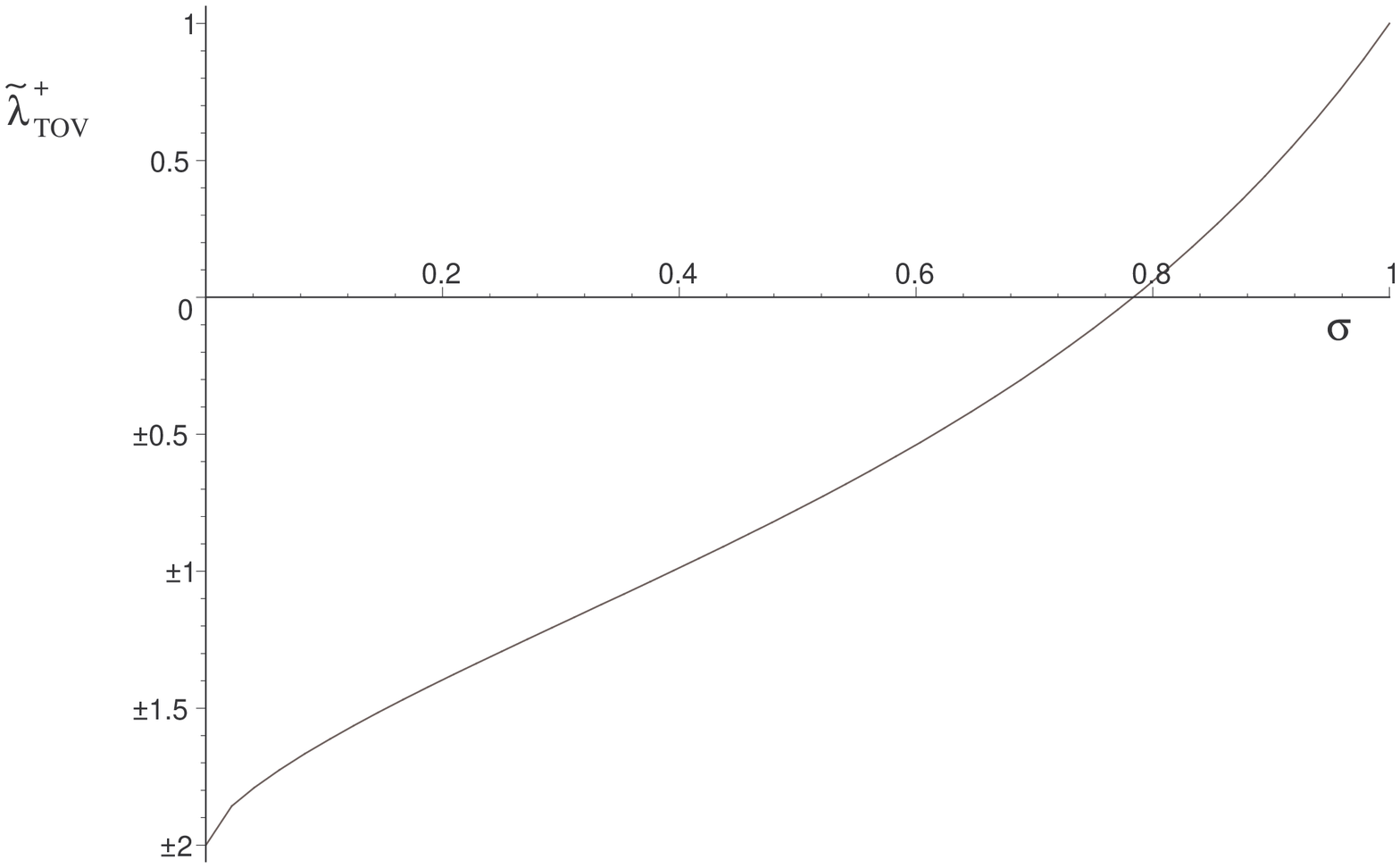}
  \caption{This graph of $\tilde{\lambda}^{+}_{TOV}(\sigma)$,
    defined in equation~\eqref{eq:TOVchar}, shows that
    $\tilde{\lambda}^{+}_{TOV}(\sigma_2)<0<s(\sigma_2)=1,$ where
    $\sigma_2\approx 0.63442.$}
  \label{fig:tovchar}
\end{figure}

Now we would like to know how $s(\sigma_2)$ relates to the
characteristics $\tilde{\lambda}^{\pm}_{TOV}$. In~\cite{st95}
Smoller and Temple show that
\begin{equation}\label{eq:TOVchar}
\tilde{\lambda}^{+}_{TOV}(\sigma) \equiv
-\frac{2-\sqrt{\bar{\sigma}^2+6\bar{\sigma}+1}}
{\sqrt{\bar{\sigma}^2+6\bar{\sigma} +1}-2\bar{\sigma}}
\sqrt{\bar{\sigma}}.
\end{equation}
Then for $\sigma=\sigma_2\approx 0.63442$ we have
$$\tilde{\lambda}^{+}_{TOV}(\sigma_2)\approx -0.45040 < 1
=s(\sigma_2),$$ see figure~\ref{fig:tovchar}.
%---Graph of TOV Characteristic-----------------------------------
%---End Figure----------------------------------------------------
Therefore, we have that
$$ \tilde{\lambda}^{-}_{FRW} <\tilde{\lambda}^{+}_{FRW} <s(\sigma_2),$$
and
$$ \tilde{\lambda}^{-}_{TOV} <\tilde{\lambda}^{+}_{TOV} <s(\sigma_2).$$
Thus both sets of characteristics cross the shock since its speed
is greater than the characteristic speeds on each side of the
shock.

\section{Comparison with the Smoller-Temple Sublight Shocks}
In~\cite{st95} Smoller and Temple show that in the limit as their
subluminal shock solutions tend to the speed of light,
$$\sigma \To \sigma_2=0.745$$
However, we have shown that an actual shock solution moving at the
speed of light the value $\sigma_2=0.63442$. This yields the
unexpected result that the solution is not equal to the limit of
Smoller-Temple subluminous solutions as they tend to the speed of
light.

% ---End Chapter 4------------------------------------------------

%% file: conc.tex
% ---Conclusion---------------------------------------------------
\label{chp:5}

We have generalized the work of Smoller and Temple given
in~\cite{st94} and~\cite{st95} in the sense that we can include
lightlike shock surfaces. In defining a more general second
fundamental form, by replacing the normal vector with a transverse
vector satisfying the jump conditions~\eqref{eq:njump}
and~\eqref{eq:njump1}, we were able to overcome the breakdown of
the standard second fundamental form for lightlike hypersurfaces.
Then we used this generalized second fundamental form, in an
analogous way to how Smoller and Temple applied the standard
second fundamental form in~\cite{st94}, to obtain a theory for the
lightlike case. In the process, we introduced a modified Gaussian
Skew coordinate system. This theory yielded the unexpected result
of having to include the condition that the metric be $C^2$ on the
spacelike subspace of $\tans$.

Then we were able to construct the an exact, spherically symmetric
shock-wave solution of the Einstein equations which propagates at
the speed of light. No quantities are moving at the speed of light
except the shock. Although, our solution was consistent with
Smoller and Temple's results in~\cite{st95}, there was the
unexpected difference in the value of $\sigma_2$. In the
non-lightlike computation given~\cite{st95} the value for $\sigma$
in the lightlike limit was
  $$\sigma_2 \approx 0.745,$$
while our lightlike computations yielded
  $$\sigma_2 \approx
  0.63442.$$
From this we conclude that the limit of the Smoller-Temple
subluminal solution as it tends to the speed of light is not equal
to our solution propagating at the speed of light.

Furthermore, we showed that in this exact solution the pressure
and density are finite on each side of the shock throughout the
solution, the sound speeds, on each side of the shock, are
constant and subluminous. Moreover, the pressure and density are
smaller at the leading edge of the shock which is consistent with
the entropy conditions in classical gas dynamics~\cite{lax, smol}.

\section{Summary of Contributions}
In the lightlike shock matching theory of Chapter~\ref{chp:3} we
applied the notion of a generalized second fundamental form given
by Barrab\`{e}s and Israel in~\cite{isr91} into the shock matching
framework of Smoller and Temple~\cite{st94,st95}. The generalized
second fundamental form $\K$ in~\cite{isr91} was given in a scalar
form, and we modified it so that $\K$ is mapping that takes
tangent vectors on the surface to tangent vectors on surface. Then
we based the analysis on a Gaussian skew type coordinate system
which we had to modify so that we could incorporate the
generalized second fundamental form into the shock matching
theory. Most of the supporting lemmas, and each of the theorems
in~\cite{st94} had to be modified and proved again in the context
of this new Modified Gaussian Skew (MGS) coordinate system and the
generalized second fundamental form. This yielded the unexpected
result of adding the extra condition in our main result,
Theorem~\ref{thm:main}, that the metric inner product on the
spacelike subspace of $\tans$ to be $C^2$ in order for the
Rankine-Hugoniot jump conditions to hold. This condition was
already satisfied for the spherically case in
Theorem~\ref{thm:sphere}.

Then in Chapter~\ref{chp:4} we constructed a new exact shock-wave
solution moving at the speed of light. In order to do this we had
to construct appropriate transverse vectors $N$, and apply our
extension of the subluminal theory to the specific example of the
matched FRW/TOV metrics. This included significant modification of
the equations expressing the Rankine-Hugoniot jump conditions
given in equation~\eqref{eq:strenj}. This yielded the difference
in the limit of the Smoller-Temple subluminal solutions as they
tended to the speed of light with the actual solution propagating
at the speed of light. Lastly, we incorporated the modifications
of the subluminal shock solutions to show that our new exact
solution was a crossing shock.

%% file: thesis_with_title.bbl
\providecommand{\bysame}{\leavevmode\hbox to3em{\hrulefill}\thinspace}
\providecommand{\MR}{\relax\ifhmode\unskip\space\fi MR }
\providecommand{\MRhref}[2]{%
  \href{http://www.ams.org/mathscinet-getitem?mr=#1}{#2}
}
\providecommand{\href}[2]{#2}
\begin{thebibliography}{10}

\bibitem{isr91}
C.~Barrab\`{e}s and W.~Israel, \emph{{Thin shells in general relativity and
  cosmology: The Lightlike limit}}, Physical Review D \textbf{43} (1991),
  no.~4, 1129--1142.

\bibitem{einst22}
Albert Einstein, \emph{{The Meaning of Relativity}}, fifth ed., Princeton
  University Press, Princeton, New Jersey, 1922.

\bibitem{evans}
Lawrence~C. Evans, \emph{{Partial Differential Equations}}, Graduate Studies in
  Mathematics, vol.~19, American Mathematical Society, Providence, Rhode
  Island, 1998.

\bibitem{stg01}
Jeff Groah, Joel Smoller, and Blake Temple, \emph{{Solving the Einstein
  Equations by Lipschitz Continuous Metrics: Shock Waves in General
  Relativity}}, Handbuch der Physik, Germany, 2001, (to appear).

\bibitem{hawkellis}
S.W. Hawking and G.F.R. Ellis, \emph{{The Large Scale Structure of
  Space-Time}}, Cambridge University Press, New York, 1973.

\bibitem{isr66}
W.~Israel, \emph{{Singular Hypersurfaces and Thin Shells in General
  Relativity}}, Nuovo Cimento B (11) \textbf{44} (1966), no.~1, 1--14.

\bibitem{lax}
Peter~D. Lax, \emph{{Hyperbolic systems of conservation laws and the
  mathematical theory of shock waves}}, Conference board of the mathematical
  sciences, Regional conference series in applied mathematics, vol.~11, Society
  for Industrial and Applied Mathematics, 1973.

\bibitem{lee}
John~M. Lee, \emph{{Riemannian Manifolds: An Introduction to Curvature}},
  Graduate Texts in Mathematics, vol. 176, Springer-Verlag, New York, 1997.

\bibitem{mtw}
Charles~W. Misner, Kip~S. Thorne, and John~A. Wheeler, \emph{{Gravitation}}, W.
  H. Freeman and Company, New York, 1973.

\bibitem{oneill}
Barrett O'Neill, \emph{{Semi-Riemannian Geometry with Applications to
  Relativity}}, Pure and Applied Mathematics, vol. 103, Academic Press, San
  Diego, 1983.

\bibitem{os39}
J.~Robert Oppenheimer and Hartland Snyder, \emph{{On Continued Gravitational
  Contraction}}, Physical Review \textbf{56} (1939), 455--459.

\bibitem{schutz}
Bernard~F. Schutz, \emph{{A First Course in General Relativity}}, Cambridge
  University Press, Cambridge, 1985.

\bibitem{smol}
Joel Smoller, \emph{{Shock Waves and Reaction-Diffusion Equations}}, second
  ed., Springer-Verlag, New York, 1983.

\bibitem{st94}
Joel Smoller and Blake Temple, \emph{{Shock-wave solutions of the Einstein
  equations: The Oppenheimer-Snyder model of gravitational collapse extended to
  the case of non-zero pressure}}, Arch. Rational Mech. Anal. \textbf{128}
  (1994), 249--297.

\bibitem{st95}
\bysame, \emph{{Astophysical shock-wave solutions of the Einstein Equations}},
  Physical Review D \textbf{51} (1995), 2733--2743.

\bibitem{stgerm99}
\bysame, \emph{{Shock-Wave Solutions of the Einstein Equations: A General
  Theory with Examples}}, May 16-22 1999, {Proceedings of the European Union
  Research Network's Third Annual Summer School Training and Mobility of
  Researchers Program, Lambrecht (Pfalz) Germany}, (to appear).

\bibitem{temp99}
Blake Temple, \emph{{Shock-Waves and Geometry}}, March 1998, Grant Proposal-NSF
  Applied Math, UC Davis Math 1998.

\bibitem{wald}
Robert~M. Wald, \emph{{General Relativity}}, University of Chicago Press,
  Chicago, 1984.

\bibitem{warner}
Frank~W. Warner, \emph{{Foundations of Differentiable Manifolds and Lie
  Groups}}, Graduate Texts in Mathematics, vol.~94, Springer-Verlag, New York,
  1983.

\bibitem{weinberg}
Steven Weinberg, \emph{{Gravitation and Cosmology: Principles and Applications
  of the General Theory of Relativity}}, John Wiley \& Sons, Inc., New York,
  1972.

\end{thebibliography}
